%% file: ms.tex
\documentclass[a4paper,12pt]{article}
\usepackage{amsmath,amssymb,amsthm,mathtools}
\usepackage{graphicx}
\usepackage{enumitem}
\usepackage{color}
\usepackage{epsfig}
\usepackage{subcaption}
\usepackage[font=small]{caption}
\usepackage[hang,flushmargin]{footmisc} 
\usepackage{float}
\usepackage{booktabs}
\usepackage[mathscr]{euscript}
\usepackage{bbm}
\usepackage{natbib}
\usepackage{setspace}
\usepackage{mathrsfs}
\usepackage[left=2.8cm,right=2.8cm,bottom=2.8cm,top=2.8cm]{geometry}

\input{macros}

\begin{document}

\heading{Multiscale Clustering}{of Nonparametric Regression Curves}
\vspace{-0.6cm}

\authors{Michael Vogt\renewcommand{\thefootnote}{1}\footnotemark[1]}{University of Bonn}{Oliver Linton\renewcommand{\thefootnote}{2}\footnotemark[2]}{University of Cambridge}
\footnotetext[1]{Address: Department of Economics and Hausdorff Center for Mathematics, University of Bonn, 53113 Bonn, Germany. Email: \texttt{michael.vogt@uni-bonn.de}.}
\footnotetext[2]{Address: Faculty of Economics, Austin Robinson Building, Sidgwick Avenue, Cambridge, CB3 9DD, UK. Email: \texttt{obl20@cam.ac.uk}.}
\vspace{-0.8cm}

\renewcommand{\thefootnote}{\arabic{footnote}}
\setcounter{footnote}{2}

\renewcommand{\abstractname}{}
\begin{abstract}
\noindent 
In a wide range of modern applications, we observe a large number of time series rather than only a single one. It is often natural to suppose that there is some group structure in the observed time series. When each time series is modelled by a nonparametric regression equation, one may in particular assume that the observed time series can be partitioned into a small number of groups whose members share the same nonparametric regression function. We develop a bandwidth-free clustering method to estimate the unknown group structure from the data. More precisely speaking, we construct multiscale estimators of the unknown groups and their unknown number which are free of classical bandwidth or smoothing parameters. In the theoretical part of the paper, we analyze the statistical properties of our estimators. Our theoretical results are derived under general conditions which allow the data to be dependent both in time series direction and across different time series. The technical analysis of the paper is complemented by a simulation study and a real-data application.
\end{abstract}

\renewcommand{\baselinestretch}{1.2}\normalsize
\renewcommand{\theequation}{\thesection.\arabic{equation}}




\setlength{\parindent}{0.75cm} 
\allowdisplaybreaks[2]

\section{Introduction}

In this paper, we are concerned with the problem of clustering nonparametric regression curves. We consider the following model setup: We observe a large number of time series $\mathcal{T}_i = \{ (Y_{it},X_{it}): 1 \le t \le T \}$ for $1 \le i \le n$. For simplicity, we synonymously speak of the $i$-th time series, the time series $i$ and the time series $\mathcal{T}_i$ in what follows. Each time series $\mathcal{T}_i$ satisfies the nonparametric regression equation
\begin{equation}\label{model-reg-intro}
Y_{it} = m_i(X_{it}) + u_{it}
\end{equation}
for $t=1,\ldots,T$, where $m_i$ is an unknown smooth function which is evaluated at the design points $X_{it}$ and $u_{it}$ denotes the error term. The $n$ time series in our sample are supposed to belong to $K_0$ different groups. More specifically, the set of time series $\{1,\ldots,n\}$ can be partitioned into $K_0$ groups $G_1,\ldots,G_{K_0}$ such that for each $k=1,\ldots,K_0$, 
\begin{equation}\label{model-group-intro}
m_i = m_j \quad \text{ for all } i,j \in G_k. 
\end{equation} 
Hence, the members of each group $G_k$ all have the same regression function. A detailed description of model \eqref{model-reg-intro}--\eqref{model-group-intro} can be found in Section \ref{sec-model}. Our modelling approach provides a parsimonious way to deal with a potentially very large number of time series $n$. It thus stands in the tradition of multiple time series analysis, an area which greatly benefited from the pioneering work of George Tiao.

An interesting statistical problem is how to construct estimators of the unknown groups $G_1,\ldots,G_{K_0}$ and their unknown number $K_0$ in model \eqref{model-reg-intro}--\eqref{model-group-intro}. For the special case that the design points $X_{it} = t/T$ represent (rescaled) time and the functions $m_i$ are nonparametric time trends, this problem has been analyzed for example in \cite{LuanLi2003} and \cite{Degras2012}. For the case that $X_{it}$ are general random design points which may differ across time series $i$, \cite{VogtLinton2017} have developed a thresholding method to estimate the unknown groups and their number. Notably, their approach can also be adapted to the case of deterministic regressors $X_{it}$, in particular to the case that $X_{it} = t/T$. The model \eqref{model-reg-intro}--\eqref{model-group-intro} with the fixed design points $X_{it} = t/T$ is closely related to models from functional data analysis. There, the aim is to cluster smooth random curves that are functions of (rescaled) time and that are observed with or without noise. A number of different clustering approaches have been proposed in the context of functional data models; see for example \cite{Abraham2003}, \cite{Tarpey2003} and \cite{Tarpey2007} for procedures based on $k$-means clustering, \cite{James2003} and \cite{Chiou2007} for model-based clustering approaches and \cite{Jacques2014} for a recent survey.

Virtually all of the proposed procedures to cluster nonparametric curves in model \eqref{model-reg-intro}--\eqref{model-group-intro} and in related functional data settings have the following drawback: they depend on a number of smoothing parameters required to estimate the nonparametric functions $m_i$. A common approach is to approximate the functions $m_i$ by a series expansion $m_i(x) \approx \sum_{j=1}^{L} \beta_{ij} \phi_j(x)$, where $\{ \phi_j: j =1,2,\ldots \}$ is a function basis and $L$ is the number of basis elements taken into account for the estimation of $m_i$. Here, $L$ plays the role of the smoothing parameter and may vary across $i$, that is, $L = L_i$. To estimate the classes $G_1,\ldots,G_{K_0}$, estimators $\widehat{\boldsymbol{\beta}}_i$ of the coefficient vectors $\boldsymbol{\beta}_i = (\beta_{i1},\ldots,\beta_{iL})^\top$ are clustered into groups by a standard clustering algorithm. Variants of this approach have for example been investigated in \cite{Abraham2003}, \cite{LuanLi2003}, \cite{Chiou2007} and \cite{Tarpey2007}. Another approach is to compute nonparametric estimators $\widehat{m}_{i} = \widehat{m}_{i,h}$ of the functions $m_i$ for some smoothing parameter $h$ (which may differ across $i$) and to calculate distances $\widehat{\rho}_{ij} = \rho(\widehat{m}_{i},\widehat{m}_{j})$ between the estimates $\widehat{m}_{i}$ and $\widehat{m}_{j}$, where $\rho(\cdot,\cdot)$ is a distance measure such as a supremum or an $L_2$-distance. A distance-based clustering algorithm is then applied to the distances $\widehat{\rho}_{ij}$. This strategy has for example been used in \cite{VogtLinton2017}.

In general, nonparametric curve estimators strongly depend on the chosen smoothing or bandwidth parameters. A clustering procedure which is based on such estimators can be expected to be strongly influenced by the choice of smoothing parameters as well. To see this issue more clearly, consider two time series $i$ and $j$ from two different groups. The corresponding regression functions $m_i$ and $m_j$ may differ on different scales. In particular, they may differ on a local/global scale, that is, they may have certain local/global features which distinguish them from each other. For example, they may be identical except for a sharp local spike, or they may have a slightly different curvature globally all over their support. Whether nonparametric estimators are able to pick up local/global features of $m_i$ and $m_j$ depends on the chosen bandwidth. When the bandwidth is large, the estimators capture global features of $m_i$ and $m_j$ but smooth out local ones. When the bandwidth is small, they pick up local features, whereas more global ones are poorly captured. As a consequence, a clustering algorithm which is based on nonparametric estimators of $m_i$ and $m_j$ will reliably detect local/global differences between the functions $m_i$ and $m_j$ only if the bandwidths are chosen appropriately. The clustering results produced by such an algorithm can thus be expected to vary considerably with the chosen bandwidths.

The main aim of this paper is to construct estimators of the unknown groups $G_1,\ldots,G_{K_0}$ and of their unknown number $K_0$ in model \eqref{model-reg-intro}--\eqref{model-group-intro} which are free of classical smoothing or bandwidth parameters. To achieve this, we construct a clustering algorithm which is based on statistical multiscale methods. In recent years, a number of multiscale techniques have been developed in the context of statistical hypothesis testing. Early examples are the SiZer approach of \cite{ChaudhuriMarron1999, ChaudhuriMarron2000} and the multiscale tests of \cite{HorowitzSpokoiny2001} and \cite{DuembgenSpokoiny2001}. More recent references include the tests in \cite{SchmidtHieber2013}, \cite{ArmstrongChan2016}, \cite{EckleBissantzDette2017} and \cite{ProkschWernerMunk2018} among others. In this paper, we develop multiscale techniques for clustering rather than testing purposes. Roughly speaking, we proceed as follows: To start with, we construct statistics which measure the distances between pairs of functions $m_i$ and $m_j$. To do so, we estimate the functions $m_i$ and $m_j$ at different resolution levels, that is, with the help of different bandwidths $h$. The resulting estimators are aggregated in supremum-type statistics which simultaneously take into account multiple bandwidth levels. We thereby obtain multiscale statistics which avoid the need to pick a specific bandwidth. To estimate the unknown classes $G_1,\ldots,G_{K_0}$, we combine the constructed multiscale statistics with a hierarchical clustering algorithm. To estimate the unknown number of classes $K_0$, we develop a thresholding rule that is applied to the dendrogram produced by the clustering algorithm. Alternatively, the multiscale statistics may be combined with other distance-based clustering algorithms. In particular, they can be used to turn the estimation strategy of \cite{VogtLinton2017} into a bandwidth-free procedure. We comment on this in more detail in Section \ref{sec-ext} of the paper.

By construction, our multiscale clustering methods allow to detect differences between the functions $m_i$ at different scales or resolution levels. An alternative way to achieve this is to employ Wavelet methods. A Bayesian Wavelet-based method to cluster nonparametric curves has been developed in \cite{Ray2006}. There, the model $Y_{it} = m_i(t/T) + u_{it}$ is considered, where $m_i$ are smooth functions of rescaled time $t/T$ and the error terms $u_{it}$ are restricted to be i.i.d.\ Gaussian noise. To the best of our knowledge, there are no Wavelet-based clustering methods available in the literature which allow to deal with the model setting \eqref{model-reg-intro}--\eqref{model-group-intro} under general conditions on the design points $X_{it}$ and the error terms $u_{it}$. Our methods and theory, in contrast, allow to do so. In particular, we do not restrict attention to the special case that $X_{it} = t/T$ but allow for general design points $X_{it}$ that may differ across $i$. Moreover, we do not restrict the error terms to be Gaussian but only impose some moderate moment conditions on them. In addition, we allow them to be dependent both across $t$ and $i$.

The problem of estimating the unknown groups and their unknown number in model \eqref{model-reg-intro}--\eqref{model-group-intro} is closely related to a developing literature in econometrics which aims to identify the unknown group structure in parametric panel regression models. The clustering problem considered in this literature can be regarded as a parametric version of our problem. In its simplest form, the panel regression model under consideration is given by the equation $Y_{it} = \boldsymbol{\beta}_i^\top X_{it} + u_{it}$ for $1 \le t \le T$ and $1 \le i \le n$, where the coefficient vectors $\boldsymbol{\beta}_i$ are allowed to vary across individuals $i$. Similarly as in our nonparametric model, the coefficients $\boldsymbol{\beta}_i$ are assumed to belong to a number of groups: there are $K_0$ groups $G_1,\ldots,G_{K_0}$ such that $\boldsymbol{\beta}_i = \boldsymbol{\beta}_j$ for all $i,j \in G_k$ and all $1 \le k \le K_0$. The problem of estimating the unknown groups and their unknown unknown number has been studied in different versions of this modelling framework in \cite{BonhommeManresa2015}, \cite{SuShiPhillips2016}, \cite{WangPhillipsSu2018} and \cite{SuJu2018} among others. Note that our clustering methods can be adapted in a straightforward way to a number of semiparametric models which are middle ground between the fully parametric panel models just discussed and our nonparametric framework. In Section \ref{sec-ext}, we discuss in more detail how to achieve this.

Our estimation methods are described in detail in Sections \ref{sec-dist-measure}--\ref{sec-est-K0}. In Section \ref{sec-dist-measure}, we construct the multiscale statistics that form the basis of our clustering methods. Section \ref{sec-est-classes} introduces the hierarchical clustering algorithm to estimate the unknown classes $G_1,\ldots,G_{K_0}$. In Section \ref{sec-est-K0}, we finally describe the procedure to estimate the unknown number of classes $K_0$. The main theoretical result of the paper is laid out in Section \ref{sec-theory}. This result characterizes the asymptotic convergence behaviour of the multiscale statistics and forms the basis to derive the theoretical properties of our clustering methods. To explore the finite sample properties of our approach and to illustrate its advantages over bandwidth-dependent clustering algorithms, we conduct a simulation study in Section \ref{sec-sim}. Moreover, we illustrate the procedure by an application from finance in Section \ref{sec-app}.

\setcounter{equation}{0}
\section{The model}\label{sec-model}

As already mentioned in the Introduction, we observe $n$ different time series $\mathcal{T}_i = \{ (Y_{it},X_{it}) : 1 \le t \le T\}$ of length $T$ for $1 \le i \le n$. In what follows, we describe in detail how the observed data $\{ \mathcal{T}_i: 1 \le i \le n \}$ are modelled. For our theoretical analysis, we regard the number of time series $n$ as a function of $T$, that is, $n = n(T)$. The time series length $T$ is assumed to tend to infinity, whereas the number of time series $n$ may be either bounded or diverging. The exact technical conditions on $T$ and $n$ are laid out in Section \ref{sec-theory}. Throughout the paper, asymptotic statements are to be understood in the sense that $T \rightarrow \infty$.

\subsection{The model for time series $\boldsymbol{\mathcal{T}_i}$}\label{subsec-model-reg}

Each time series $\mathcal{T}_i$ in our sample is modelled by the nonparametric regression equation
\begin{equation}\label{model-eq1}
Y_{it} = m_i(X_{it}) + u_{it}
\end{equation}
for $1 \le t \le T$, where $m_i$ is an unknown smooth function and $u_{it}$ denotes the error term. We focus attention on the case that the design points $X_{it}$ are random as this is the technically more involved case. Our methods can be adapted to deterministic design points $X_{it}$ with some minor modifications. To keep the exposition as simple as possible, we assume that the regressors $X_{it}$ are real-valued. As discussed in Section \ref{sec-ext}, our methods and theory carry over to the multivariate case in a straightforward way. We further suppose that the regressors $X_{it}$ have compact support, which w.l.o.g.\ is equal to $[0,1]$ for each $i$. The error terms $u_{it}$ in \eqref{model-eq1} are assumed to have the additive component structure
\begin{equation}\label{model-eq2}
u_{it} = \alpha_i + \gamma_t + \varepsilon_{it},
\end{equation}
where $\varepsilon_{it}$ are standard regression errors that satisfy $\ex[\varepsilon_{it}|X_{it}] = 0$ and the terms $\alpha_i$ and $\gamma_t$ are so-called fixed effects. The expression $\alpha_i$ is an error component which is specific to the $i$-th time series $\mathcal{T}_i$. It can be interpreted as capturing unobserved characteristics of the time series $\mathcal{T}_i$ which are stable over time. Suppose for instance that the observations of $\mathcal{T}_i$ are sampled from some subject $i$. In this case, $\alpha_i$ can be regarded as controlling for time-invariant unobserved characteristics of subject $i$, such as intelligence or certain unknown genetic factors. Similarly, the term $\gamma_t$ captures unobserved time-specific effects like calendar effects or trends that are common across time series $i$. In many applications, the regressors may be correlated with unobserved subject- or time-specific characteristics. To take this into account, we allow the errors $\alpha_i$ and $\gamma_t$ to be correlated with the regressors in an arbitrary way. Specifically, defining $\mathcal{X}_{n,T} = \{X_{it}: 1 \le i \le n, \, 1 \le t \le T \}$, we permit that $\ex[\alpha_i|\mathcal{X}_{n,T}] \ne 0$ and $\ex[\gamma_t|\mathcal{X}_{n,T}] \ne 0$. The error terms $\varepsilon_{it}$ are allowed to be dependent across $t$ but are assumed to be independent across $i$. The fixed effects $\alpha_i$, in contrast, may be correlated across $i$ in an arbitrary way. Hence, by including $\alpha_i$ and $\gamma_t$ in the error structure, we allow for some restricted types of cross-sectional dependence in the errors $u_{it}$. As a result, we accommodate for both time series dependence and certain forms of cross-sectional dependence in the error terms of our model. The exact conditions on the dependence structure are stated in \ref{C-mixing} in Section \ref{sec-theory}.

\subsection{The group structure}\label{subsec-model-groups}

We impose the following group structure on the time series $\mathcal{T}_i$ in our sample: There are $K_0$ groups of time series $G_1,\ldots,G_{K_0}$ with $\dot\bigcup_{k=1}^{K_0} G_k = \{1,\ldots,n\}$ such that for each $1 \le k \le K_0$, 
\begin{equation}\label{model-group-eq}
m_i = m_j \quad \text{ for all } i,j \in G_k. 
\end{equation} 
Put differently, for each $1 \le k \le K_0$,  
\begin{equation}\label{model-group-eq-alt}
m_i = g_k \quad \text{ for all } i \in G_k, 
\end{equation}
where $g_k$ is the group-specific regression function associated with the class $G_k$. According to \eqref{model-group-eq-alt}, the time series of a given class $G_k$ all have the same regression curve $g_k$. To make sure that time series which belong to different classes have different regression curves, we suppose that $g_k \ne g_{k^\prime}$ for $k \ne k^\prime$. The exact technical conditions on the functions $g_k$ are summarized in \ref{C-m} in Section \ref{sec-theory}. For simplicity, we assume that the number of groups $K_0$ is fixed. It is however straightforward to allow $K_0$ to grow with the number of time series $n$. We comment on this in more detail in Section \ref{sec-ext}. The groups $G_k = G_{k,n}$ depend on the cross-section dimension $n$ in general. For ease of notation, we however suppress this dependence on $n$ throughout the paper.

\subsection{Identification of the functions $\boldsymbol{m_i}$}\label{subsec-model-ident}

Plugging \eqref{model-eq2} into \eqref{model-eq1}, we obtain the model equation 
\begin{equation}\label{model-eq}
Y_{it} = m_i(X_{it}) + \alpha_i + \gamma_t + \varepsilon_{it},
\end{equation}
where $\ex[\varepsilon_{it}|X_{it}] = 0$. If we drop the fixed effects $\alpha_i$ and $\gamma_t$ from \eqref{model-eq}, we are left with the standard regression equation $Y_{it} = m_i(X_{it}) + \varepsilon_{it}$. Obviously, $m_i$ is identified in this case since $m_i( \, \cdot \, ) = \ex[Y_{it}|X_{it} = \, \cdot \, ]$. In the full model \eqref{model-eq}, in contrast, $m_i$ is not identified. In particular, we can rewrite \eqref{model-eq} as $Y_{it} = \{ m_i(X_{it}) + a_i \} + \{ \alpha_i - a_i \} + \gamma_t + \varepsilon_{it}$, where $a_i$ is an arbitrary real constant. In order to get identification, we need to impose certain constraints which pin down the expectation $\ex[m_i(X_{it})]$ for any $i$ and $t$. We in particular work with the identification constraint that 
\begin{equation}\label{C-identify}
\ex[m_i(X_{it})] = 0 \quad \text{for } 1 \le t \le T \text{ and } 1 \le i \le n. 
\end{equation}
Under this constraint, it is straightforward to show that the functions $m_i$ are identified. In particular, we can derive the following formal result whose proof is given in the Supplementary Material for completeness. 
\begin{prop}\label{prop-identify} 
Let the constraint \eqref{C-identify} be satisfied and suppose that the regularity conditions \ref{C-mixing}--\ref{C-m} from Section \ref{sec-theory} are fulfilled. Then the functions $m_i$ in model \eqref{model-eq} are identified. More precisely, let $m_i$ and $\widetilde{m}_i$ be two functions for some $i \in \{1,\ldots,n\}$ which satisfy the model equation \eqref{model-eq} for any $t$ and which are normalized such that $\ex[m_i(X_{it})] = \ex[\widetilde{m}_i(X_{it})] = 0$ for any $t$. Then $m_i(x) = \widetilde{m}_i(x)$ must hold for all $x \in [0,1]$.   
\end{prop}

Apart from a couple of technicalities, conditions \ref{C-mixing}--\ref{C-m} contain the following two assumptions which are essential for the identification result of Proposition \ref{prop-identify}: 
\begin{enumerate}[label=(\alph*),leftmargin=0.75cm]
\item The time series $\{ X_{it}: t = 1,2,\ldots \}$ is strictly stationary with $X_{it} \sim f_i$ for each $i$.
\item The density $f_i$ is the same for all time series $i$ in a given group $G_k$, that is, $f_i = f_j$ for all $i,j \in G_k$ and any $k$. 
\end{enumerate}
Under (a) and (b), the identification constraint \eqref{C-identify} amounts to a harmless normalization of the functions $m_i$. On the other hand, it is in general not possible to satisfy \eqref{C-identify} without the assumptions (a) and (b): Suppose that (a) is violated and that for some $i$, $X_{it} \sim f_{it}$ with a density $f_{it}$ that differs across $t$. In this case, the constraint \eqref{C-identify} requires that $\int m_i(x) f_{it}(x) dx = 0$ for all $t$. In general, it is however not possible to satisfy the equation $\int m_i(x) f_{it}(x) dx = 0$ simultaneously for all $t$ if the density $f_{it}$ differs across $t$. An analogous problem arises when (b) is violated and the density $f_i$ varies across $i \in G_k$. According to these considerations, the normalization constraint \eqref{C-identify} requires us to impose assumptions (a) and (b). Hence, in order to identify the functions $m_i$ in the presence of a general fixed effects error structure, we need the regressors to satisfy (a) and (b). If we dropped the fixed effects from the model, we could of course do without these assumptions. There is thus a certain trade-off between a general fixed effects error structure and weaker conditions on the regressors.

\setcounter{equation}{0}  
\section{The multiscale distance statistic}\label{sec-dist-measure}

Let $i$ and $j$ be two time series from our sample. In what follows, we construct a test statistic $\widehat{d}_{ij}$ for the null hypothesis $H_0: m_i(x) = m_j(x)$ for all $x \in [0,1]$, that is, for the null hypothesis that $i$ and $j$ belong to the same group $G_k$ for some $1 \le k \le K_0$. We design the statistic $\widehat{d}_{ij}$ in such a way that it does not depend on a specific bandwidth or smoothing parameter. The statistic $\widehat{d}_{ij}$ will serve as a distance measure between the functions $m_i$ and $m_j$ in our clustering algorithm later on.

\subsection{Construction of the multiscale statistic}\label{subsec-dist-measure}

\noindent \textsc{Step 1.} As a first preliminary step, we define a nonparametric estimator $\widehat{m}_{i,h}$ of the function $m_i$, where $h$ denotes the bandwidth. To do so, suppose for a moment that the fixed effects $\alpha_i$ and $\gamma_t$ are known, which implies that the variables $Y_{it}^* = Y_{it} - \alpha_i - \gamma_t$ are known as well. In this case, we can work with the model equation $Y_{it}^* = m_i(X_{it}) + \varepsilon_{it}$ and estimate the function $m_i$ by applying standard nonparametric regression techniques to the sample $\{(Y_{it}^*,X_{it}): 1 \le t\le T \}$. Since $\alpha_i$ and $\gamma_t$ are unobserved in practice, we replace the unknown variables $Y_{it}^*$ by the approximations $\widehat{Y}_{it}^* = Y_{it} - \overline{Y}_i - \overline{Y}_t^{(i)} + \overline{\overline{Y}}^{(i)}$, where
\begin{equation}\label{sample-av-eq}
\overline{Y}_i = \frac{1}{T} \sum_{t=1}^{T} Y_{it}, \quad \overline{Y}_t^{(i)} = \frac{1}{n-1} \sum_{\substack{j=1 \\ j \ne i}}^{n} Y_{jt} \quad \text{and} \quad \overline{\overline{Y}}^{(i)} = \frac{1}{(n-1)T} \sum_{\substack{j=1 \\ j \ne i}}^{n} \sum_{t=1}^T Y_{jt}. 
\end{equation}
With these approximations at hand, we can estimate $m_i$ by applying kernel regression techniques to the constructed sample $\{ (\widehat{Y}_{it}^*,X_{it}): 1 \le t \le T \}$. In particular, we define a local linear kernel estimator of $m_i$ by 
\begin{equation}\label{ll-smoother-eq1}
\widehat{m}_{i,h}(x) = \frac{\sum\nolimits_{t=1}^T \wght_{it}(x,h) \widehat{Y}_{it}^*}{\sum\nolimits_{t=1}^T \wght_{it}(x,h)}, 
\end{equation}
where the weights $\wght_{it}(x,h)$ have the form 
\begin{equation}\label{ll-smoother-eq2}
\wght_{it}(x,h) = \kernel_h(X_{it} - x) \Big\{ S_{i,2}(x,h) - \Big(\frac{X_{it} - x}{h}\Big) S_{i,1}(x,h) \Big\} 
\end{equation}
with $S_{i,\ell}(x,h) = T^{-1} \sum\nolimits_{t=1}^T \kernel_h(X_{it} - x) (\frac{X_{it} - x}{h})^\ell$ for $\ell = 0,1,2$ and $K$ is a kernel function with $\kernel_h(\varphi) = h^{-1} \kernel(\varphi/h)$. Throughout the paper, we assume that the kernel $\kernel$ has compact support $[-C_K,C_K]$ and we set $C_K = 1$ for ease of notation. 
\vspace{10pt}

\noindent \textsc{Step 2.} As an intermediate step in our construction, we set up a bandwidth-dependent test statistic for a somewhat simpler hypothesis than $H_0$. Specifically, we consider the hypothesis $H_{0,x}: m_i(x) = m_j(x)$ for a fixed point $x \in [0,1]$. A test statistic for this problem is given by 
\begin{equation}\label{def-psi}
\widehat{\psi}_{ij}(x,h) = \sqrt{T h} \, \frac{\big( \widehat{m}_{i,h}(x)  - \widehat{m}_{j,h}(x) \big)}{\sqrt{\widehat{\nu}_{ij}(x,h)}},
\end{equation}
where
\begin{equation}\label{def-psi-norm}
\widehat{\nu}_{ij}(x,h) = \left\{ \frac{\widehat{\sigma}_{i,h}^2(x)}{\widehat{f}_{i,h}(x)} + \frac{\widehat{\sigma}_{j,h}^2(x)}{\widehat{f}_{j,h}(x)} \right\} s(x,h)
\end{equation}
is a scaling factor which normalizes the variance of $\widehat{\psi}_{ij}(x,h)$ to be approximately equal to $1$ for sufficiently large $T$. In formula \eqref{def-psi-norm}, $s(x,h) = \{ \int_{-x/h}^{(1-x)/h} \kernel^2(u) [\kappa_2(x,h) - \kappa_1(x,h) u ]^2 du \} / \{ \kappa_0(x,h) \kappa_2(x,h) - \kappa_1(x,h)^2 \}^2$ is a kernel constant with $\kappa_\ell(x,h) = \int_{-x/h}^{(1-x)/h} u^\ell \kernel(u) du$ for $0 \le \ell \le 2$. Moreover, $\widehat{f}_{i,h}(x) = \{\kappa_0(x,h) T\}^{-1}$ $\sum\nolimits_{t=1}^T \kernel_h(X_{it} - x)$ is a boundary-corrected kernel density estimator of $f_i$, where $f_i$ denotes the density of the regressor $X_{it}$ as in Section \ref{subsec-model-ident}, and $\widehat{\sigma}_{i,h}^2(x) = \{ \sum\nolimits_{t=1}^T \kernel_h(X_{it}-x) [ \widehat{Y}_{it}^* - \widehat{m}_{i,h}(X_{it}) ]^2 \} / \{ \sum\nolimits_{t=1}^T \kernel_h(X_{it}-x) \}$ is an estimator of the conditional error variance $\sigma_i^2(x) = \ex[\varepsilon_{it}^2|X_{it}=x]$. If the error terms $\varepsilon_{it}$ are homoskedastic, that is, if $\sigma_i^2(x) \equiv \sigma_i^2 = \ex[\varepsilon_{it}^2]$ for any $x$, we can replace $\widehat{\sigma}_{i,h}^2(x)$ by the simpler estimator $\widehat{\sigma}_{i,h}^2 = T^{-1} \sum\nolimits_{t=1}^T \{ \widehat{Y}_{it}^* - \widehat{m}_{i,h}(X_{it}) \}^2$.

For some of the discussion later on, it is convenient to decompose the statistic $\widehat{\psi}_{ij}(x,h)$ into a bias part $\widehat{\psi}_{ij}^B(x,h)$ and a variance part $\widehat{\psi}_{ij}^V(x,h)$. Standard calculations for kernel estimators yield that 
\begin{equation}\label{psi-expansion} 
\widehat{\psi}_{ij}(x,h) = \widehat{\psi}_{ij}^B(x,h) + \widehat{\psi}_{ij}^V(x,h) + \text{lower order terms}, 
\end{equation}
where  
\[ \widehat{\psi}_{ij}^B(x,h) = \sqrt{Th}\frac{\int_{-x/h}^{(1-x)/h} \{ w_i(u,x,h) m_i(x + hu) - w_j(u,x,h) m_j(x + hu) \} K(u) du}{\sqrt{\widehat{\nu}_{ij}(x,h)}} \]
with $w_i(u,x,h) = \{ \ex[S_{i,2}(x,h)] - \ex[S_{i,1}(x,h)] u \} f_i(x + hu)/ \{\ex[S_{i,0}(x,h)] \ex[S_{i,2}(x,h)] - \ex[S_{i,1}(x,h)]^2 \}$
and 
\[ \widehat{\psi}_{ij}^V(x,h) = \sqrt{T h} \, \frac{\big( \widehat{m}_{i,h}^V(x) - \widehat{m}_{j,h}^V(x) \big)}{\sqrt{\widehat{\nu}_{ij}(x,h)}} \]
with $\widehat{m}_{i,h}^V(x) = \{\sum\nolimits_{t=1}^T \wght_{it}(x,h) (\varepsilon_{it} - \overline{\varepsilon}_t^{(i)} - \overline{m}_t^{(i)}) \}\{ \sum\nolimits_{t=1}^T \wght_{it}(x,h) \}$ as well as $\overline{\varepsilon}_t^{(i)} = (n-1)^{-1} \sum\nolimits_{j=1, j \ne i}^n \varepsilon_{jt}$ and $\overline{m}_t^{(i)} = (n-1)^{-1} \sum\nolimits_{j=1, j \ne i}^n m_j(X_{jt})$. Under the regularity conditions from Section \ref{sec-theory}, it can be shown that $\widehat{\psi}_{ij}^V(x,h) \convd \normal(0,V_{ij})$, where the asymptotic variance $V_{ij}$ is exactly equal to $1$ in the case that $n \rightarrow \infty$ and is approximately equal to $1$ if $n$ is large but bounded. Moreover, under these conditions, the bias term $\widehat{\psi}_{ij}^B(x,h)$ vanishes for any pair of time series $i$ and $j$ that belong to the same class $G_k$, that is, $\widehat{\psi}_{ij}^B(x,h) = 0$ for any $i, j \in G_k$ and $1 \le k \le K_0$.

The variance part $\widehat{\psi}_{ij}^V(x,h)$ captures the stochastic fluctuations of the statistic $\widehat{\psi}_{ij}(x,h)$, whereas $\widehat{\psi}_{ij}^B(x,h)$ can be regarded as a signal which indicates a deviation from the null $H_{0,x}$. The strength of the signal $\widehat{\psi}_{ij}^B(x,h)$ depends on the choice of the bandwidth $h$. To better understand how the signal varies with the bandwidth $h$, suppose that the two functions $m_i$ and $m_j$ differ on the interval $I(x,h_0) = [x-h_0,x+h_0]$ but are the same outside $I(x,h_0)$. The parameter $h_0$ specifies how local the differences between $m_i$ and $m_j$ are. Put differently, it specifies the scale on which $m_i$ and $m_j$ differ: For small/large values of $h_0$, the interval $I(x,h_0)$ is small/large compared to the overall support $[0,1]$, which means that $m_i$ and $m_j$ differ on a local/global scale. Usually, the signal $\widehat{\psi}_{ij}^B(x,h)$ is strongest for bandwidths $h$ close to $h_0$ and becomes weak for bandwidths $h$ that are substantially smaller or larger than $h_0$. The heuristic reason for this is as follows: If $h$ is much larger than $h_0$, the differences between $m_i$ and $m_j$ get smoothed out by the kernel methods that underlie the statistic $\widehat{\psi}_{ij}(x,h)$. If $h$ is much smaller than $h_0$, in contrast, we do not take into account all data points which convey information on the difference between $m_i$ and $m_j$. As a result, the signal $\widehat{\psi}_{ij}^B(x,h)$ gets rather weak. Hence, if the bandwidth $h$ is much smaller/larger than the scale $h_0$ on which $m_i$ and $m_j$ mainly differ, the statistic $\widehat{\psi}_{ij}(x,h)$ is not able to pick up the differences between $m_i$ and $m_j$ and thus to detect a deviation from the null $H_{0,x}$.
\vspace{10pt}

\noindent \textsc{Step 3.} Let us now turn to the problem of testing the hypothesis $H_0: m_i(x) = m_j(x)$ for all $x \in [0,1]$. A simple bandwidth-dependent test statistic for $H_0$ is the supremum statistic
\[ \widehat{d}_{ij}(h) = \sup_{x \in [0,1]} \big| \widehat{\psi}_{ij}(x,h) \big|. \]
Obviously, this statistic suffers from the same problem as the statistic $\widehat{\psi}_{ij}(x,h)$: It is not able to pick up local/global differences between the functions $m_i$ and $m_j$ in a reliable way if the bandwidth $h$ is chosen too large/small. Its performance can thus be expected to strongly depend on the chosen bandwidth.

A simple strategy to get rid of the dependence on the bandwidth $h$ is as follows: We compute the statistic $\widehat{d}_{ij}(h)$ not only for a single bandwidth $h$ but for a wide range of different bandwidths. We in particular consider all bandwidths $h$ in the set $\mathcal{H} = \{ h: h_{\min} \le h \le h_{\max} \}$, where $h_{\min}$ and $h_{\max}$ denote some minimal and maximal bandwidth values that are specified later on. This leaves us with a whole family of statistics $\{ \widehat{d}_{ij}(h) : h \in \mathcal{H} \}$. By taking the supremum over all these statistics, we obtain the rudimentary multiscale statistic
\begin{equation}\label{MS-stat-noadd}
\widetilde{d}_{ij} = \sup_{h \in \mathcal{H}} \widehat{d}_{ij}(h) = \sup_{h \in \mathcal{H}} \, \sup_{x \in [0,1]} \big| \widehat{\psi}_{ij}(x,h) \big|. 
\end{equation}
This statistic does not depend on a specific bandwidth $h$ that needs to be selected. It rather takes into account a wide range of different bandwidths $h \in \mathcal{H}$ simultaneously. It should thus be able to detect differences between the functions $m_i$ and $m_j$ on multiple scales simultaneously. Put differently, it should be able to pick up both local and global differences between $m_i$ and $m_j$.

Inspecting the statistic $\widetilde{d}_{ij}$ more closely, it can be seen to have the following drawback: It does not take into account all scales $h \in \mathcal{H}$ in an equal fashion. Its stochastic behaviour is rather dominated by the statistics $\widehat{\psi}_{ij}(x,h)$ that correspond to small scales $h$. To see this, let us examine the statistic $\widetilde{d}_{ij}$ under the null hypothesis $H_0$, that is, in the case that $i$ and $j$ belong to the same group $G_k$. In this case, $\widehat{\psi}_{ij}(x,h) = \widehat{\psi}_{ij}^V(x,h) + \text{lower order terms}$, since the bias term $\widehat{\psi}_{ij}^B(x,h)$ in \eqref{psi-expansion} is equal to $0$ for all $x$ and $h$ as already noted in Step 2 above. Hence, the statistic $\widehat{\psi}_{ij}(x,h)$ is approximately equal to the variance term $\widehat{\psi}_{ij}^V(x,h)$, which captures its stochastic fluctations. Neglecting terms of lower order, we obtain that under $H_0$, $\widehat{\psi}_{ij}(x,h) = \widehat{\psi}_{ij}^V(x,h)$ and thus  
\[ \widetilde{d}_{ij} = \sup_{h \in \mathcal{H}} \widehat{d}_{ij}(h) \qquad \text{with} \qquad \widehat{d}_{ij}(h) = \sup_{x \in [0,1]} |\widehat{\psi}_{ij}^V(x,h)|. \]
For a given bandwidth $h$, the statistics $\widehat{\psi}_{ij}^V((2\ell-1)h,h)$ for $\ell=1,\ldots,\lfloor 1/2h \rfloor$ are (approximately) standard normal and independent (for sufficiently large $T$). Since the maximum over $\lfloor 1/2h \rfloor$ independent standard normal random variables is $\lambda(2h) + o_p(1)$ as $h \rightarrow 0$ with $\lambda(r) = \sqrt{2 \log(1/r)}$, it holds that $\max_\ell \widehat{\psi}_{ij}^V((2\ell-1)h,h)$ is approximately of size $\lambda(2h)$ for small bandwidths $h$. Moreover, since the statistics $\widehat{\psi}_{ij}^V(x,h)$ with $(2\ell-1) h < x < (2\ell+1) h$ are correlated with $\widehat{\psi}_{ij}^V((2\ell-1)h,h)$ and $\widehat{\psi}_{ij}^V((2\ell+1)h,h)$, the supremum $\sup_x \widehat{\psi}_{ij}^V(x,h)$ approximately behaves as the maximum $\max_\ell \widehat{\psi}_{ij}^V((2\ell-1)h,h)$. Taken together, these considerations suggest that 
\begin{equation}\label{approx-without-correction}
\widehat{d}_{ij}(h) \approx \max_{1 \le \ell \le \lfloor 1/2h \rfloor} \big|\widehat{\psi}_{ij}^V((2\ell-1)h,h)\big| \approx \lambda(2h) 
\end{equation}
for small bandwidth values $h$. According to \eqref{approx-without-correction}, the statistic $\widehat{d}_{ij}(h)$ tends to be much larger in size for small than for large bandwidths $h$. As a consequence, the stochastic behaviour of $\widetilde{d}_{ij}$ tends to be dominated by the statistics $\widehat{d}_{ij}(h)$ which correspond to small bandwidths $h$.

To fix this problem, we follow \cite{DuembgenSpokoiny2001} and replace the statistic $\widetilde{d}_{ij}$ by the modified version 
\begin{equation}\label{MS-stat-add}
\widehat{d}_{ij} = \sup_{h \in \mathcal{H}} \, \sup_{x \in [0,1]} \big\{ |\widehat{\psi}_{ij}(x,h)| - \lambda(2h) \big\}, 
\end{equation}
where $\lambda(r) = \sqrt{2 \log(1/r)}$. For each given bandwidth $h$, we thus subtract the additive correction term $\lambda(2h)$ from the statistics $\widehat{\psi}_{ij}(x,h)$. The idea behind this additive correction is as follows: We can write $\widehat{d}_{ij} = \sup_{h \in \mathcal{H}} \{ \widehat{d}_{ij}(h) - \lambda(2h) \}$ with $\widehat{d}_{ij}(h) = \sup_{x \in [0,1]} |\widehat{\psi}_{ij}(x,h)|$. According to the heuristic considerations from above, when $i$ and $j$ belong to the same class, the statistic $\widehat{d}_{ij}(h)$ is approximately of size $\lambda(2h)$ for small values of $h$. Hence, we correct $\widehat{d}_{ij}(h)$ by subtracting its approximate size under the null hypothesis $H_0$. This calibrates the statistics $\widehat{d}_{ij}(h)$ in such a way that their stochastic fluctuations are comparable across scales $h$. We thus put them on a more equal footing and prevent small scales from dominating the stochastic behaviour of the multiscale statistic. As a result, the statistic $\widehat{d}_{ij}$ should be able to detect differences between the functions $m_i$ and $m_j$ on multiple scales simultaneously without being dominated by a particular scale. It should thus be a reliable test statistic for $H_0$, no matter whether the differences between $m_i$ and $m_j$ are on local or global scales.

To make the statistic $\widehat{d}_{ij}$ defined in \eqref{MS-stat-add} computable in practice, we replace the supremum over $x \in [0,1]$ and $h \in \mathcal{H}$ by the maximum over all points $(x,h)$ in a suitable grid $\mathcal{G}_T$. The final version of the multiscale statistic is thus defined as 
\begin{equation}\label{def-MS-stat}
\widehat{d}_{ij} = \max_{(x,h) \in \mathcal{G}_T} \big\{ |\widehat{\psi}_{ij}(x,h)| - \lambda(2h) \big\}. 
\end{equation}
In this definition, $\mathcal{G}_T$ may be any subset of $\mathcal{G} = \{ (x,h) \, | \, h_{\min} \le h \le h_{\max} \text{ and } x \in [0,1] \}$ with the following properties: (a) $\mathcal{G}_T$ becomes dense in $\mathcal{G}$ as $T \rightarrow \infty$, (b) $|\mathcal{G}_T| \le C T^\beta$ for some arbitrarily large but fixed constants $C,\beta > 0$, where $|\mathcal{G}_T|$ denotes the cardinality of $\mathcal{G}_T$, and (c) $h_{\min} \ge c T^{-(1-\delta)}$ and $h_{\max} \le C T ^{-\delta}$ for some arbitrarily small but fixed $\delta > 0$ and some positive constants $c$ and $C$. According to conditions (a) and (b), the number of points $(x,h)$ in $\mathcal{G}_T$ should grow to infinity as $T \rightarrow \infty$, however it should not grow faster than $CT^\beta$ for some arbitrarily large constants $C, \beta > 0$. This is a fairly weak restriction as it allows the set $\mathcal{G}_T$ to be extremely large as compared to the sample size $T$. As an example, we may use the Wavelet multiresolution grid $\mathcal{G}_T = \{ (x,h) = (2^{-\nu} r, 2^{-\nu}) \, | \, 1 \le r \le 2^{\nu}-1$ and $h_{\min} \le 2^{-\nu} \le h_{\max} \}$. Condition (c) is quite weak as well, allowing us to choose the bandwidth window $[h_{\min},$ $h_{\max}]$ extremely large. In particular, we can choose the minimal bandwidth $h_{\min}$ to converge to zero almost as quickly as the time series length $T$ and thus to be extremely small. Moreover, the maximal bandwidth $h_{\max}$ is allowed to converge to zero very slowly, in particular much more slowly than the optimal bandwidths for estimating the functions $m_i$, which are of the order $T^{-1/5}$ for all $i$ under our technical conditions from Section \ref{sec-theory}. Hence, $h_{\max}$ can be chosen very large.

\subsection{Tuning parameter choice}\label{subsec-dist-measure-bw}

The multiscale statistic $\widehat{d}_{ij}$ does not depend on a specific bandwidth $h$ that needs to be selected. It is thus free of a classical bandwidth or smoothing parameter. However, it is of course not completely free of tuning parameters. It obviously depends on the minimal and maximal bandwidths $h_{\min}$ and $h_{\max}$. Importantly, $h_{\min}$ and $h_{\max}$ are much more harmless tuning parameters than a classical bandwidth $h$. In particular, (a) they are much simpler to choose and (b) the multiscale methods are much less sensitive to their exact choice than conventional methods are to the choice of bandwidth. In what follows, we discuss the reasons for (a) and (b) in detail and give some guidelines how to choose $h_{\min}$ and $h_{\max}$ appropriately in practice. These guidelines are in particular used to implement our methods in the simulations of Section \ref{sec-sim} and the empirical application of Section \ref{sec-app}.

Ideally, we would like to make the interval $[h_{\min},h_{\max}]$ as large as possible, thus taking into account as many scales $h$ as possible. From a technical perspective, we can pick any bandwidths $h_{\min}$ and $h_{\max}$ with $h_{\min} \ge c T^{-(1-\delta)}$ and $h_{\max} \le C T^{-\delta}$ for some small $\delta > 0$. Hence, our theory allows us to choose $h_{\min}$ and $h_{\max}$ extremely small and large, respectively. Heuristically speaking, the bandwidth $h_{\min}$ can be considered very small if the effective sample size $T h_{\min}$ for estimating the functions $m_i$ is very small, say $T h_{\min} \le 10$. Likewise, $h_{\max}$ can be regarded as extremely large if the effective sample size $T h_{\max}$ is very large compared to the full sample size $T$, say $T h_{\max} \approx T/4$ or $T h_{\max} \approx T/3$. Hence, in practice, we have a pretty good idea of what it means for $h_{\min}$ and $h_{\max}$ to be very small and large, respectively. It is thus clear in which range we need to pick the bandwidths $h_{\min}$ and $h_{\max}$ in practice.

As long as the bandwidth window  $[h_{\min},h_{\max}]$ is chosen reasonably large, the exact choice of $h_{\min}$ and $h_{\max}$ can be expected to have little effect on the overall behaviour of the multiscale statistic $\widehat{d}_{ij}$. To see why, write $\widehat{\psi}_{ij}(x,h) = \widehat{\psi}_{ij}^B(x,h) + \widehat{\psi}_{ij}^V(x,h) + \text{lower order terms}$ as in \eqref{psi-expansion}, where the variance term $\widehat{\psi}_{ij}^V(x,h)$ captures the stochastic fluctuations of $\widehat{\psi}_{ij}(x,h)$ and the bias term $\widehat{\psi}_{ij}^B(x,h)$ is a signal which picks up differences between the functions $m_i$ and $m_j$ locally around $x$. Neglecting terms of lower order, the multiscale statistic $\widehat{d}_{ij}$ from \eqref{MS-stat-add} can be written as
\[ \widehat{d}_{ij} = \sup_{h \in [h_{\min},h_{\max}]} \, \sup_{x \in [0,1]} \big\{ |\widehat{\psi}_{ij}^B(x,h) + \widehat{\psi}_{ij}^V(x,h)| - \lambda(2h) \big\}. \]
If the bandwidth window $[h_{\min}, h_{\max}]$ is chosen sufficiently large, it will contain all the scales $h^*$ on which the two functions $m_i$ and $m_j$ mainly differ. As discussed in Section \ref{subsec-dist-measure}, the signals $\widehat{\psi}_{ij}^B(x,h)$ should be strongest for bandwidths $h$ which are close to the scales $h^*$. Hence, as long as the window $[h_{\min}, h_{\max}]$ is chosen large enough to contain all the scales $h^*$, the size of the overall signal of the multiscale statistic $\widehat{d}_{ij}$ should be hardly affected by the exact choice of $h_{\min}$ and $h_{\max}$. Moreover, the size of the stochastic fluctuations of $\widehat{d}_{ij}$ should not be strongly influenced either: The stochastic part of $\widehat{d}_{ij}$ can be expressed as 
\[ \sup_{h \in [h_{\min},h_{\max}]} \widehat{V}_{ij}(h) \quad \text{with} \quad \widehat{V}_{ij}(h) = \sup_{x \in [0,1]} \big\{ |\widehat{\psi}_{ij}^V(x,h)| - \lambda(2h) \big\}, \]
where $\widehat{V}_{ij}(h)$ captures the stochastic fluctuations corresponding to bandwidth $h$. According to our heuristic considerations from Section \ref{subsec-dist-measure}, the variables $\widehat{V}_{ij}(h)$ are comparable in size across bandwidths $h$. Moreover, for $h$ and $h^\prime$ close to each other, $\widehat{V}_{ij}(h)$ and $\widehat{V}_{ij}(h^\prime)$ are strongly correlated. For these reasons, the size of the stochastic part $\sup_{h \in [h_{\min},h_{\max}]} \widehat{V}_{ij}(h)$ should not change much when we make the very large bandwidth window $[h_{\min},h_{\max}]$ somewhat larger or smaller.

In view of these heuristic considerations,  we suggest to choose $h_{\min}$ in practice such that the effective sample size $T h_{\min}$ is small, say $\le 10$, and $h_{\max}$ such that the effective sample size $T h_{\max}$ is large compared to $T$, say $T h_{\max} \ge T/4$.

\subsection{Properties of the multiscale statistic}

We now discuss some theoretical properties of the multiscale statistic $\widehat{d}_{ij}$ which are needed to derive the formal properties of the clustering methods developed in the following sections. Specifically, we compare the maximal multiscale distance between two time series $i$ and $j$ from the same class,
\[ \max_{1 \le k \le K_0} \max_{i, j \in G_k} \, \widehat{d}_{ij}, \]
with the minimal distance between two time series $i$ and $j$ from two different classes, 
\[ \min_{1 \le k < k^\prime \le K_0} \min_{\substack{ i \in G_k, \\ j \in G_{k^\prime} }} \widehat{d}_{ij}. \]
In Section \ref{sec-theory}, we show that under appropriate regularity conditions, 
\begin{align}
\max_{1 \le k \le K_0} \max_{i, j \in G_k} \, \widehat{d}_{ij} & = O_p\big( \sqrt{\log n + \log T} \big) \label{dhat-prop1a} \\
\min_{1 \le k < k^\prime \le K_0} \min_{\substack{ i \in G_k, \\ j \in G_{k^\prime} }} \widehat{d}_{ij} & \ge c_0 \sqrt{Th_{\max}} + o_p \big( \sqrt{Th_{\max}} \big), \label{dhat-prop2a}
\end{align}
where $c_0$ is a sufficiently small positive constant. These two statements imply that 
\begin{align}
\max_{1 \le k \le K_0} \max_{i, j \in G_k} \, \widehat{d}_{ij} \big/ \sqrt{Th_{\max}} & = o_p(1) \label{dhat-prop1b} \\
\min_{1 \le k < k^\prime \le K_0} \min_{\substack{ i \in G_k, \\ j \in G_{k^\prime} }} \widehat{d}_{ij} \big/ \sqrt{Th_{\max}} & \ge c_0 + o_p(1). \label{dhat-prop2b}
\end{align}
According to \eqref{dhat-prop1b} and \eqref{dhat-prop2b}, the maximal distance between time series of the same class converges to zero when normalized by $\sqrt{Th_{\max}}$, whereas the minimal distance between time series of two different classes remains bounded away from zero. Asymptotically, the distance measures $\widehat{d}_{ij}$ thus contain enough information to detect which time series belong to the same class. Technically speaking, we can make the following statement for any fixed positive constant $c < c_0$: with probability tending to $1$, any time series $i$ and $j$ with $\widehat{d}_{ij} \le c$ belong to the same class, whereas those with $\widehat{d}_{ij} > c$ belong to two different classes. The hierarchical clustering algorithm introduced in the next section exploits this information in the distances $\widehat{d}_{ij}$.

\setcounter{equation}{0}
\section{Estimation of the unknown groups}\label{sec-est-classes}

Let $S \subseteq \{1,\ldots,n\}$ and $S^\prime \subseteq \{1,\ldots,n\}$ be two sets of time series from our sample. We define a dissimilarity measure between $S$ and $S^\prime$ by setting 
\begin{equation}\label{dissimilarity}
\widehat{\Delta}(S,S^\prime) = \max_{\substack{i \in S, \\ j \in S^\prime}} \widehat{d}_{ij}. 
\end{equation}
This is commonly called a complete linkage measure of dissimilarity. Alternatively, we may work with an average or a single linkage measure. To partition the set of time series $\{1,\ldots,n\}$ into groups, we combine the multiscale dissimilarity measure $\widehat{\Delta}$ with a hierarchical agglomerative clustering (HAC) algorithm which proceeds as follows: 
\vspace{10pt}

\noindent \textsc{Step $0$ (Initialization):} Let $\widehat{G}_i^{[0]} = \{ i \}$ denote the $i$-th singleton cluster for $1 \le i \le n$ and define $\{\widehat{G}_1^{[0]},\ldots,\widehat{G}_n^{[0]} \}$ to be the initial partition of time series into clusters. 
\vspace{5pt}

\noindent \textsc{Step $r$ (Iteration):} Let $\widehat{G}_1^{[r-1]},\ldots,\widehat{G}_{n-(r-1)}^{[r-1]}$ be the $n-(r-1)$ clusters from the previous step. Determine the pair of clusters $\widehat{G}_{k}^{[r-1]}$ and $\widehat{G}_{k^\prime}^{[r-1]}$ for which 
\[ \widehat{\Delta}(\widehat{G}_{k}^{[r-1]},\widehat{G}_{k^\prime}^{[r-1]}) = \min_{1 \le \ell < \ell^\prime \le n-(r-1)} \widehat{\Delta}(\widehat{G}_{\ell}^{[r-1]},\widehat{G}_{\ell^\prime}^{[r-1]}) \]  
and merge them into a new cluster. 
\vspace{10pt}

\noindent Iterating this procedure for $r = 1,\ldots,n-1$ yields a tree of nested partitions $\{\widehat{G}_1^{[r]},\ldots$ $\ldots,\widehat{G}_{n-r}^{[r]}\}$, which can be graphically represented by a dendrogram. Roughly speaking, the HAC algorithm merges the $n$ singleton clusters $\widehat{G}_i^{[0]} = \{ i \}$ step by step until we end up with the cluster $\{1,\ldots,n\}$. In each step of the algorithm, the closest two clusters are merged, where the distance between clusters is measured in terms of the dissimilarity $\widehat{\Delta}$. We refer the reader to \cite{Ward1963} for an early reference on HAC clustering and to Section 14.3.12 in \cite{HastieTibshiraniFriedman2009} for an overview of hierarchical clustering methods.

We now examine the properties of our HAC algorithm. In particular, we investigate how the partitions $\{\widehat{G}_1^{[r]},\ldots,\widehat{G}_{n-r}^{[r]}\}$ for $r = 1,\ldots,n-1$ are related to the true class structure $\{G_1,\ldots,G_{K_0}\}$. From \eqref{dhat-prop1b} and \eqref{dhat-prop2b}, it immediately follows that the multiscale statistics $\widehat{d}_{ij}$ have the following property: 
\begin{equation}\label{prop-d} 
\pr \Big( \max_{1 \le k \le K_0} \max_{i,j \in G_k} \widehat{d}_{ij} < \min_{1 \le k < k^\prime \le K_0} \min_{\substack{ i \in G_k, \\ j \in G_{k^\prime} }} \widehat{d}_{ij} \Big) \rightarrow 1.  
\end{equation}
To formulate the results on the HAC algorithm, we do not restrict attention to the multiscale statistics $\widehat{d}_{ij}$ from \eqref{def-MS-stat} but let $\widehat{d}_{ij}$ denote any statistics with the high-level property \eqref{prop-d}. We further make use of the following notation: Let $\mathcal{A} = \{ A_1,\ldots,A_r \}$ and $\mathcal{B} = \{ B_1,\ldots,B_{r^\prime} \}$ be two partitions of the set $\{1,\ldots,n\}$, that is, $\mathbin{\dot{\bigcup}}_{\ell = 1}^r A_\ell = \{1,\ldots,n\}$ and $\mathbin{\dot{\bigcup}}_{\ell = 1}^{r^\prime} B_\ell = \{1,\ldots,n\}$. We say that $\mathcal{A}$ is a refinement of $\mathcal{B}$ if each $A_\ell \in \mathcal{A}$ is a subset of some $B_{\ell^\prime} \in \mathcal{B}$. With this notation at hand, the properties of the HAC algorithm can be summarized as follows: 
\vspace{3pt}

\begin{theorem}\label{theo-alg}
Suppose that the statistics $\widehat{d}_{ij}$ satisfy condition \eqref{prop-d}. Then  
\begin{enumerate}[label=\textnormal{(\alph*)},leftmargin=0.75cm]
\item \label{theo-alg-stat1} $\pr \Big(\big\{ \widehat{G}_1^{[n-K_0]},\ldots,\widehat{G}_{K_0}^{[n-K_0]} \big\} = \big\{ G_1,\ldots,G_{K_0} \big\} \Big) \rightarrow 1$,
\item \label{theo-alg-stat2} $\pr \Big( \big\{ \widehat{G}_1^{[n-K]},\ldots,\widehat{G}_{K}^{[n-K]} \big\} \text{ is a refinement of } \big\{ G_1,\ldots,G_{K_0} \big\} \Big) \rightarrow 1 \text{ for any } K > K_0$,  
\item \label{theo-alg-stat3} $\pr \Big( \big\{ G_1,\ldots,G_{K_0} \big\} \text{ is a refinement of } \big\{ \widehat{G}_1^{[n-K]},\ldots,\widehat{G}_{K}^{[n-K]} \big\} \Big) \rightarrow 1 \text{ for any } K < K_0$. 
\end{enumerate}
\end{theorem}
\vspace{3pt}

\noindent The proof of Theorem \ref{theo-alg} is trivial and thus omitted, the statements \ref{theo-alg-stat1}--\ref{theo-alg-stat3} being immediate consequences of condition \eqref{prop-d}. By \ref{theo-alg-stat1}, the partition $\{ \widehat{G}_1,\ldots,\widehat{G}_{K_0} \}$ with $\widehat{G}_k = \widehat{G}_k^{[n-K_0]}$ for $1 \le k \le K_0$ is a consistent estimator of the true class structure $\{ G_1,\ldots, G_{K_0} \}$ in the following sense: $\{ \widehat{G}_1,\ldots,\widehat{G}_{K_0} \}$ coincides with $\{ G_1,\ldots,G_{K_0} \}$ with probability tending to $1$. Hence, if the number of classes $K_0$ were known, we could consistently estimate the true class structure by $\{ \widehat{G}_1,\ldots,\widehat{G}_{K_0}\}$. The partitions $\{ \widehat{G}_1^{[n-K]},\ldots,\widehat{G}_{K}^{[n-K]} \}$ with $K \ne K_0$ can of course not serve as consistent estimators of the true class structure. According to \ref{theo-alg-stat2} and \ref{theo-alg-stat3}, there is nevertheless a close link between these partitions and the unknown class structure. In particular, by \ref{theo-alg-stat2}, for any $K > K_0$,  the estimated clusters $\widehat{G}_1^{[n-K]},\ldots,\widehat{G}_K^{[n-K]}$ are subsets of the unknown classes with probability tending to $1$. Conversely, by \ref{theo-alg-stat3}, for any $K < K_0$, the unknown classes are subsets of the estimated clusters with probability tending to $1$.

\setcounter{equation}{0}
\section{Estimation of the unknown number of groups}\label{sec-est-K0}

\subsection{The estimation method}

Let $\widehat{\Delta}(S,S^\prime)$ be the dissimilarity measure from \eqref{dissimilarity} and define the shorthand $\widehat{\Delta}(S) = \widehat{\Delta}(S,S)$. Moreover, let $\{ \thres_{n,T} \}$ be any sequence with the property that 
\begin{equation}\label{prop-thres}
\sqrt{\log n + \log T} \ll \thres_{n,T} \ll \sqrt{Th_{\max}}, 
\end{equation}
where the notation $a_{n,T} \ll b_{n,T}$ means that $a_{n,T} = o(b_{n,T})$. Combining properties \eqref{dhat-prop1a} and \eqref{dhat-prop2a} of the multiscale distance statistics $\widehat{d}_{ij}$ with the statements of Theorem \ref{theo-alg}, we immediately obtain the following: For any $K < K_0$, 
\begin{equation}\label{thres-iq1}
\pr \Big( \max_{1 \le k \le K} \widehat{\Delta} \big( \widehat{G}_k^{[n-K]} \big) \le \thres_{n,T} \Big) \rightarrow 0, 
\end{equation}
whereas for $K = K_0$, 
\begin{equation}\label{thres-iq2}
\pr \Big( \max_{1 \le k \le K_0}  \widehat{\Delta} \big( \widehat{G}_k^{[n-K_0]} \big) \le \thres_{n,T} \Big) \rightarrow 1. 
\end{equation}
Taken together, \eqref{thres-iq1} and \eqref{thres-iq2} motivate to estimate the unknown number of classes $K_0$ by the smallest number $K$ for which the criterion 
\begin{equation*}
\max_{1 \le k \le K} \widehat{\Delta} \big( \widehat{G}_k^{[n-K]} \big) \le \thres_{n,T} 
\end{equation*}
is satisfied. Formally speaking, we estimate $K_0$ by 
\[ \widehat{K}_0 = \min \Big\{ K = 1,2,\ldots \Big| \max_{1 \le k \le K} \widehat{\Delta} \big( \widehat{G}_k^{[n-K]} \big) \le \thres_{n,T} \Big\}. \]
$\widehat{K}_0$ can be shown to be a consistent estimator of $K_0$ in the sense that $\pr(\widehat{K}_0 = K_0) \rightarrow 1$. More precisely, we can prove the following result. 
\begin{theorem}\label{theo-thres}  
Suppose that the multiscale statistics $\widehat{d}_{ij}$ defined in \eqref{def-MS-stat} have the properties \eqref{dhat-prop1a} and \eqref{dhat-prop2a}. Moreover, let $\{ \thres_{n,T} \}$ be any threshold sequence with the property \eqref{prop-thres}. Then it holds that
$\pr(\widehat{K}_0 = K_0) \rightarrow 1$. 
\end{theorem}
\noindent The proof of Theorem \ref{theo-thres} is straightforward: As already noted, the properties \eqref{dhat-prop1a} and \eqref{dhat-prop2a} of the multiscale distance statistics and the statements of Theorem \ref{theo-alg} immediately imply \eqref{thres-iq1} and \eqref{thres-iq2}. From \eqref{thres-iq1}, it further follows that $\pr(\widehat{K}_0 < K_0) = o(1)$, whereas \eqref{thres-iq2} yields that $\pr(\widehat{K}_0 > K_0) = o(1)$. As a consequence, we obtain that $\pr(\widehat{K}_0 = K_0) \rightarrow 1$.

\begin{figure}[t]
\centering
\includegraphics[width=0.5\textwidth]{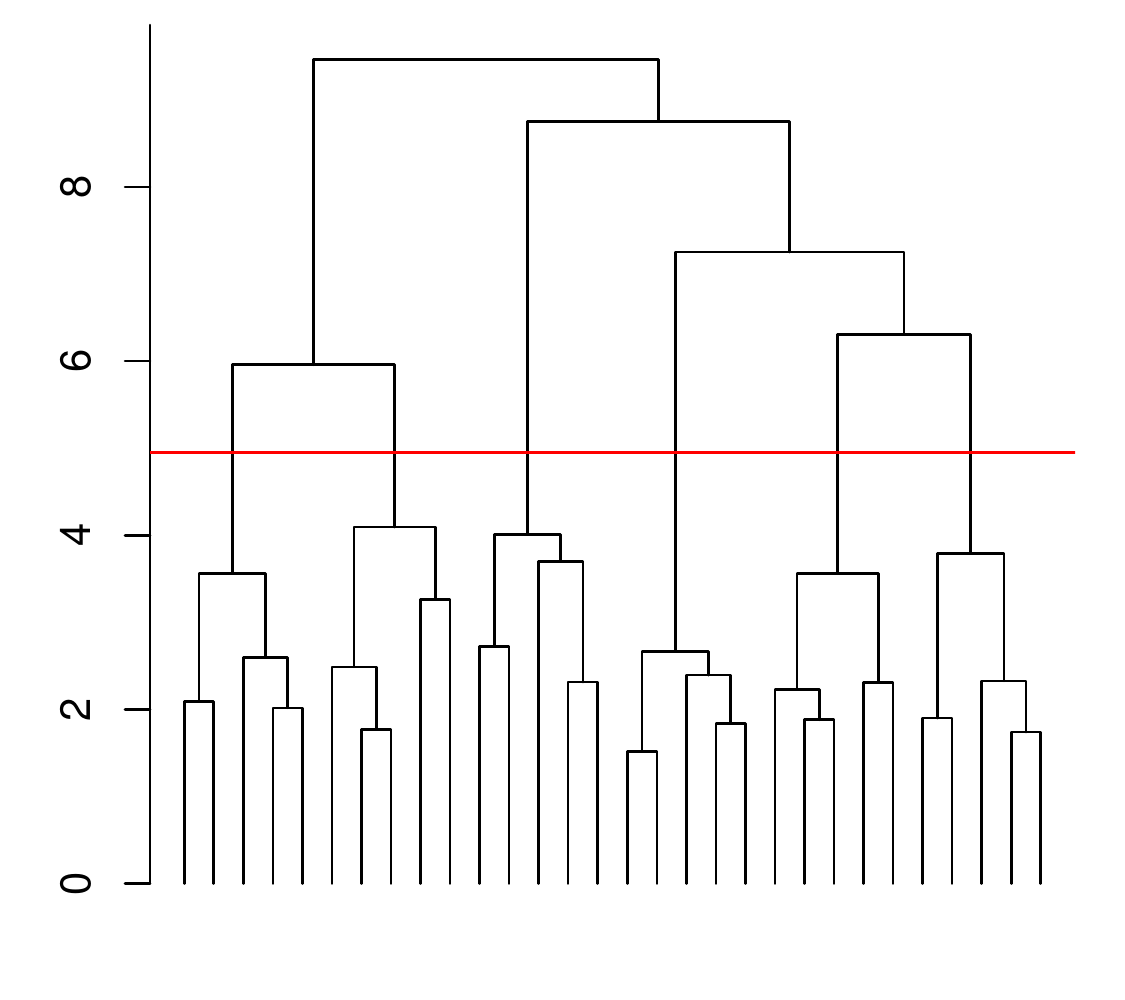} 
\vspace{-0.5cm}

\caption{Example of a dendrogram produced by the HAC algorithm. The red horizontal line indicates the dissimilarity level $\thres_{n,T}$. The estimator $\widehat{K}_0$ can be computed by counting the vertical lines that intersect the red horizontal threshold. In the above example, $\widehat{K}_0$ is equal to $6$.}\label{fig-dendrogram}
\end{figure}

The estimator $\widehat{K}_0$ can be interpreted in terms of the dendrogram produced by the HAC algorithm. It specifies a simple cutoff rule for the dendrogram: The value 
\[ \max_{1 \le k \le K} \widehat{\Delta} \big( \widehat{G}_k^{[n-K]} \big) = \min_{1 \le k < k^\prime \le K+1} \widehat{\Delta} \big( \widehat{G}_k^{[n-(K+1)]}, \widehat{G}_{k^\prime}^{[n-(K+1)]} \big) \]
is the dissimilarity level at which two clusters are merged to obtain a partition with $K$ clusters. In the dendrogram, the clusters are usually indicated by vertical lines and the dissimilarity level at which two clusters are merged is marked by a horizontal line which connects the two vertical lines representing the clusters. To compute the estimator $\widehat{K}_0$, we simply have to cut the dendrogram at the dissimilarity level $\thres_{n,T}$ and count the vertical lines that intersect the horizontal cut at the level $\thres_{n,T}$. See Figure \ref{fig-dendrogram} for an illustration.

\subsection{Choice of the threshold level $\boldsymbol{\thres_{n,T}}$}\label{subsec-choice-thres}

As shown in Theorem \ref{theo-thres}, $\widehat{K}_0$ is a consistent estimator of $K_0$ for any threshold sequence $\{\thres_{n,T}\}$ with the property that $\sqrt{\log n + \log T} \ll \thres_{n,T} \ll \sqrt{Th_{\max}}$. From an asymptotic perspective, we thus have a lot of freedom to choose the thres\-hold $\thres_{n,T}$. In finite samples, a totally different picture arises. There, different choices of $\thres_{n,T}$ may result in markedly different estimates of $K_0$. Selecting the threshold level $\thres_{n,T}$ in a suitable way is thus a crucial issue in finite samples.

In what follows, we give some heuristic discussion on how to pick the threshold level $\thres_{n,T}$ appropriately in practice. To do so, we suppose that the technical conditions from Section \ref{sec-theory} are fulfilled. In addition, we make the simplifying assumption that $\alpha_i = \gamma_t = 0$ for all $i$ and $t$, that is, we drop the fixed effects from the model. Moreover, we suppose that the errors $\varepsilon_{it}$ are homoskedastic and that the error variances $\sigma_i^2 = \ex[\varepsilon_{it}^2]$ are the same within groups. As already discussed in Section \ref{subsec-model-ident}, the densities $f_i$ of the regressors $X_{it}$ are supposed to be the same within groups as well. Slightly abusing notation, we write $\sigma_k^2$ and $f_k$ to denote the group-specific error variance and regressor density in the $k$-th class $G_k$. We can now make the following heuristic observations: 
\begin{enumerate}[label=(\alph*),leftmargin=0.75cm] 

\item \label{obs1} Consider any pair of time series $i$ and $j$ that belong to the same class $G_k$. As in \eqref{psi-expansion}, we decompose $\widehat{\psi}_{ij}(x,h)$ into a bias and a variance part according to $\widehat{\psi}_{ij}(x,h) = \widehat{\psi}_{ij}^B(x,h) + \widehat{\psi}_{ij}^V(x,h) + \text{lower order terms}$. As already noted in Section \ref{subsec-dist-measure}, $\widehat{\psi}_{ij}^B(x,h) = 0$ for $i,j \in G_k$, which implies that 
\begin{equation}\label{heuristic-approx-1}
\widehat{\psi}_{ij}(x,h) \approx \widehat{\psi}_{ij}^V(x,h) = \sqrt{Th} \big\{ \widehat{m}_{i,h}^V(x) - \widehat{m}_{j,h}^V(x) \big\} \big/ \{ \widehat{\nu}_{ij}(x,h) \}^{1/2},   
\end{equation}
where $\widehat{m}_{i,h}^V(x) = \{ \sum\nolimits_{t=1}^T \wght_{it}(x,h) \varepsilon_{it} \} / \{ \sum\nolimits_{t=1}^T \wght_{it}(x,h) \}$ under our simplifying assumptions. Standard arguments for kernel smoothers suggest that
\begin{align} 
\widehat{m}_{i,h}^V(x) 
 & \approx \big\{ f_k(x) \big[ \kappa_0(x,h) \kappa_2(x,h) - \kappa_1(x,h)^2 \big] \big\}^{-1} \nonumber \\
 & \qquad \times \frac{1}{T} \sum\limits_{t=1}^T \kernel_h(X_{it}-x) \Big[ \kappa_2(x,h) - \kappa_1(x,h) \Big( \frac{X_{it}-x}{h} \Big) \Big] \varepsilon_{it}, \label{heuristic-approx-2}
\end{align} 
where $\kappa_\ell(x,h) = \int_{-x/h}^{(1-x)/h} u^\ell K(u) du$ for $0 \le \ell \le 2$. Since by construction, $\widehat{\nu}_{ij}(x,h)$ is an estimator of $\nu_{ij}(x,h) = 2 \{ \sigma_k^2/f_k(x) \} s(x,h)$ with $s(x,h)$ introduced in \eqref{def-psi-norm}, we can combine \eqref{heuristic-approx-1} and \eqref{heuristic-approx-2} to obtain the approximation $\widehat{\psi}_{ij}(x,h) \approx \widehat{\psi}_i(x,h) - \widehat{\psi}_j(x,h)$ with
\begin{align*}
\widehat{\psi}_i(x,h) 
 & = \big\{ 2 \rho(x,h) \sigma_k^2 f_k(x) \big\}^{-{1/2}} \\  
 & \qquad \times \frac{1}{\sqrt{Th}} \sum\limits_{t=1}^T \kernel\Big(\frac{X_{it} - x}{h}\Big) \Big[ \kappa_2(x,h) - \kappa_1(x,h) \Big( \frac{X_{it}-x}{h} \Big) \Big] \varepsilon_{it}, 
\end{align*}
where we use the shorthand  $\rho(x,h) = \int_{-x/h}^{(1-x)/h} \kernel^2(u) [\kappa_2(x,h) - \kappa_1(x,h) u ]^2 du$. For each $i$, we stack the random variables $\widehat{\psi}_i(x,h)$ with $(x,h) \in \mathcal{G}_T$ in the vector
\[ \widehat{\boldsymbol{\psi}}_i = 
\left( \widehat{\psi}_i \big(x_1^1,h_1\big), \dots, \widehat{\psi}_i \big(x_1^{N_1},h_1\big), \dots\dots, \widehat{\psi}_i \big(x_p^1,h_p\big), \dots, \widehat{\psi}_i \big(x_p^{N_p},h_p\big) \right)^\top, 
\]
where $\mathcal{G}_T = \bigcup_{\nu=1}^p \mathcal{G}_{T,\nu}$ and $\mathcal{G}_{T,\nu} = \{ (x_\nu^{\ell},h_\nu): 1 \le \ell \le N_\nu\}$ is the set of points corresponding to the bandwidth level $h_\nu$. Moreover, we write $\boldsymbol{\lambda} = (\boldsymbol{\lambda}_1, \dots, \boldsymbol{\lambda}_p)^\top$ with $\boldsymbol{\lambda}_\nu = (\lambda(2h_\nu),\ldots,\lambda(2h_\nu))$ being a vector of length $N_\nu$ for each $\nu$ and we introduce the notation $|z| = (|z_1|,\ldots,|z_q|)^\top$ and $(z)_{\infty} = \max_{1 \le \ell \le q} z_\ell$ for $z \in \reals^q$. With this notation at hand, we obtain that 
\[ \widehat{d}_{ij} \approx \big( \, |\widehat{\boldsymbol{\psi}}_i - \widehat{\boldsymbol{\psi}}_j| - \boldsymbol{\lambda} \, \big)_{\infty} \]
for any pair of time series $i$ and $j$ that belong to the same class.

\item \label{obs2} For any fixed number of points $z_1,\ldots,z_q \in (0,1)$ and related bandwidths $h_{z_\ell}$ with $h_{\min} \le h_{z_\ell} \le h_{\max}$ for $1 \le \ell \le q$, the random vector $[ \, \widehat{\psi}_i(z_1,h_{z_1}),\ldots,\widehat{\psi}_i(z_q,h_{z_q}) \, ]^\top$ is asymptotically normal. Hence, the random vector $\widehat{\boldsymbol{\psi}}_i$ can be treated as approximately Gaussian for sufficiently large sample sizes. More specifically, since 
\begin{align}
\cov & \big( \widehat{\psi}_i(x,h), \widehat{\psi}_i(x^\prime,h^\prime) \big) \nonumber \\*
     & \approx \big\{ 2 \sqrt{\rho(x,h) \rho(x^\prime,h^\prime)} \big\}^{-1} \sqrt{\frac{h}{h^\prime}} \Big\{ \int_{-x/h}^{(1-x)/h} \kernel(u) \big[ \kappa_2(x,h) - \kappa_1(x,h) u \big] \nonumber \\*
 & \qquad \times \kernel\Big(\frac{hu + x-x^\prime}{h^\prime}\Big)  \Big[ \kappa_2(x^\prime,h^\prime) - \kappa_1(x^\prime,h^\prime) \Big(\frac{hu + x - x^\prime}{h^\prime} \Big) \Big] du \Big\}, \label{cov-psi} 
\end{align}
we can approximate the random vector $\widehat{\boldsymbol{\psi}}_i$ by a Gaussian vector with the covariance structure specified on the right-hand side of \eqref{cov-psi}. 
Moreover, since the vectors $\widehat{\boldsymbol{\psi}}_i$ are independent across $i$ under our assumptions, we can approximate the distribution of 
\[ \max_{i,j \in S} \big( \, |\widehat{\boldsymbol{\psi}}_i - \widehat{\boldsymbol{\psi}}_j| - \boldsymbol{\lambda} \, \big)_{\infty} \]
by that of 
\[ \max_{i,j \in S} \big( \, |\boldsymbol{\zeta}_i - \boldsymbol{\zeta}_j| - \boldsymbol{\lambda} \, \big)_{\infty} \]
for any $S \subseteq \{1,\ldots,n\}$, where $\boldsymbol{\zeta}_i$ are independent Gaussian random vectors with the covariance structure from \eqref{cov-psi}.  

\end{enumerate}

\noindent Ideally, we would like to tune the threshold level $\thres_{n,T}$ such that $\widehat{K}_0 = K_0$ with high probability. Put differently, we would like to choose $\thres_{n,T}$ such that it is slightly larger than $\max_{1 \le k \le K_0} \widehat{\Delta} (\widehat{G}_k^{[n-K_0]})$ with high probability. With the help of the observations \ref{obs1} and \ref{obs2} as well as some further heuristic arguments, this can be achieved as follows: Since the partition $\{ \widehat{G}_1^{[n-K_0]},\ldots,\widehat{G}_{K_0}^{[n-K_0]} \}$ consistently estimates the class structure $\{G_1,\ldots,G_{K_0}\}$, we have that 
\begin{equation}\label{heur-tau-1}
\max_{1 \le k \le K_0} \widehat{\Delta} (\widehat{G}_k^{[n-K_0]}) \approx \max_{1 \le k \le K_0} \widehat{\Delta} (G_k). 
\end{equation}
By observation \ref{obs1}, we further obtain that  
\begin{align}
\max_{1 \le k \le K_0} \widehat{\Delta} (G_k) 
 & = \max_{1 \le k \le K_0} \Big\{ \max_{i,j \in G_k} \widehat{d}_{ij} \Big\} \nonumber \\
 & \approx \max_{1 \le k \le K_0} \Big\{ \max_{i,j \in G_k} \big( \, |\widehat{\boldsymbol{\psi}}_i - \widehat{\boldsymbol{\psi}}_j| - \boldsymbol{\lambda} \, \big)_{\infty} \Big\}, \label{heur-tau-2}
\end{align}
and by \ref{obs2},
\begin{equation}\label{heur-tau-3}
\max_{1 \le k \le K_0} \Big\{ \max_{i,j \in G_k} \big( |\widehat{\boldsymbol{\psi}}_i - \widehat{\boldsymbol{\psi}}_j| - \boldsymbol{\lambda} \big)_{\infty} \Big\} \stackrel{d}{\approx} \max_{1 \le k \le K_0} \Big\{ \max_{i,j \in G_k} \big( \, |\boldsymbol{\zeta}_i - \boldsymbol{\zeta}_j| - \boldsymbol{\lambda} \, \big)_{\infty} \Big\}, 
\end{equation}
where $Z \stackrel{d}{\approx} Z^\prime$ means that $Z$ is approximately distributed as $Z^\prime$. Since the right-hand side of \eqref{heur-tau-3} depends on the unknown groups $G_1,\ldots,G_{K_0}$, we apply the trivial bound
\begin{align}
\max_{1 \le k \le K_0} \Big\{ \max_{i,j \in G_k} & \big( \, |\boldsymbol{\zeta}_i - \boldsymbol{\zeta}_j| - \boldsymbol{\lambda} \, \big)_{\infty} \Big\} \nonumber \\ & \le B_n := \max_{1 \le i,j \le n} \big( \, |\boldsymbol{\zeta}_i - \boldsymbol{\zeta}_j| - \boldsymbol{\lambda} \, \big)_{\infty} \label{heur-tau-4}
\end{align}
and define $q_n(\alpha)$ to be the $\alpha$-quantile of $B_n$. Taken together, \eqref{heur-tau-1}--\eqref{heur-tau-4} suggest that
\[ \max_{1 \le k \le K_0} \widehat{\Delta} (\widehat{G}_k^{[n-K_0]}) \le q_n(\alpha) \] 
holds with high probability if we pick $\alpha$ close to $1$. In particular, if the random variable $\max_{1 \le k \le K_0} \widehat{\Delta}(\widehat{G}_k^{[n-K_0]})$ is not only approximately but exactly distributed as $\max_{1 \le k \le K_0} \max_{i,j \in G_k} ( \, |\boldsymbol{\zeta}_i - \boldsymbol{\zeta}_j| - \boldsymbol{\lambda} \, )_{\infty}$, then 
\[ \pr \Big( \max_{1 \le k \le K_0} \widehat{\Delta} (\widehat{G}_k^{[n-K_0]}) \le q_n(\alpha) \Big) \ge \alpha. \] 
According to these considerations, $\thres_{n,T} = q_n(\alpha)$ with $\alpha$ close to $1$ should be an appropriate threshold level. Throughout the simulations and applications, we set $\alpha = 0.95$.

\setcounter{equation}{0}
\section{Theoretical results}\label{sec-theory}

In this section, we derive the statements \eqref{dhat-prop1a} and \eqref{dhat-prop2a} under appropriate regularity conditions. These statements characterize the convergence behaviour of the multiscale statistics $\widehat{d}_{ij}$ and underlie Theorems \ref{theo-alg} and \ref{theo-thres} which describe the theoretical properties of our clustering methods. To prove \eqref{dhat-prop1a} and \eqref{dhat-prop2a}, we impose the following conditions. 
\begin{enumerate}[label=(C\arabic*),leftmargin=1cm]

\item \label{C-mixing} The time series processes $\mathcal{P}_i = \{ (X_{it},\varepsilon_{it}): t = 1,2,\ldots \}$ are independent across $i$. Moreover, they are strictly stationary and strongly mixing for each $i$. Let $\alpha_i(\ell)$ for $\ell=1,2,\ldots$ be the mixing coefficients corresponding to the $i$-th time series $\mathcal{P}_i$. It holds that $\alpha_i(\ell) \le \alpha(\ell)$ for all $i$, where the coefficients $\alpha(\ell)$ decay exponentially fast to zero as $\ell \rightarrow \infty$. 

\item \label{C-dens} For each $1 \le i \le n$, the random variables $X_{it}$ have a density $f_i$ with the following properties: (a) $f_i$ has bounded support, which w.l.o.g.\ equals $[0,1]$ for all $i$, (b) $f_i$ is bounded away from zero and infinity on $[0,1]$ uniformly over $i$, that is, $0 < c \le f_i(x) \le C < \infty$ for all $x \in [0,1]$ with some constants $c$ and $C$ that neither depend on $x$ nor on $i$, (c) $f_i$ is twice continuously differentiable on $[0,1]$ with first and second derivatives that are bounded away from infinity in absolute value uniformly over $i$. Moreover, the variables $(X_{it},X_{it+\ell})$ have a joint density $f_{i,\ell}$ which is bounded away from infinity uniformly over $i$, that is, $f_{i,\ell}(x,x^\prime) \le C < \infty$ for all $i$, $x$, $x^\prime$ and $\ell$, where the constant $C$ neither depends on $i$, $x$, $x^\prime$ nor on $\ell$.

\item \label{C-sigma} The error terms $\varepsilon_{it}$ are homoskedastic, that is, $\sigma_i^2 = \ex[\varepsilon_{it}^2] = \ex[\varepsilon_{it}^2|X_{it}=x]$ for all $x \in [0,1]$. The error variances $\sigma_i^2$ are uniformly bounded away from zero and infinity, that is, $0 < c \le \sigma_i^2 \le C < \infty$ for all $i$, where the constants $c$ and $C$ do not depend on $i$. 

\item \label{C-group-homo} The densities $f_i$ and the error variances $\sigma_i^2$ are the same within groups. That is, for any $k$ with $1 \le k \le K_0$, it holds that $f_i = f_j$ and $\sigma_i^2 = \sigma_j^2$ for all $i,j \in G_k$. 

\item \label{C-moments} There exist a real number $\theta > 4$ and a natural number $\ell^*$ such that for any $\ell \in \integers$ with $|\ell| \ge \ell^*$ and some constant $C < \infty$, 
\begin{align*}
 & \max_{1 \le i \le n} \sup_{x \in [0,1]} \ex \big[ |\varepsilon_{it}|^\theta \big| X_{it} = x \big] \le C < \infty \\
 & \max_{1 \le i \le n} \sup_{x,x^\prime \in [0,1]} \ex \big[ |\varepsilon_{it} \varepsilon_{it+\ell}| \big| X_{it} = x, X_{it+\ell} = x^\prime \big] \le C < \infty.
\end{align*} 

\item \label{C-m} The group-specific regression functions $g_k$ are twice continuously differentiable on $[0,1]$ for $1 \le k \le K_0$ with Lipschitz continuous second derivatives $g_k^{\prime\prime}$, that is, $|g_k^{\prime\prime}(v) - g_k^{\prime\prime}(w)| \le L |v - w|$ for any $v,w \in [0,1]$ and some constant $L$. Moreover, for any pair of indices $(k,k^\prime)$ with $1 \le k < k^\prime \le K_0$, the functions $g_k$ and $g_{k^\prime}$ are different in the sense that $g_k(x) \ne g_{k^\prime}(x)$ for some point $x \in [0,1]$. 

\item \label{C-nT} It holds that  
\begin{equation}\label{condition-n}
n = n(T) \le C \frac{( T^{1/2} \land Th_{\min} )^{\frac{\theta - \delta}{2}}}{T^{1+\delta}} 
\end{equation}
for some small $\delta > 0$ and a sufficiently large constant $C > 0$, where we use the notation $a \land b = \min \{a,b\}$ and $\theta$ is defined in \ref{C-moments}.
 
\item \label{C-h} The minimal and maximal bandwidths have the form $h_{\min} = a T^{-B}$ and $h_{\max} = A T^{-b}$ with some positive constants $a$, $A$, $b$ and $B$, where $0 < b \le B < 1$.

\item \label{C-ker} The kernel $\kernel$ is non-negative, bounded and integrates to one. Moreover, it is symmetric about zero, has compact support $[-1,1]$ and fulfills the Lipschitz condition that $|\kernel(v) - \kernel(w)| \le L|v-w|$ for some $L$ and all $v,w \in \reals$. 

\end{enumerate}

\begin{remark}
{\normalfont We briefly comment on the above assumptions. 
\begin{enumerate}[label=(\roman*),leftmargin=1cm] 

\item \ref{C-mixing} imposes some weak dependence conditions on the variables $(X_{it},\varepsilon_{it})$ across $t$ in the form of mixing assumptions. Note that we do not necessarily require exponentially decaying mixing rates as assumed in \ref{C-mixing}. These could alternatively be replaced by sufficiently high polynomial rates. We nevertheless make the stronger assumption of exponential mixing to keep the proofs as clear as possible. \ref{C-mixing} further restricts the regressors $X_{it}$ and the errors $\varepsilon_{it}$ to be independent across $i$. Some restricted types of cross-sectional dependence in the data are however possible via the fixed effect error terms $\alpha_i$ and $\gamma_t$. 

\item The homoskedasticity assumption in \ref{C-sigma} as well as the condition in \ref{C-group-homo} that the error variances $\sigma_i^2$ are the same within groups are not necessarily needed but are imposed for simplicity. The restriction in \ref{C-group-homo} that the densities $f_i$ are the same within groups, in contrast, is required for identification purposes as already discussed in Section \ref{subsec-model-ident}. 

\item \ref{C-dens}, \ref{C-moments} and \ref{C-m} are standard moment, boundedness and smoothness conditions to derive uniform convergence results for the kernel estimators on which the multiscale statistics $\widehat{d}_{ij}$ are based; see \cite{Hansen2008} for similar assumptions. 

\item \ref{C-nT} imposes restrictions on the growth of the number of time series $n$. Loosely speaking, it says that $n$ is not allowed to grow too quickly in comparison to $T$. More specifically, let $h_{\min} = a T^{-B}$ with some $B \le 1/2$ and $h_{\max} = A T^{-b}$ with some $b > 0$. In this case, \eqref{condition-n} simplifies to $n \le C T^{(\theta - 4 - 5 \delta)/4}$ with some small $\delta > 0$. This shows that the growth restriction \eqref{condition-n} on $n$ is closely related to the moment conditions on the error terms $\varepsilon_{it}$ in \ref{C-moments}. In particular, the larger the value of $\theta$, that is, the stronger the moment conditions on $\varepsilon_{it}$, the faster $n$ may grow in comparison to $T$. If $\theta = 8$, for example, then $n$ may grow (almost) as quickly as $T$. If $\theta$ can be picked arbitrarily large, that is, if all moments of $\varepsilon_{it}$ exist, then $n$ may grow as quickly as any polynomial of $T$, that is, $n \le CT^\rho$ with $\rho > 0$ as large as desired. 

\item \ref{C-h} imposes some conditions on the minimal and maximal bandwidths $h_{\min}$ and $h_{\max}$. Specifically, it requires that $h_{\min} \ge c T^{-(1-\delta)}$ and $h_{\max} \le C T^{-\delta}$ for some small $\delta > 0$ and positive constants $c$ and $C$. These conditions are fairly weak as already discussed in Section \ref{sec-dist-measure}: According to them, we can choose $h_{\min}$ to converge to zero extremely fast, in particular much faster than the optimal bandwidths for estimating the functions $m_i$, which are of the order $T^{-1/5}$ for any $i$ under the smoothness conditions \ref{C-dens} and \ref{C-m}. Similarly, we can let $h_{\max}$ converge to zero much more slowly than the optimal bandwidths. Hence, we can choose the interval $[h_{\min},h_{\max}]$ to be very large, allowing for both substantial under- and oversmoothing.

\item Finally, it is worth noting that our assumptions do not impose any restrictions on the class sizes $|G_k|$. The sizes $|G_k|$ may thus be very different across the classes $G_k$. In particular, they may be fixed for some classes and grow to infinity at different rates for others.

\end{enumerate}}
\end{remark}

Under the regularity conditions just discussed, we can derive the following result whose proof is provided in the Supplementary Material. 
\begin{theorem}\label{theo-dist}
Under \ref{C-mixing}--\ref{C-ker}, it holds that 
\begin{align}
\max_{1 \le k \le K_0} \max_{i, j \in G_k} \, \widehat{d}_{ij} & = O_p\big( \sqrt{\log n + \log T} \big) \label{theo-dist-stat1} \\
\min_{1 \le k < k^\prime \le K_0} \min_{\substack{ i \in G_k, \\ j \in G_{k^\prime} }} \widehat{d}_{ij} & \ge c_0 \sqrt{Th_{\max}} + o_p \big( \sqrt{Th_{\max}} \big), \label{theo-dist-stat2} 
\end{align}
where $c_0$ is a fixed positive constant that does not depend on $T$ \textnormal{(}nor on $n = n(T)$\textnormal{)}. 
\end{theorem}

\setcounter{equation}{0}
\section{Simulations}\label{sec-sim}

In this section, we carry out some simulations to illustrate the advantages of our multiscale approach over clustering methods that depend on a specific bandwidth. When the grid $\mathcal{G}_T$ of location-scale points $(x,h)$ comprises only one bandwidth value $h$, our multiscale approach reduces to a bandwidth-dependent procedure. Specifically, the resulting procedure consists in applying a hierarchical clustering algorithm to the supremum distances $\widehat{d}_{ij}(h) = \max_{x \in \mathcal{X}} |\widehat{\psi}_{ij}(x,h)|$, where $\mathcal{X}$ is the set of locations under consideration and $h$ is the chosen bandwidth.\footnote{Note that the additive correction term $\lambda(2h)$ can be dropped as it is a fixed constant when only one bandwidth value $h$ is considered.} In what follows, we compare our multiscale approach with this bandwidth-dependent procedure for se\-ve\-ral bandwidth values $h$.

\begin{figure}[t]
\includegraphics[width=\textwidth]{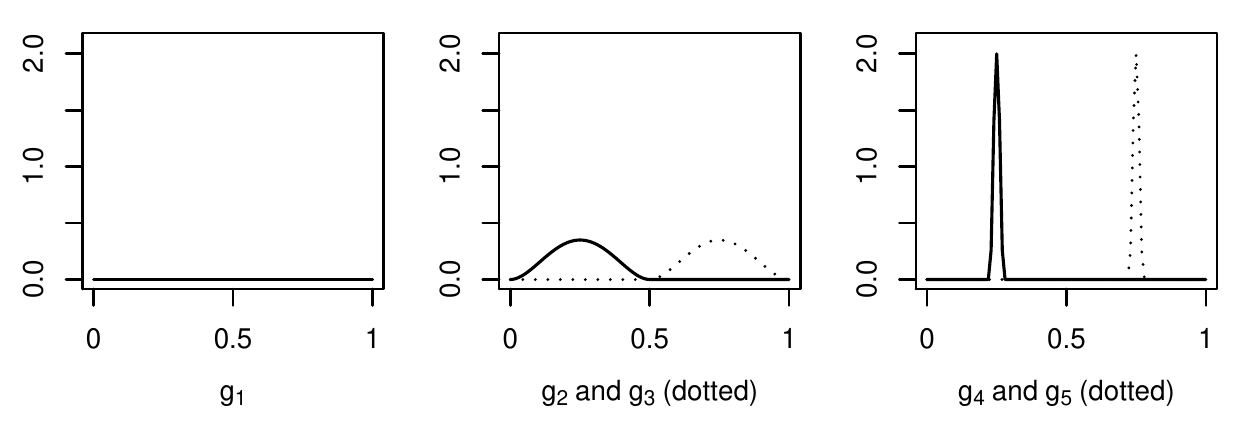} 
\caption{Plot of the functions $g_k$ for $1 \le k \le 5$.}\label{fig1-sim}
\end{figure}

We consider the following setup for the simulations: The data are drawn from the model
\begin{equation}\label{model-sim}
Y_{it} = m_i(X_{it}) + \varepsilon_{it} \quad (1 \le t \le T, \, 1 \le i \le n), 
\end{equation}
where $T = 1000$ and $n = 100$. The time series $i \in \{1,\ldots,n\}$ belong to $K_0=5$ different groups $G_1,\ldots,G_{K_0}$ of the same size. In particular, we set $G_k = \{ (k-1)n/5 + 1,\ldots, kn/5 \}$ for $1 \le k \le K_0 = 5$. The group-specific regression functions $g_k: [0,1] \rightarrow \mathbb{R}$ are given by $g_1(x) = 0$ and 
\begin{align*} 
g_2(x) = 0.35 \, b\big(x,\textstyle{\frac{1}{4}},\textstyle{\frac{1}{4}}\big) & \qquad g_4(x) = 2 \, b\big(x,\textstyle{\frac{1}{4}},\textstyle{\frac{1}{40}}\big) \\
g_3(x) = 0.35 \, b\big(x,\textstyle{\frac{3}{4}},\textstyle{\frac{1}{4}}\big) & \qquad g_5(x) = 2 \, b\big(x,\textstyle{\frac{3}{4}},\textstyle{\frac{1}{40}}\big), 
\end{align*} 
where $b(x,x_0,h) = 1(|x-x_0|/h \le 1) \, \{1 - ((x-x_0)/h)^2\}^2$. Figure \ref{fig1-sim} provides a graphical illustration of the functions $g_k$ for $1 \le k \le 5$. The error process $\mathcal{E}_i = \{ \varepsilon_{it}: 1 \le t \le T \}$ has an autoregressive (AR) structure for each $i$, in particular $\varepsilon_{it} = a \varepsilon_{it-1} + \eta_{it}$ for $1 \le t \le T$, where $a$ is the AR parameter and the innovations $\eta_{it}$ are i.i.d.\ normal with $\ex[\eta_{it}] = 0$ and $\ex[\eta_{it}^2] = \nu^2$. We consider two different values for the AR parameter $a$, in particular $a = -0.25$ and $a=0.25$. The innovation variance $\nu^2$ is chosen as $\nu^2 = 1-a^2$, which implies that $\var(\varepsilon_{it}) = 1$. The regressors $X_{it}$ are drawn independently from a uniform distribution on $[0,1]$ for each $i$. As can be seen, there is no time series dependence in the regressors, and we do not include fixed effects $\alpha_i$ and $\gamma_t$ in the model. We do not take into account these complications because the main aim of the simulations is to display the advantages of our multiscale approach over bandwidth-dependent procedures. These advantages can be seen most clearly in a simple stylized simulation setup as the one under consideration.

To implement our multiscale approach, we use the location-scale grid $\mathcal{G}_T = \{ (x,h) : x \in \mathcal{X} \text{ and } h \in \mathcal{H} \}$, where $\mathcal{X} = \{ x : x = r/100 \text{ for } r = 5,\ldots,95 \}$ is the set of locations and $\mathcal{H} = \{ h : 0.025 \le h \le 0.25 \text{ with } h = 0.025 k \text{ for } k = 1,2,\ldots \}$ is the set of bandwidths. The bandwidth-dependent algorithm is implemented with the same set of locations $\mathcal{X}$ and five different bandwidth values $h$, in particular $h \in \{ 0.025,0.05,0.1,0.2,0.25 \}$. The number of classes $K_0 = 5$ is estimated as described in Section \ref{sec-est-K0} both when the multiscale and the bandwidth-dependent algorithm is used. The threshold parameter $\pi_{n,T}$ is set to $\pi_{n,T} = q_n(\alpha)$ with $\alpha = 0.95$. To produce our simulation results, we draw $S=1000$ samples from model \eqref{model-sim} and compute the estimates of the classes $G_1,\ldots,G_{K_0}$ and their number $K_0$ for each simulated sample both for the multiscale and the bandwidth-dependent algorithm.

\begin{figure}[p]
\phantom{upper margin}
\vspace{-1.25cm}

\begin{subfigure}[t]{\textwidth}
\caption{Histograms of the number of classification errors $\#F$}\label{fig2a-sim}
\vspace{0.05cm}

\includegraphics[width=\textwidth]{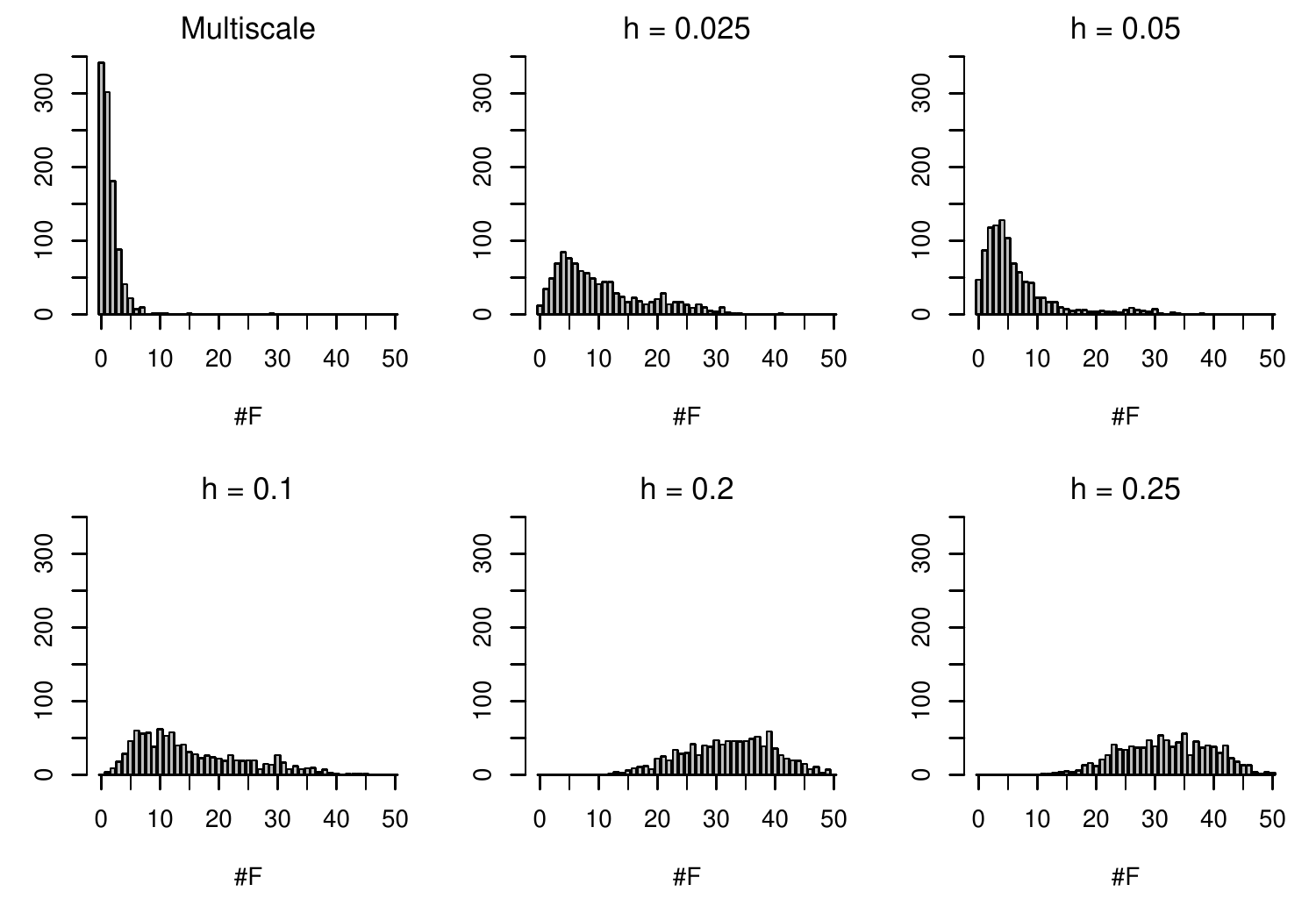} 
\end{subfigure}
\vspace{-0.4cm}

\begin{subfigure}[t]{\textwidth}
\caption{Histograms of the estimated number of clusters $\widehat{K}_0$}\label{fig2b-sim}
\vspace{0.05cm}

\includegraphics[width=\textwidth]{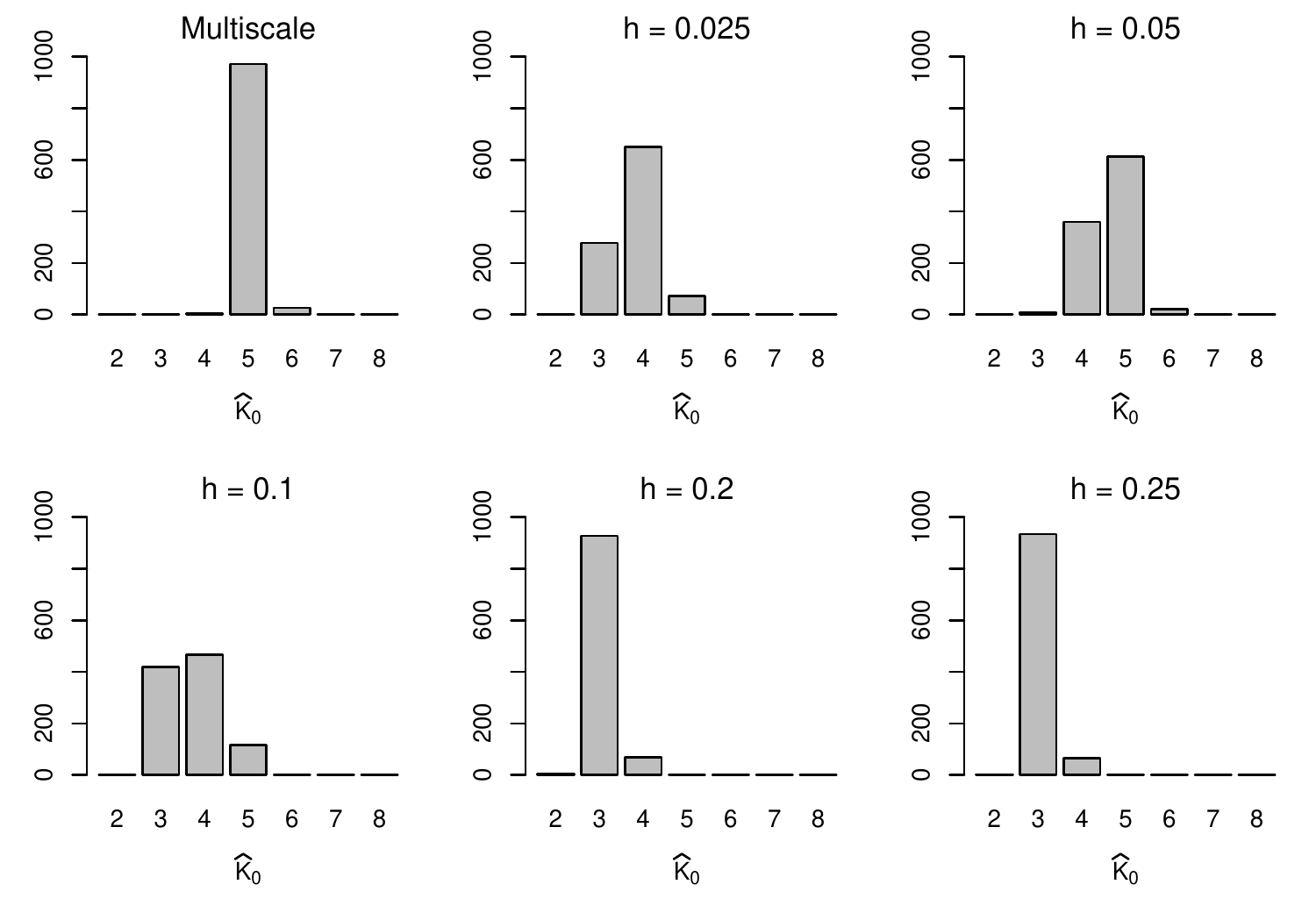} 
\end{subfigure}
\vspace{-0.3cm}

\caption{Simulation results for the design with the negative AR parameter $a = -0.25$. In both subfigures (a) and (b), the upper left panel shows the results for our multiscale approach and the other panels those for the bandwidth-dependent competitor with different bandwidths $h$.}\label{fig2-sim}
\end{figure}

\begin{figure}[p]
\phantom{upper margin}
\vspace{-1.25cm}

\begin{subfigure}[t]{\textwidth}
\caption{Histograms of the number of classification errors $\#F$}\label{fig3a-sim}
\vspace{0.05cm}

\includegraphics[width=\textwidth]{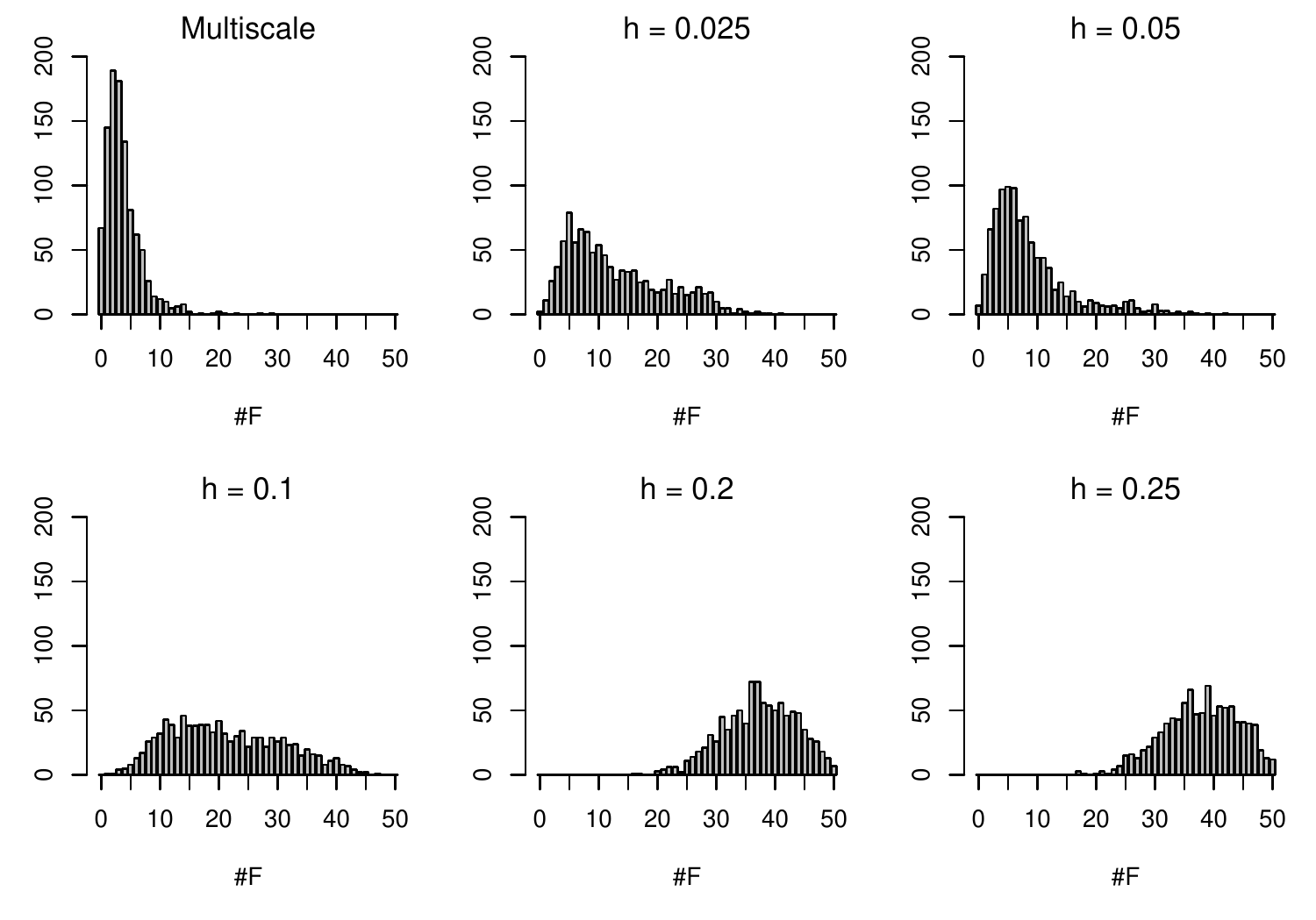} 
\end{subfigure}
\vspace{-0.4cm}

\begin{subfigure}[t]{\textwidth}
\caption{Histograms of the estimated number of clusters $\widehat{K}_0$}\label{fig3b-sim}
\vspace{0.05cm}

\includegraphics[width=\textwidth]{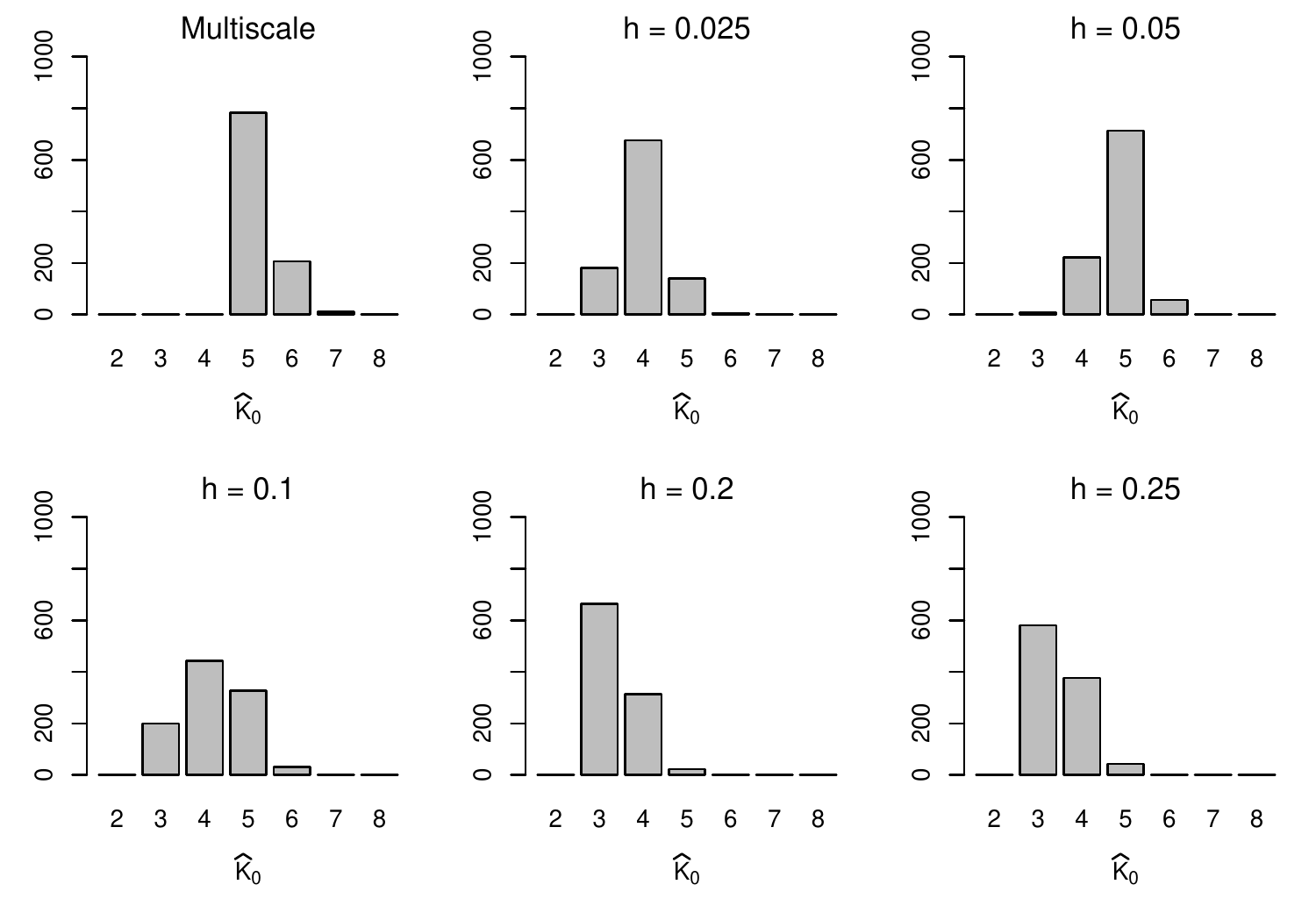} 
\end{subfigure}
\vspace{-0.3cm}

\caption{Simulation results for the design with the positive AR parameter $a = 0.25$. In both subfigures (a) and (b), the upper left panel shows the results for our multiscale approach and the other panels those for the bandwidth-dependent competitor with different bandwidths $h$.}\label{fig3-sim}
\end{figure}

The simulation results for the scenario with the negative AR parameter $a = -0.25$ are reported in Figure \ref{fig2-sim} and those for the scenario with the positive parameter $a = 0.25$ in Figure \ref{fig3-sim}. We first have a closer look at the results in Figure \ref{fig2-sim}. To produce Figure \ref{fig2a-sim}, we treat $K_0$ as known and compute the number of classification errors $\#F$, that is, the number of wrongly classified indices $i$ for each of the $S=1000$ simulated samples.\footnote{Precisely speaking, $\#F$ is defined as follows: Let $\pi$ be some permutation of the class labels $\{1,\ldots,K_0\}$ and denote the set of all possible permutations by $\Pi$. Moreover, denote the group membership of index $i$ by $\rho(i)$, i.e.\ set $\rho(i) = k$ if $i \in G_k$. Similarly, let $\widehat{\rho}_{\pi}(i)$ be the estimated group membership of index $i$, where the estimated classes are labelled according to the permutation $\pi$. More specifically, set $\widehat{\rho}_{\pi}(i) = \pi(k)$ if $i \in \widehat{G}_k$. With this notation at hand, we define $\#F = \min_{\pi \in \Pi} \sum\nolimits_{i=1}^n \ind(\rho(i) \ne \widehat{\rho}_{\pi}(i))$.} The upper left panel of Figure \ref{fig2a-sim} shows the histogram of these $S=1000$ values for our multiscale approach. The other panels of Figure \ref{fig2a-sim} present the corresponding histograms for the bandwidth-dependent algorithm with the five different bandwidth values $h$ under consideration. As can be seen very clearly, our multiscale approach performs much better than the bandwidth-dependent competitor for any of the considered bandwidths. Figure \ref{fig2b-sim} shows the simulation results for the estimated number of classes $\widehat{K}_0$. The upper left panel depicts the histogram of the $S=1000$ values of $\widehat{K}_0$ produced by the multiscale approach. As one can see, the estimate $\widehat{K}_0$ equals the true number of classes $K_0 = 5$ in about 95\% of the cases (that is, in about $950$ out of $S=1000$ simulations). The performance of the bandwidth-dependent algorithm is considerably worse, which becomes apparent upon inspecting the other panels of Figure \ref{fig2b-sim}. The results in Figure \ref{fig3-sim} for the scenario with the positive AR parameter $a = 0.25$ give a very similar picture. In particular, our multiscale approach shows a much better performance than the bandwidth-dependent competitor for any of the considered bandwidths. Comparing Figures \ref{fig2-sim} and \ref{fig3-sim}, one can further see that the estimation precision is a bit better for the negative than the positive AR parameter (both for the multiscale and the bandwidth-dependent approach). This is not very surprising but simply reflects the fact that it is more difficult for the procedures to handle positive rather than negative correlation in the error terms.

Overall, our multiscale approach clearly outperforms the bandwidth-dependent algorithm in the simulation setup under consideration. Heuristically, this can be explained as follows: The setup comprises two very different types of signals. The signals $g_4$ and $g_5$ are very local in nature; they differ from a flat line only by a sharp, very local spike. The signals $g_2$ and $g_3$, in contrast, are much more global in nature; they differ from a flat line on a large part of the support $[0,1]$, but they are much smaller in magnitude than $g_4$ and $g_5$. A bandwidth-dependent clustering algorithm is hardly able to distinguish these signals reliably from each other. When a small bandwidth value is used, local features of the functions (the spikes in $g_4$ and $g_5$) can be detected reliably, but more global features (the slight curvature in $g_2$ and $g_3$) are hard to see. Hence, when implemented with a small bandwidth, the algorithm is barely able to detect the global differences between the functions. When implemented with a large bandwidth, in contrast, it is hardly able to capture the local differences. Our multiscale approach, in contrast, is able to produce appropriate estimates since it analyzes the data on various scales simultaneously.

Even though we have considered a quite stylized setup in our simulations, the advantages of our multiscale approach that become visible in this setup can be expected to persist in real-data applications. In practice, it is usually not known whether the group-specific regression functions $g_k$ $(1 \le k \le K_0)$ differ on a local or global scale. Hence, it is usually not clear at all which bandwidth is appropriate for implementing a bandwidth-dependent clustering algorithm. If the bandwidth is not picked suitably, the clustering results may not be very accurate. Moreover, when the functions $g_k$ differ on multiple scales, a clustering approach which is based on a single bandwidth $h$ can be expected to perform not very well, regardless of the specific value of $h$. Our multiscale approach, in contrast, can be expected to produce reliable clustering results, no matter whether the functions $g_k$ differ on a local, global or multiple scales.

\setcounter{equation}{0}
\section{Application}\label{sec-app}

In what follows, we revisit the application example from \cite{VogtLinton2017}. The aim of this example is to investigate the effect of trading venue fragmentation on market quality in the European stock market. For each stock $i$ in the FTSE 100 and FTSE 250 index, we observe a time series $\mathcal{T}_i = \{(Y_{it},X_{it}): 1 \le t \le T\}$ of weekly data from May 2008 to June 2011, where $Y_{it}$ is a measure of market quality and $X_{it}$ a measure of fragmentation for stock $i$ at time $t$. More specifically, $Y_{it}$ denotes the logarithmic volatility level of stock $i$ at time $t$, where volatility is measured by the so-called high-low range, which is defined as the difference between the highest and the lowest price of the stock at time $t$ divided by the latter. As a measure of fragmentation, we use the so-called Herfindahl index. The Herfindahl index of stock $i$ at time $t$ is defined as the sum of the squared market shares of the venues where the stock is traded at time $t$. It thus takes values between $0$ and $1$. If $X_{it}$ takes a value close to $0$, there is strong fragmentation in stock $i$ at time $t$, that is, stock $i$ is traded at many different venues at time $t$. A value of $X_{it}$ close to $1$, in contrast, indicates little fragmentation, that is, stock $i$ is traded only at a few venues at time $t$. The measures $Y_{it}$ and $X_{it}$ are constructed from data provided by Fidessa and Datastream. More details on the underlying data set and on variable construction can be found in \cite{LintonVogt2015} and \cite{KoerberLintonVogt2016}.

For each stock $i$, we model the relationship between $Y_{it}$ and $X_{it}$ by the nonparametric regression equation
\begin{equation}\label{model-app}
Y_{it} = m_i(X_{it}) + u_{it}, 
\end{equation}
where the error term has the fixed effects structure $u_{it} = \alpha_i + \gamma_t + \varepsilon_{it}$. The function $m_i$ captures the effect of trading-venue fragmentation on market quality for stock $i$. It is quite plausible to suppose that there are groups of stocks for which this effect is fairly similar. We thus impose a formal group structure on the stocks in our sample. In particular, we suppose that there are $K_0$ groups of stocks $G_1,\ldots,G_{K_0}$ such that $m_i = g_k$ for all $i \in G_k$ and all $1 \le k \le K_0$, where $g_k$ denotes the group-specific regression function associated with group $G_k$. Hence, we model the effect of fragmentation on market quality to be the same for all stocks in a given group.

\begin{figure}[p]
\includegraphics[width=0.925\textwidth]{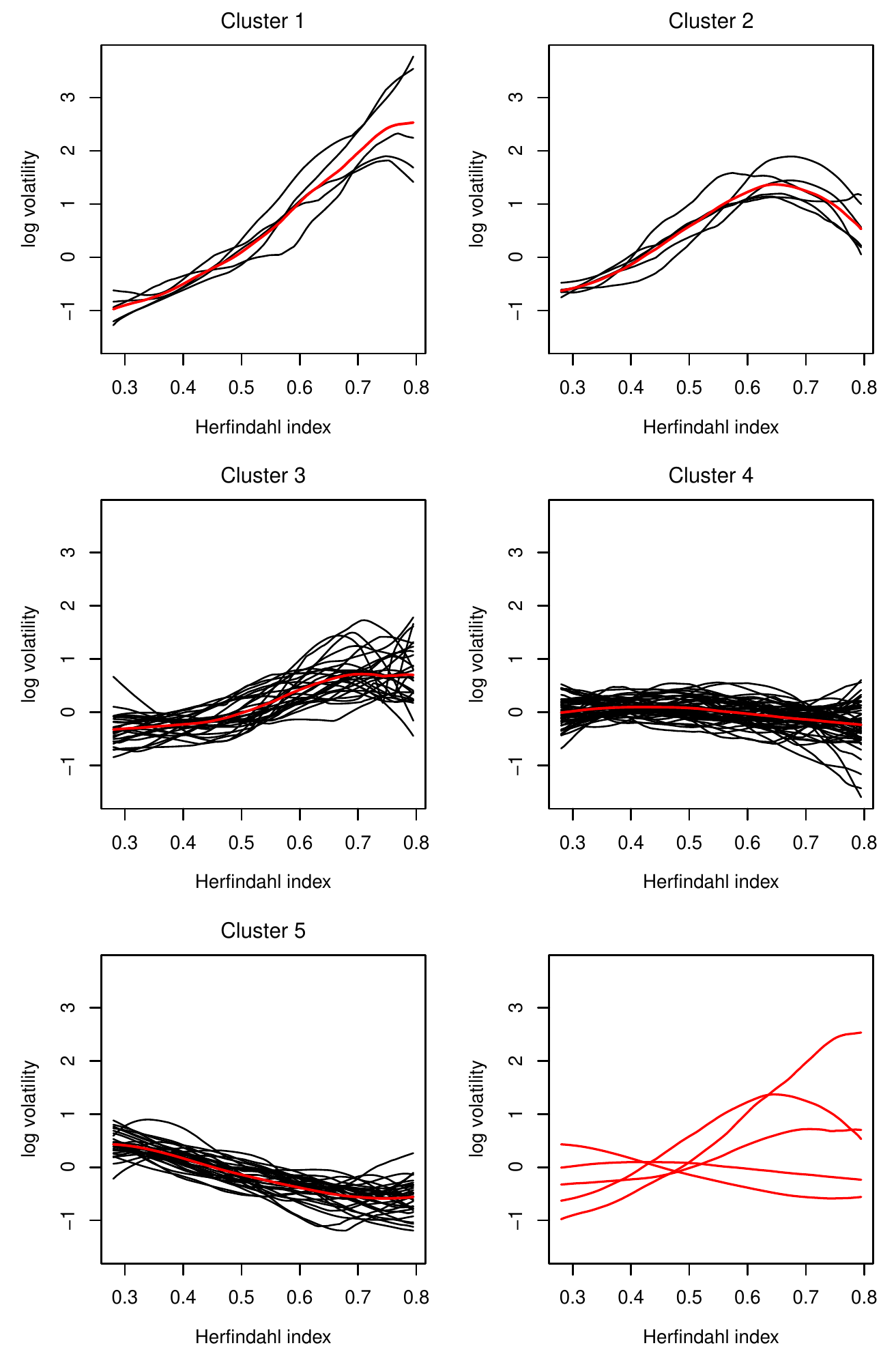} 
\caption{Estimated clusters in the application example of Section \ref{sec-app}. Each panel corresponds to one cluster. The black lines are the estimated regression curves $\widehat{m}_{i,h}$ that belong to the respective cluster. The red lines are estimates of the group-specific regression functions. These are plotted once again together in the lower right panel of the figure.}\label{fig-app}
\end{figure}

We now use our multiscale clustering methods to estimate the unknown groups $G_1,\ldots,G_{K_0}$ along with their unknown number $K_0$ from the data sample at hand. As in \cite{VogtLinton2017}, we drop stocks from the sample for which data points are missing. Moreover, we eliminate stocks $i$ with a very small empirical support $\mathcal{S}_i$ of the fragmentation data $\{ X_{it}: 1 \le t \le T \}$. In particular, we only take into account stocks $i$ for which the support $\mathcal{S}_i$ contains the interval $[0.275,0.8]$. We thus use exactly the same data set as in \cite{VogtLinton2017}, which comprises $n=125$ time series of length $T=151$ weeks. To implement our multiscale methods, we employ the location-scale grid $\mathcal{G}_T = \{ (x,h) : x \in \mathcal{X} \text{ and } h \in \mathcal{H} \}$, where $\mathcal{X} = \{ x : x = r/100 \text{ for } r = 1,\ldots,99 \}$ is the set of locations and $\mathcal{H} = \{ h : 0.175 \le h \le 0.5 \text{ with } h = 0.025 k \text{ for } k = 1,2,\ldots \}$ is the set of bandwidths. Note that $h = 0.175$ is the smallest possible bandwidth we can use: If we pick $h$ smaller than $0.175$, we cannot compute the statistics $\widehat{\psi}_{ij}(x,h)$ for all stocks $i$ and locations $x \in \mathcal{X}$ any more because for some $i$ and $x$, there are less than two data points in the bandwidth window. The threshold parameter for the estimation of $K_0$ is set to $\pi_{n,T} = q_n(\alpha)$ with $\alpha = 0.95$.

\enlargethispage{0.2cm}
The estimation results are presented in Figure \ref{fig-app}. Each panel of the figure corresponds to one of the estimated groups $\widehat{G}_k^{[\widehat{K}_0]}$ for $1 \le k \le \widehat{K}_0$, where the estimated number of groups is $\widehat{K}_0 = 5$. In particular, each panel depicts the estimated curves $\widehat{m}_{i,h}$ that belong to some cluster $\widehat{G}_k^{[\widehat{K}_0]}$. The red curve in each panel is an estimate $\widehat{g}_{k,h}$ of the group-specific regression function $g_k$. More specifically, we define 
\[ \widehat{g}_{k,h}(x) = \frac{1}{\widehat{G}_k^{[\widehat{K}_0]}} \sum\limits_{i \in \widehat{G}_k^{[\widehat{K}_0]}} \widehat{m}_{i,h}(x), \]
that is, we simply average the fits $\widehat{m}_{i,h}$ with $i \in \widehat{G}_k^{[\widehat{K}_0]}$. The estimates $\widehat{g}_{k,h}(x)$ are once again plotted together in the lower right panel of Figure \ref{fig-app}. Whereas we do not need to specify a bandwidth $h$ to compute the multiscale estimates $\widehat{G}_k^{[\widehat{K}_0]}$ of the unknown groups for $1 \le k \le \widehat{K}_0$, the kernel smoothers $\widehat{m}_{i,h}$ of course depend on a specific bandwidth $h$. As these smoothers are only computed for illustrative purposes, in particular for the graphical illustration of the results in Figure \ref{fig-app}, we use the same bandwidth $h$ for all stocks $i$. In particular, we choose the bandwidth adhoc as $h = 0.25$ for all $i$, which produces a good visual impression of the results.

In order to interpret the results in Figure \ref{fig-app}, we regard volatility as a bad, meaning that higher volatility implies lower market quality. As can be seen, the effect of fragmentation on volatility is quite moderate for most stocks: Most of the curve fits in Cluster 4 are close to a flat line, whereas those in Clusters 3 and 5 slightly slope upwards and downwards, respectively. In contrast to this, the fits in Cluster 1 and to a lesser extent also those in Cluster 2 exhibit a strong increase. This indicates that higher fragmentation is accompanied by lower volatility and thus higher market quality for these stocks. To summarize, fragmentation appears to substantially improve market quality only for a small share of stocks (in particular for those in Clusters 1 and 2), whereas the effect of fragmentation is quite moderate for the great bulk of stocks (in particular for those in Clusters 3, 4 and 5). These findings are in line with those in \cite{VogtLinton2017}. Indeed, the clusters produced by our multiscale method are fairly similar to those obtained there. Hence, our multiscale approach confirms the results of the bandwidth-dependent algorithm from \cite{VogtLinton2017}, but without the need to go through the complicated bandwidth-selection procedure from there which may very well perform less accurate in other applications.

\setcounter{equation}{0} 
\section{Extensions and modifications}\label{sec-ext}

\subsection{Extension to the multivariate case and to other model settings}

Throughout the paper, we have restricted attention to real-valued regressors $X_{it}$. Our approach extends to $\reals^d$-valued regressors $X_{it} = (X_{it,1},\ldots,X_{it,d})^\top$ in a straightforward way. The clustering methods described in Sections \ref{sec-est-classes} and \ref{sec-est-K0} remain the same in the multivariate case, only the multiscale statistics $\widehat{d}_{ij}$ need to be adjusted. To do so, we simply need to (i) replace the involved kernel estimators by multivariate versions and (ii) modify the scaling factors $\widehat{\nu}_{ij}(x,h)$ appropriately to normalize the variance of the statistics $\widehat{\psi}_{ij}(x,h)$. We neglect the details as these modifications are very straightforward.

The kernel smoothers on which the multiscale statistics $\widehat{d}_{ij}$ are based suffer from the usual curse of dimensionality. Hence, our fully nonparametric approach is only useful in practice as long as the dimension $d$ of the regressors is moderate. If $d$ is large, it makes sense to resort to structured nonparametric or semiparametric approaches. As an example, consider the partially linear model 
\begin{equation}\label{model-partiallinear}
Y_{it} = m_i(X_{it}) + \boldsymbol{\beta}^\top Z_{it} + u_{it},
\end{equation}
where $X_{it}$ is real-valued, $Z_{it} = (Z_{it,1},\ldots,Z_{it,d})^\top$ is an $\reals^d$-valued vector and the error terms $u_{it}$ have the fixed effects structure $u_{it} = \alpha_i + \gamma_t + \varepsilon_{it}$ with $\ex[\varepsilon_{it}|X_{it},Z_{it}] = 0$. In this model, $Z_{it}$ is a vector of controls which enters the equation \eqref{model-partiallinear} linearly for simplicity. In particular, $\boldsymbol{\beta} = (\beta_{1},\ldots,\beta_{d})^\top$ is an unknown parameter vector which is assumed to be the same for all $i$. Suppose we are mainly interested in the effect of $X_{it}$ on the response $Y_{it}$, which is captured by the functions $m_i$. As in Section \ref{sec-model}, we may model this effect by imposing a group structure on the curves $m_i$: We may suppose that there exist classes $G_1,\ldots,G_{K_0}$ and associated functions $g_1,\ldots,g_{K_0}$ such that $m_i = g_k$ for all $i \in G_k$ and $1 \le k \le {K_0}$. In order to apply our estimation methods in this context, we merely need to adjust the multiscale statistics $\widehat{d}_{ij}$. In particular, we need to replace the local linear smoothers $\widehat{m}_{i,h}(x)$ by appropriate estimators of $m_i$ and adjust the scaling factors $\widehat{\nu}_{ij}(x,h)$. The functions $m_i$ may for example be estimated with the help of the methods developed in \cite{Robinson1988}. Once the multiscale statistics $\widehat{d}_{ij}$ have been adjusted to the partially linear model setting \eqref{model-partiallinear}, estimators of the unknown classes and their unknown number can be obtained as described in Sections \ref{sec-est-classes} and \ref{sec-est-K0}. We conjecture that the two main Theorems \ref{theo-alg} and \ref{theo-thres} on the multiscale clustering methods remain to hold true in the context of the partially linear model \eqref{model-partiallinear}. However, extending our theoretical results to model \eqref{model-partiallinear} is by no means trivial but would require a substantial deal of additional work.

Another interesting model setting to which our methods can be extended is the following: Suppose that 
\begin{equation}\label{model-ext-treatment-1}
Y_{it}=m_i(X_{it})+ \alpha_i + \gamma_t + \varepsilon_{it},
\end{equation}
where $X_{it}$ is a continuous treatment effect that is applied to unit $i$ during periods $\tau_{i} \subset \{1,\ldots ,T\}$ and is not present during pre- and post-treatment periods. Moreover, suppose that there is a matched control group that never receives the treatment and so satisfies  
\begin{equation}\label{model-ext-treatment-2}
Y_{j(i),t}=\alpha_{j(i)}+\gamma_{t}+\varepsilon _{j(i),t},
\end{equation}
where $j(i)$ denotes the unit in the control group which is matched with $i$. This setting is similar to that considered in \cite{Bonevaetal2018} who evaluate the effects of the UK government's corporate bond purchase scheme on market quality measures such as liquidity. The treatment in this study is continuously distributed and applied during an $18$-month-period to a subset of all UK listed corporate bonds. The authors assume a linear homogeneous treatment effect in the baseline model and apply difference-in-difference methods to estimate the effect. However, one could easily allow for more general nonlinear and heterogeneous effects $m_i$ as in equation \eqref{model-ext-treatment-1} and impose a group structure on them. Notice that for $t \in \tau_{i}$ and $s \in \tau_{i}^c$, we have  
\begin{equation}\label{model-ext-treatment-dd}
(Y_{it} - Y_{j(i),t}) -(Y_{is} - Y_{j(i),s}) = m_i(X_{it}) + \varepsilon_{it} - \varepsilon_{j(i),t} - \varepsilon_{is} + \varepsilon_{j(i),s}, 
\end{equation}
which is essentially a nonparametric regression equation for each $i$. We could thus apply our methods to the difference-in-difference equation \eqref{model-ext-treatment-dd}.

\subsection{Alternatives to hierarchical clustering}

 In order to estimate the unknown class structure in model \eqref{model-eq1}--\eqref{model-eq2}, we have combined the multiscale statistics $\widehat{d}_{ij}$ with a hierarchical clustering algorithm. It is also possible to combine them with other distance-based clustering approaches. In parti\-cular, they can be employed as distance statistics in the thresholding algorithm of \cite{VogtLinton2017}. To do so, we replace the $L_2$-type distance statistics $\widehat{\Delta}_{ij}$ from \cite{VogtLinton2017} by the multiscale statistics $\widehat{d}_{ij}$ and construct the threshold estimators of the unknown groups $G_1,\ldots,G_{K_0}$ and of their unknown number $K_0$ exactly as described in Section 2.2 of \cite{VogtLinton2017}. This leads to estimators $\widetilde{K}_0$ and $\widetilde{G}_1,\ldots,\widetilde{G}_{\widetilde{K}_0}$, which unlike those constructed in \cite{VogtLinton2017} are free of classical bandwidth parameters. 

Under regularity conditions very similar to those from Section \ref{sec-theory}, we can derive some basic theoretical properties of the estimators $\widetilde{K}_0$ and $\widetilde{G}_1,\ldots,\widetilde{G}_{\widetilde{K}_0}$: Suppose that the threshold parameter $\tau_{n,T}$ of the procedure fulfills Condition 6 from Section 3.2 of \cite{VogtLinton2017}, that is, $\tau_{n,T} \searrow 0$ such that $\max_{i,j \in G_k} \widehat{d}_{ij} \le \tau_{n,T}$ with probability tending to $1$ for all $k$. Then it can be shown that $\pr(\widetilde{K}_0 = K_0) \rightarrow 1$ as well as $\pr( \{ \widetilde{G}_1,\ldots,\widetilde{G}_{\widetilde{K}_0} \} = \{G_1,\ldots,G_{K_0}\} ) \rightarrow 1$. 

To implement the estimators $\widetilde{K}_0$ and $\widetilde{G}_1,\ldots,\widetilde{G}_{\widetilde{K}_0}$ in practice, we need to choose the threshold level $\tau_{n,T}$. In view of Condition 6 from \cite{VogtLinton2017}, we would like to tune $\tau_{n,T}$ such that $\max_{i,j \in G_k} \widehat{d}_{ij} \le \tau_{n,T}$ holds with high probability for all $k$. According to our heuristic arguments from Section \ref{subsec-choice-thres}, this may be achieved by setting $\tau_{n,T} = q_n(\alpha)$ with $\alpha$ close to $1$. We thus suggest to choose the threshold parameter $\tau_{n,T}$ in the same way as the dissimilarity level $\thres_{n,T}$ at which we cut the dendrogram to estimate $K_0$.

\subsection{Letting $\boldsymbol{K_0}$ grow with the sample size}

Throughout the paper, we have assumed that the number of classes $K_0$ is fixed. We now allow $K_0$ to grow with the number of time series $n$, that is, we admit of $K_0 = K_{0,n} \rightarrow \infty$ as $n \rightarrow \infty$. To deal with this situation, we require the group-specific regression functions $g_k$ to fulfill the following additional condition:
\begin{enumerate}[label=(C\arabic*),leftmargin=1.25cm]
\setcounter{enumi}{9}
\item \label{C-add} The functions $g_k$ as well as their first and second derivatives are uniformly bounded in absolute value, that is, $|g_k^{(\ell)}(x)| \le C$ for all $x \in [0,1]$ and $\ell = 0,1,2$, where $g_k^{(\ell)}$ denotes the $\ell$-th derivative of $g_k$ and the constant $C < \infty$ does not depend on $k$. Moreover, 
\begin{equation}\label{C-add-eq}
\min_{1 \le k < k^{\prime} \le K_0} \max_{\{x \, : \, (x,h_{\max}) \in \mathcal{G}_T\}} |g_k(x) - g_{k^\prime}(x)| \gg \frac{\sqrt{\log n + \log T} + \sqrt{T h_{\max}^5}}{\sqrt{Th_{\max}}}. 
\end{equation}
\end{enumerate}
As before, the expression $a_{n,T} \gg b_{n,T}$ means that $b_{n,T} = o(a_{n,T})$ and the notation $a_{n,T} \ll b_{n,T}$ is used analogously. \eqref{C-add-eq} essentially says that the regression functions $g_k$ and $g_{k^\prime}$ of two different classes do not approach each other too quickly as $n \rightarrow \infty$. If condition \ref{C-add} is fulfilled, a slightly modified version of Theorem \ref{theo-dist} can be proven. In particular, with the help of the technical arguments from the Supplementary Material, it is not difficult to show that 
\begin{align*}
\max_{1 \le k \le K_0} \max_{i, j \in G_k} \, \widehat{d}_{ij} & = O_p\big( \sqrt{\log n + \log T} \big) \\
\min_{1 \le k < k^\prime \le K_0} \min_{\substack{ i \in G_k, \\ j \in G_{k^\prime} }} \widehat{d}_{ij} & \gg \sqrt{\log n + \log T} + \sqrt{T h_{\max}^5}.
\end{align*}
These two statements immediately imply that Theorem \ref{theo-alg} remains to hold true. Moreover, Theorem \ref{theo-thres} remains valid as well if the threshold level $\pi_{n,T}$ satisfies a strengthened version of condition \eqref{prop-thres}, namely the condition that $\sqrt{\log n + \log T} \ll \thres_{n,T} \ll \sqrt{Th_{\max}} \, \min_{1 \le k < k^{\prime} \le K_0}  \max_{\{x \, : \, (x,h_{\max}) \in \mathcal{G}_T\}} |g_k(x) - g_{k^\prime}(x)|$.

\bibliographystyle{ims}
{\small
\setlength{\bibsep}{0.45em}
\bibliography{bibliography}}

\newpage
\input{supplement}

\end{document}

%% file: macros.tex
\newcommand{\reals}{\mathbb{R}}
\newcommand{\integers}{\mathbb{Z}}

\newcommand{\pr}{\mathbb{P}}        
\newcommand{\ex}{\mathbb{E}}        
\newcommand{\var}{\textnormal{Var}} 
\newcommand{\cov}{\textnormal{Cov}} 

\newcommand{\normal}{N}        

\newcommand{\ind}{1} 
\newcommand{\kernel}{K} 
\newcommand{\wght}{W} 
\newcommand{\thres}{\pi} 

\newcommand{\convd}{\stackrel{d}{\longrightarrow}}              
\newcommand{\convp}{\stackrel{P}{\longrightarrow}}              

\theoremstyle{plain}

\newtheorem{theorem}{Theorem}[section]
\newtheorem{prop}{Proposition}[section]

\newtheorem{propS}{Proposition}
\newtheorem{lemmaS}[propS]{Lemma}

\newtheorem{remark}{Remark}[section]

\makeatletter
\newcommand{\lefteqno}{\let\veqno\@@leqno}
\makeatother

\newcommand{\heading}[2]
{  \setcounter{page}{1}
   \begin{center}

   \phantom{Distance to upper boundary}
   \vspace{0.5cm}

   {\LARGE \textbf{#1}}
   \vspace{0.4cm}
 
   {\LARGE \textbf{#2}}
   \end{center}
}

\newcommand{\authors}[4]
{  \parindent0pt
   \begin{center}
      \begin{minipage}[c][2cm][c]{5cm}
      \begin{center} 
      {\large #1} 
      \vspace{0.1cm}
      
      #2 
      \end{center}
      \end{minipage}
      \begin{minipage}[c][2cm][c]{5cm}
      \begin{center} 
      {\large #3}
      \vspace{0.1cm}

      #4 
      \end{center}
      \end{minipage}
   \end{center}
}



%% file: supplement.tex
\setcounter{page}{1}
\renewcommand{\baselinestretch}{1.0}\normalsize

\begin{center}
{\LARGE \textbf{Supplement to}} \\[0.3cm]
{\LARGE \textbf{``Multiscale Clustering of}} \\[0.3cm]
{\LARGE \textbf{Nonparametric Regression Curves''}}
\end{center}
\vspace{-0.6cm}

\authors{Michael Vogt}{University of Bonn}{Oliver Linton}{University of Cambridge}
\vspace{-0.8cm}

\renewcommand{\abstractname}{}
\begin{abstract}
\noindent In this supplement, we provide the technical details omitted in the paper. In Section \ref{secS-identification}, we prove Proposition \ref{prop-identify} which concerns identification of the functions $m_i$. Sections \ref{secS-uniform1} and \ref{secS-uniform2} contain some auxiliary results needed for the proof of Theorem \ref{theo-dist}. In Section \ref{secS-uniform1}, we in particular derive a general uniform convergence result which is applied to the kernel smoothers $\widehat{m}_{i,h}$ in Section \ref{secS-uniform2}. The final Section \ref{secS-theo-dist} contains the proof of Theorem \ref{theo-dist}. 
Throughout the supplement, we use the following notation: The symbol $C$ denotes a universal real constant which may take a different value on each occurrence. In addition, the symbols $C_0, C_1, \ldots$ are used to denote specific real constants that are defined in the course of the supplement. Unless stated differently, the constants $C, C_0, C_1,\ldots$ depend neither on the dimensions $n$ and $T$, nor on the indices $i \in \{1,\ldots,n\}$ and $t \in \{1,\ldots,T\}$, nor on the location-bandwidth points $(x,h) \in \mathcal{G}_T$. To emphasize that the constants $C, C_0, C_1,\ldots$ do not depend on any of these parameters, we refer to them as absolute constants in many places.
\end{abstract}

\renewcommand{\baselinestretch}{1.2}\normalsize
\setlength{\parindent}{0.75cm} 
\renewcommand{\theequation}{S.\arabic{equation}}
\setcounter{equation}{0}
\def\thesection{S.\arabic{section}}
\setcounter{section}{0}
\allowdisplaybreaks[4]

\section{Proof of Proposition \ref{prop-identify}}\label{secS-identification}

Let $\overline{Y}_i$, $\overline{Y}_t^{(i)}$ and $\overline{\overline{Y}}^{(i)}$ be the sample averages introduced in \eqref{sample-av-eq}, that is, 
\[ \overline{Y}_i = \frac{1}{T} \sum_{t=1}^{T} Y_{it}, \quad \overline{Y}_t^{(i)} = \frac{1}{n-1} \sum_{\substack{j=1 \\ j \ne i}}^{n} Y_{jt} \quad \text{and} \quad \overline{\overline{Y}}^{(i)} = \frac{1}{(n-1)T} \sum_{\substack{j=1 \\ j \ne i}}^{n} \sum_{t=1}^T Y_{jt}. \]
Define $\overline{\varepsilon}_i$, $\overline{\varepsilon}_t^{(i)}$ and $\overline{\overline{\varepsilon}}^{(i)}$ analogously and set $\overline{m}_i = T^{-1} \sum_{t=1}^{T}m_i(X_{it})$, $\overline{m}_t^{(i)} = (n-1)^{-1}$ $\sum_{j = 1, j \ne i}^n m_j(X_{jt})$ and $\overline{\overline{m}}^{(i)} = (\{n-1\}T)^{-1} \sum_{j=1, j \ne i}^{n} \sum_{t=1}^T m_j(X_{jt})$. Straightforward calculations yield that 
\begin{align}
Y_{it} - \overline{Y}_i - \overline{Y}_t^{(i)} + \overline{\overline{Y}}^{(i)} 
   = m_{i}(X_{it}) & - \overline{m}_i - \overline{m}_t^{(i)} + \overline{\overline{m}}^{(i)} \nonumber \\
 & + \varepsilon_{it} - \overline{\varepsilon}_i - \overline{\varepsilon}_t^{(i)} + \overline{\overline{\varepsilon}}^{(i)}. \label{within-eq}
\end{align}  
Hence, by adding/subtracting the sample averages $\overline{Y}_i$, $\overline{Y}_t^{(i)}$ and $\overline{\overline{Y}}^{(i)}$ from $Y_{it}$, we can eliminate the fixed effects $\alpha_i$ and $\gamma_t$ from the model equation \eqref{model-eq}. We now consider the transformed model equation \eqref{within-eq} for arbitrary but fixed indices $i$ and $t$ and examine the following two cases separately: (a) $n = n(T) \rightarrow \infty$ as $T \rightarrow \infty$, and (b) $n = n(T)$ remains bounded as $T \rightarrow \infty$. 
\begin{enumerate}[label=(\alph*),leftmargin=0.75cm]

\item Under the normalization constraint \eqref{C-identify} and the assumptions of Proposition \ref{prop-identify}, it holds that for any fixed $i$ and $t$, $\overline{\varepsilon}_i = O_p(T^{-1/2})$ and $\overline{m}_i = O_p(T^{-1/2})$, $\overline{\varepsilon}_t^{(i)} = O_p(n^{-1/2})$ and $\overline{m}_t^{(i)} = O_p(n^{-1/2})$ as well as $\overline{\overline{\varepsilon}}^{(i)} = O_p(\{nT\}^{-1/2})$ and $\overline{\overline{m}}^{(i)} = O_p(\{nT\}^{-1/2})$. Using these facts in equation \eqref{within-eq} for a fixed pair of indices $i$ and $t$, we obtain that 
\begin{equation}\label{within-eq-1}
Y_{it}^{\infty} = m_i(X_{it}) + \varepsilon_{it} \quad \text{a.s.}, 
\end{equation}
where $Y_{it}^{\infty}$ denotes the limit of $\widehat{Y}_{it}^* = Y_{it} - \overline{Y}_i - \overline{Y}_t^{(i)} + \overline{\overline{Y}}^{(i)}$ in probability, that is, $\widehat{Y}_{it}^* \convp Y_{it}^{\infty}$. From \eqref{within-eq-1}, it follows that $\ex[ Y_{it}^{\infty} | X_{it}] = m_i(X_{it})$ almost surely, which identifies $m_i$. 

\item Now suppose that $n = n(T)$ remains bounded as $T \rightarrow \infty$. Let us assume for simplicity that $n = n(T)$ is non-decreasing in $T$, implying that $n$ is a fixed number for sufficiently large $T$. (Without this assumption, we would have to consider a subsequence of time series lengths $T_k$ for $k = 1,2,\ldots$ such that $n(T_k)$ is non-decreasing.) Similar to the previous case, we have that $\overline{\varepsilon}_i = O_p(T^{-1/2})$ and $\overline{m}_i = O_p(T^{-1/2})$ as well as $\overline{\overline{\varepsilon}}^{(i)} = O_p(T^{-1/2})$ and $\overline{\overline{m}}^{(i)} = O_p(T^{-1/2})$. Using these facts in equation \eqref{within-eq}, we arrive at 
\begin{equation}\label{within-eq-2}
Y_{it}^{\infty} = m_i(X_{it}) + \varepsilon_{it} - \big\{ \overline{m}_t^{(i)} + \overline{\varepsilon}_t^{(i)} \big\}  \quad \text{a.s.},
\end{equation}
where $Y_{it}^{\infty}$ is defined as before and, slightly abusing notation, we let $\overline{\varepsilon}_t^{(i)} = (N-1)^{-1} \sum_{j = 1, j \ne i}^N \varepsilon_{jt}$ and $\overline{m}_t^{(i)} = (N-1)^{-1} \sum_{j = 1, j \ne i}^N m_j(X_{jt})$ with $N = \lim_{T \rightarrow \infty} n(T)$. Since $\ex[\overline{\varepsilon}_t^{(i)}|X_{it}] = \ex[\overline{\varepsilon}_t^{(i)}] = 0$ and $\ex[\overline{m}_t^{(i)}|X_{it}] = \ex[\overline{m}_t^{(i)}] = 0$ under the normalization constraint \eqref{C-identify} and the assumptions of Proposition \ref{prop-identify}, we get that $\ex[ Y_{it}^{\infty} | X_{it}] = m_i(X_{it})$ almost surely, which once again identifies $m_i$. 

\end{enumerate}

\section{A general result on uniform convergence}\label{secS-uniform1}

In this and the subsequent section, we derive some uniform convergence results needed for the proof of Theorem \ref{theo-dist}. The multiscale statistics $\widehat{d}_{ij}$ are composed of kernel estimators whose building blocks are kernel averages of the form 
\begin{equation}\label{kernel-av}
\Phi_i(x,h) = \frac{1}{T} \sum\limits_{t=1}^{T} \kernel_h(X_{it} - x) \Big( \frac{X_{it} - x}{h} \Big)^\ell Z_{it,T},
\end{equation}
where $\ell$ is a fixed natural number and $X_{it}$ are the regressor variables from model \eqref{model-eq1}. Moreover, $Z_{it,T}$ are general real-valued random variables that may depend on the sample size parameter $T$. For each $i$, the variables $(Z_{it,T},X_{it})$ form a triangular array $\mathcal{A}_i = \{ \mathcal{A}_{i,T} \}_{T=1}^\infty$, where $\mathcal{A}_{i,T} = \{ (Z_{it,T},X_{it}): 1 \le t \le T \}$. We make the following assumptions on the random variables $(Z_{it,T},X_{it})$:
\begin{enumerate}[label=(P\arabic*),leftmargin=1cm]
\item \label{P1} For each $i$ and $T$, the collection of random variables $\mathcal{A}_{i,T}$ is strongly mixing. The mixing coefficients $\alpha_{i,T}(\ell)$ of $\mathcal{A}_{i,T}$ are such that $\alpha_{i,T}(\ell) \le n \, \alpha(\ell)$ for all $i$, $T$ and $\ell$, where the coefficients $\alpha(\ell)$ decay exponentially fast to zero as $\ell \rightarrow \infty$.  
\item \label{P2} There exist a real number $\theta > 2$ and a natural number $\ell^*$ such that for any $\ell \in \integers$ with $|\ell| \ge \ell^*$ and some absolute constant $C < \infty$, 
\begin{align*}
 & \max_{1 \le t \le T} \max_{1 \le i \le n} \sup_{x \in [0,1]} \ex \big[ |Z_{it,T}|^\theta \big| X_{it} = x \big] \le C < \infty \\
 & \max_{1 \le t \le T} \max_{1 \le i \le n} \sup_{x,x^\prime \in [0,1]} \ex \big[ |Z_{it,T} Z_{it+\ell,T}| \big| X_{it} = x, X_{it+\ell} = x^\prime \big] \le C < \infty.
\end{align*} 
\end{enumerate}
The following lemma characterizes the convergence behaviour of the kernel average $\Phi_i(x,h)$ uniformly over $i$, $x$ and $h$. 
\begin{propS}\label{propS1} 
Let \ref{P1} and \ref{P2} be satisfied. Moreover, assume that \ref{C-dens} and \ref{C-nT}--\ref{C-ker} are fulfilled. Then it holds that
\[ \pr \Big( \max_{1 \le i \le n} \max_{(x,h) \in \mathcal{G}_T} \sqrt{Th} \big| \Phi_i(x,h) - \ex \Phi_i(x,h) \big| > C_0 \sqrt{\gamma_{n,T}} \Big) = o(1), \]
where $\gamma_{n,T} = \log n + \log T$ and $C_0$ is a sufficiently large absolute constant. 
\end{propS}

\noindent \textbf{Proof of Proposition \ref{propS1}.} 
To prove the proposition, we modify standard arguments to derive uniform convergence rates for kernel estimators, which can be found e.g.\ in \cite{Masry1996}, \cite{Bosq1998} or \cite{Hansen2008}. These arguments were originally designed to derive the convergence rates of kernel averages such as $\Phi_i(x,h) - \ex \Phi_i(x,h)$ uniformly over $x$ but pointwise in $h$ and $i$. In contrast to this, we aim to derive the convergence rate of $\Phi_i(x,h) - \ex \Phi_i(x,h)$ uniformly over $x$, $h$ and $i$. Related results can be found e.g.\ in \cite{Einmahl2005} and \cite{VogtLinton2017} (see in particular Lemma S.1 therein).

We now turn to the proof of the proposition. For simplicity of notation, we let $\ell = 0$ in \eqref{kernel-av}, the arguments being completely analogous for $\ell \ne 0$. To start with, we define
\begin{align*}
Z_{it,T}^{\le} & = Z_{it,T} \, \ind\big(|Z_{it,T}| \le (nT)^{\frac{1}{\theta-\delta}} \big) \\
Z_{it,T}^{>} & = Z_{it,T} \, \ind\big(|Z_{it,T}| > (nT)^{\frac{1}{\theta-\delta}} \big), 
\end{align*}
where $\delta > 0$ is an absolute constant that can be chosen as small as desired. Moreover, we write
\[ \sqrt{Th} \big\{ \Phi_i(x,h) - \ex \Phi_i(x,h) \big\} = \sum\limits_{t=1}^{T} \mathcal{Z}_{it,T}^{\le}(x,h) + \sum\limits_{t=1}^{T} \mathcal{Z}_{it,T}^{>}(x,h), \] 
where 
\begin{align*}
\mathcal{Z}_{it,T}^{\le}(x,h) & = \frac{1}{\sqrt{Th}} \Big\{ \kernel \Big(\frac{X_{it}-x}{h}\Big) Z_{it,T}^{\le} - \ex \Big[ \kernel \Big(\frac{X_{it}-x}{h}\Big) Z_{it,T}^{\le} \Big] \Big\} \\
\mathcal{Z}_{it,T}^{>}(x,h) & = \frac{1}{\sqrt{Th}} \Big\{ \kernel \Big(\frac{X_{it}-x}{h}\Big) Z_{it,T}^{>} - \ex \Big[ \kernel \Big(\frac{X_{it}-x}{h}\Big) Z_{it,T}^{>} \Big] \Big\}.
\end{align*}
With this notation at hand, we get that 
\[ \pr \Big( \max_{1 \le i \le n} \max_{(x,h) \in \mathcal{G}_T} \sqrt{Th} \big| \Phi_i(x,h) - \ex \Phi_i(x,h) \big| > C_0 \sqrt{\gamma_{n,T}} \Big) \le P^{\le} + P^{>}, \]
where
\begin{align*}
P^{\le} & = \pr \Big( \max_{1 \le i \le n} \max_{(x,h) \in \mathcal{G}_T} \Big| \sum\limits_{t=1}^{T} \mathcal{Z}_{it,T}^{\le}(x,h) \Big| > \frac{C_0}{2} \sqrt{\gamma_{n,T}} \Big) \\ 
P^{>} & = \pr \Big( \max_{1 \le i \le n} \max_{(x,h) \in \mathcal{G}_T} \Big| \sum\limits_{t=1}^{T} \mathcal{Z}_{it,T}^{>}(x,h) \Big| > \frac{C_0}{2} \sqrt{\gamma_{n,T}} \Big). 
\end{align*}
In what follows, we show that $P^{\le} = o(1)$ and $P^{>} = o(1)$, which implies the statement of Proposition \ref{propS1}.

We first have a closer look at $P^{>}$. It holds that 
\[ P^{>} \le \sum\limits_{i=1}^n \pr \Big( \max_{(x,h) \in \mathcal{G}_T} \Big| \sum\limits_{t=1}^{T} \mathcal{Z}_{it,T}^{>}(x,h) \Big| > \frac{C_0}{2} \sqrt{\gamma_{n,T}} \Big) \le P^{>}_1 + P^{>}_2, \]
where
\begin{align*}
P_1^{>} & = \sum\limits_{i=1}^n \pr \Big( \max_{(x,h) \in \mathcal{G}_T} \Big| \frac{1}{\sqrt{Th}} \sum\limits_{t=1}^{T} \kernel \Big(\frac{X_{it}-x}{h}\Big) Z_{it,T}^{>} \Big| > \frac{C_0}{4} \sqrt{\gamma_{n,T}} \Big) \\
P_2^{>} & = \sum\limits_{i=1}^n \pr \Big( \max_{(x,h) \in \mathcal{G}_T} \Big| \frac{1}{\sqrt{Th}} \sum\limits_{t=1}^{T} \ex \Big[ \kernel \Big(\frac{X_{it}-x}{h}\Big) Z_{it,T}^{>} \Big] \Big| > \frac{C_0}{4} \sqrt{\gamma_{n,T}} \Big). 
\end{align*} 
With the help of \ref{P2}, we obtain that
\begin{align*} 
P_1^{>} & \le \sum\limits_{i=1}^n \pr \Big( |Z_{it,T}| > (nT)^{\frac{1}{\theta-\delta}} \text{ for some } 1 \le t \le T \Big) \\*
 & \le \sum\limits_{i=1}^n \sum\limits_{t=1}^{T} \pr \Big( |Z_{it,T}| > (nT)^{\frac{1}{\theta-\delta}} \Big) \\*[0.2cm] & \le C (nT) \big/ (nT)^{\frac{\theta}{\theta-\delta}}  \\*[0.2cm] & = o(1). 
\end{align*}
Once again exploiting \ref{P2}, we can further infer that   
\begin{align*} 
\Big| \frac{1}{\sqrt{Th}} \sum\limits_{t=1}^{T} \ex \Big[ \kernel \Big(\frac{X_{it}-x}{h}\Big) Z_{it,T}^{>} \Big] \Big|
 & \le \frac{1}{\sqrt{Th}} \sum\limits_{t=1}^{T} \ex \Big[ \kernel \Big(\frac{X_{it}-x}{h}\Big) \frac{|Z_{it,T}|^{\theta}}{(nT)^{\frac{\theta-1}{\theta-\delta}}} \Big] \\ & \le C \sqrt{Th} \big/ (nT)^{\frac{\theta-1}{\theta-\delta}} \\
 & = o\big(\sqrt{\gamma_{n,T}}\big),
\end{align*}
which immediately implies that $P_2^> = 0$ for sufficiently large $T$. Putting everything together, we arrive at the result that $P^{>} = o(1)$.

We now turn to the analysis of $P^{\le}$. In what follows, we show that 
\begin{equation}\label{res-exp-bound} 
\max_{1 \le i \le n} \max_{(x,h) \in \mathcal{G}_T} \pr \Big( \Big| \sum\limits_{t=1}^{T} \mathcal{Z}_{it,T}^{\le}(x,h) \Big| > \frac{C_0}{2} \sqrt{\gamma_{n,T}} \Big) \le C T^{-r}, 
\end{equation}
where the constant $r > 0$ can be chosen as large as desired. From \eqref{res-exp-bound}, it immediately follows that $P^{\le} = o(1)$, since 
\[ P^{\le} \le \sum\limits_{i=1}^n \sum\limits_{(x,h) \in \mathcal{G}_T} \pr \Big( \Big| \sum\limits_{t=1}^{T} \mathcal{Z}_{it,T}^{\le}(x,h) \Big| > \frac{C_0}{2} \sqrt{\gamma_{n,T}} \Big). \]
To complete the proof of Proposition \ref{propS1}, it thus remains to verify \eqref{res-exp-bound}. To do so, we split the term $\sum\nolimits_{t=1}^{T} \mathcal{Z}_{it,T}^{\le}(x,h)$ into blocks as follows: 
\[ \sum\limits_{t=1}^{T} \mathcal{Z}_{it,T}^{\le}(x,h) = \sum\limits_{s=1}^{\lceil N_{T} \rceil} B_{2s-1} + \sum\limits_{s=1}^{\lfloor N_{T} \rfloor} B_{2s} \] 
with
\[ B_s = B_{is}(x,h) = \sum\limits_{t = (s-1) L_T + 1}^{\min\{s L_T,T\}} \mathcal{Z}_{it,T}^{\le}(x,h), \]
where $L_T = L_{T,h} = \sqrt{Th / \gamma_{n,T}} \, (nT)^{-1/(\theta-\delta)}$ is the block length and $2 N_{T}$ with $N_{T} = \lceil T/L_T \rceil / 2$ is the number of blocks. Note that under condition \eqref{condition-n}, it holds that $cT^\xi \le L_{T,h} \le CT^{1-\xi}$ for any $h$ with $h_{\min} \le h \le h_{\max}$ and some sufficiently small $\xi > 0$, where $c$, $C$ and $\xi$ are absolute constants that in particular do not depend on $h$. With this notation at hand, we obtain that
\begin{align}
\pr \Big( \Big| \sum\limits_{t=1}^{T} \mathcal{Z}_{it,T}^{\le}(x,h) \Big| > \frac{C_0}{2} \sqrt{\gamma_{n,T}} \Big) 
 & \le \pr \Big( \Big| \sum\limits_{s=1}^{\lceil N_{T} \rceil} B_{2s-1} \Big| > \frac{C_0}{4} \sqrt{\gamma_{n,T}} \Big) \nonumber \\
 & \quad + \pr \Big( \Big| \sum\limits_{s=1}^{\lfloor N_{T} \rfloor} B_{2s} \Big| > \frac{C_0}{4} \sqrt{\gamma_{n,T}} \Big). \label{eq-exp-bound-1}
\end{align}
As the two terms on the right-hand side of \eqref{eq-exp-bound-1} can be treated analogously, we focus attention to the first one. By Bradley's strong approximation theorem (see Theorem 3 in \cite{Bradley1983}), we can construct a sequence of random variables $B_1^*, B_3^*, \ldots$ such that (i) $B_1^*, B_3^*, \ldots$ are independent, (ii) $B_{2s-1}$ and $B_{2s-1}^*$ have the same distribution for each $s$, and (iii) for $0 < \mu \le \| B_{2s-1} \|_{\infty}$, $\pr (|B_{2s-1}^* - B_{2s-1}| > \mu) \le 18 (\| B_{2s-1} \|_{\infty}/\mu)^{1/2} n \, \alpha(L_T)$. With the variables $B_{2s-1}^*$, we can construct the bound
\begin{equation}\label{eq-exp-bound-2}
\pr \Big( \Big| \sum\limits_{s=1}^{\lceil N_{T} \rceil} B_{2s-1} \Big| > \frac{C_0}{4} \sqrt{\gamma_{n,T}} \Big) \le P_1^* + P_2^*, 
\end{equation}
where
\begin{align*}
P_1^* & = \pr \Big( \Big| \sum\limits_{s=1}^{\lceil N_{T} \rceil} B_{2s-1}^* \Big| > \frac{C_0}{8} \sqrt{\gamma_{n,T}} \Big) \\
P_2^* & = \pr \Big( \Big| \sum\limits_{s=1}^{\lceil N_{T} \rceil} \big( B_{2s-1} - B_{2s-1}^* \big) \Big| > \frac{C_0}{8} \sqrt{\gamma_{n,T}} \Big). 
\end{align*}
Using (iii) together with the fact that the mixing coefficients $\alpha(\cdot)$ decay to zero exponentially fast, it is not difficult to see that $P_2^* \le C T^{-r}$, where the constant $r > 0$ can be picked as large as desired. To deal with $P_1^*$, we make use of the following three facts:
\begin{enumerate}[label=(\alph*),leftmargin=0.75cm]
\item For a real-valued random variable $B$ and $\lambda > 0$, Markov's inequality yields that $\pr( \pm B > \delta ) \le \ex \exp(\pm \lambda B) / \exp(\lambda \delta)$.  
\item Since $|B_{2s-1}| \le \{ C L_T (nT)^{1/(\theta-\delta)} \} / \sqrt{Th}$, it holds that $\lambda_{n,T} |B_{2s-1}| \le 1/2$, where we set $\lambda_{n,T} = \sqrt{Th} / \{ 2C L_T (nT)^{1/(\theta-\delta)} \}$. As $\exp(x) \le 1 + x + x^2$ for $|x| \le 1/2$, we get that
\begin{align*}
\ex \Big[ \exp \big( \pm \lambda_{n,T} B_{2s-1} \big) \Big] 
 & \le 1 + \lambda_{n,T}^2 \ex \big[ (B_{2s-1})^2 \big] 
   \le \exp \big( \lambda_{n,T}^2 \ex \big[ (B_{2s-1})^2 \big] \big)
\end{align*}
along with
\[ \ex \Big[ \exp \big( \pm \lambda_{n,T} B_{2s-1}^* \big) \Big] \le \exp \big( \lambda_{n,T}^2 \ex \big[ (B_{2s-1}^*)^2 \big] \big). \]
\item Standard calculations for kernel estimators yield that $\sum\nolimits_{s=1}^{\lceil N_{T} \rceil} \ex \big[ (B_{2s-1}^*)^2 \big] \le C_2$.   
\end{enumerate}
Using (a)--(c), we obtain that
\[ P_1^* \le \pr \Big( \sum\limits_{s=1}^{\lceil N_{T} \rceil} B_{2s-1}^* > \frac{C_0}{8} \sqrt{\gamma_{n,T}} \Big) + \pr \Big( -\sum\limits_{s=1}^{\lceil N_{T} \rceil} B_{2s-1}^* > \frac{C_0}{8} \sqrt{\gamma_{n,T}} \Big), \]
where
\begin{align*} 
 & \pr \Big( \pm \sum\limits_{s=1}^{\lceil N_{T} \rceil} B_{2s-1}^* > \frac{C_0}{8} \sqrt{\gamma_{n,T}} \Big) \\
 & \le \exp \Big( -\frac{C_0}{8} \lambda_{n,T} \sqrt{\gamma_{n,T}} \Big) \, \ex \Big[ \exp \Big( \pm \lambda_{n,T} \sum\limits_{s=1}^{\lceil N_{T} \rceil} B_{2s-1}^* \Big) \Big] \\
 & \le \exp \Big( -\frac{C_0}{8} \lambda_{n,T} \sqrt{\gamma_{n,T}} \Big) \, \prod\limits_{s=1}^{\lceil N_{T} \rceil} \ex \Big[ \exp \big( \pm \lambda_{n,T} B_{2s-1}^* \big) \Big] \\
 & \le \exp \Big( -\frac{C_0}{8} \lambda_{n,T} \sqrt{\gamma_{n,T}} \Big) \, \prod\limits_{s=1}^{\lceil N_{T} \rceil} \exp \Big( \lambda_{n,T}^2 \ex \big[ (B_{2s-1}^*)^2 \big] \Big) \\
 & = \exp \Big( -\frac{C_0}{8} \lambda_{n,T} \sqrt{\gamma_{n,T}} \Big) \exp \Big( \lambda_{n,T}^2 \sum\limits_{s=1}^{\lceil N_{T} \rceil} \ex \big[ (B_{2s-1}^*)^2 \big] \Big) \\
 & \le \exp \Big( -\frac{C_0}{8} \lambda_{n,T} \sqrt{\gamma_{n,T}} + C_2 \lambda_{n,T}^2 \Big). 
\end{align*}
From the definition of $\lambda_{n,T}$, it follows that $\lambda_{n,T} = C_3 \sqrt{\gamma_{n,T}}$ with some absolute constant $C_3 > 0$. Hence, 
\begin{align*}
P_1^* & \le 2 \exp \Big( -\frac{C_0}{8} \lambda_{n,T} \sqrt{\gamma_{n,T}} + C_2 \lambda_{n,T}^2 \Big) \\ 
     & = 2 \exp \Big( - \frac{C_0 C_3}{8} \{ \log n + \log T \} + C_2 C_3^2 \{ \log n + \log T \} \Big) \le C T^{-r}, 
\end{align*}
where the constant $r > 0$ can be made arbitrarily large by picking $C_0$ large enough. To summarize, we have shown that $P_1^* \le C T^{-r}$ and $P_2^* \le C T^{-r}$ with some arbitrarily large $r > 0$. This together with the bounds from \eqref{eq-exp-bound-2} and \eqref{eq-exp-bound-1} yields \eqref{res-exp-bound}, which in turn completes the proof. \qed

\section{Auxiliary results on uniform convergence}\label{secS-uniform2}

We now use Proposition \ref{propS1} from the previous section to derive the uniform convergence rates of some kernel estimators of interest. To start with, we consider the kernel averages
\begin{align}
S_{i,\ell}(x,h) & = \frac{1}{T} \sum\limits_{t=1}^{T} \kernel_h(X_{it} - x) \Big( \frac{X_{it}-x}{h}\Big)^{\ell} \label{def-S} \\
S_{i,\ell}^+(x,h) & = \frac{1}{T} \sum\limits_{t=1}^{T} \kernel_h(X_{it} - x) \Big| \frac{X_{it}-x}{h}\Big|^{\ell} \label{def-S-plus} \\
S_{i,\ell}^{\varepsilon}(x,h) & = \frac{1}{T} \sum\limits_{t=1}^{T} \kernel_h(X_{it} - x) \Big( \frac{X_{it}-x}{h}\Big)^{\ell} \varepsilon_{it} \label{def-S-eps} \\
S_{i,\ell}^m(x,h) & = \frac{1}{T} \sum\limits_{t=1}^{T} \kernel_h(X_{it} - x) \Big( \frac{X_{it}-x}{h}\Big)^{\ell} \{ m_i(X_{it}) - m_i(x) \} \label{def-S-m}
\end{align}
for $0 \le \ell \le 3$. 
\begin{lemmaS}\label{lemmaS1}
Under \ref{C-mixing}, \ref{C-dens} and \ref{C-moments}--\ref{C-ker}, it holds that 
\begin{align}
\max_{1 \le i \le n} \max_{(x,h) \in \mathcal{G}_T} \sqrt{Th} \, \big| S_{i,\ell}(x,h) - \ex [S_{i,\ell}(x,h)] \big| & = O_p \big( \sqrt{\gamma_{n,T}} \big) \label{lemmaS1-statement1} \\  
\max_{1 \le i \le n} \max_{(x,h) \in \mathcal{G}_T} \sqrt{Th} \, \big| S_{i,\ell}^+(x,h) - \ex [S_{i,\ell}^+(x,h)] \big| & = O_p \big( \sqrt{\gamma_{n,T}} \big) \label{lemmaS1-statement2} \\
\max_{1 \le i \le n} \max_{(x,h) \in \mathcal{G}_T} \sqrt{Th} \, \big| S_{i,\ell}^{\varepsilon}(x,h) \big| & = O_p \big( \sqrt{\gamma_{n,T}} \big) \label{lemmaS1-statement3} \\
\max_{1 \le i \le n} \max_{(x,h) \in \mathcal{G}_T} \sqrt{Th} \, \big| S_{i,\ell}^m(x,h) - \ex [S_{i,\ell}^m(x,h)] \big| & = O_p \big( \sqrt{\gamma_{n,T}} \big) \label{lemmaS1-statement4}
\end{align}
with $\gamma_{n,T} = \log n + \log T$. 
\end{lemmaS} 
\noindent \textbf{Proof of Lemma \ref{lemmaS1}.} The terms $S_{i,\ell}(x,h)$ and $S_{i,\ell}^\varepsilon(x,h)$ can be written in the form $T^{-1} \sum\nolimits_{t=1}^{T} \kernel_h(X_{it} - x) \{(X_{it}-x)/h\}^{\ell} Z_{it,T}$ with  $Z_{it,T} = 1$ and $Z_{it,T} = \varepsilon_{it}$, respectively. In addition, $S_{i,\ell}^m(x,h)$ can be expressed as $S_{i,\ell}^m(x,h) = T^{-1} \sum\nolimits_{t=1}^{T} \kernel_h(X_{it} - x) \{(X_{it}-x)/h\}^{\ell} Z_{it,T}^A - m_i(x) \, T^{-1} \sum\nolimits_{t=1}^{T} \kernel_h(X_{it} - x) \{(X_{it}-x)/h\}^{\ell} Z_{it,T}^B$ with $Z_{it,T}^A = m_i(X_{it})$ and $Z_{it,T}^B = 1$. Hence, the statements \eqref{lemmaS1-statement1}, \eqref{lemmaS1-statement3} and \eqref{lemmaS1-statement4} are simple consequences of Proposition \ref{propS1}. Moreover, it is trivial to modify the proof of Proposition \ref{propS1} to apply to the expression $S_{i,\ell}^+(x,h)$ and thus to derive statement \eqref{lemmaS1-statement2}. \qed
\vspace{10pt}

The terms $S_{i,\ell}(x,h)$, $S_{i,\ell}^{\varepsilon}(x,h)$ and $S_{i,\ell}^m(x,h)$ are the building blocks of the local linear kernel averages
\begin{align} 
Q_i(x,h) & = \frac{1}{T} \sum\limits_{t=1}^{T} \wght_{it}(x,h) \label{def-Q} \\
Q_i^{\varepsilon}(x,h) & = \frac{1}{T} \sum\limits_{t=1}^{T} \wght_{it}(x,h) \varepsilon_{it} \label{def-Q-eps} \\
Q_i^m(x,h) & = \frac{1}{T} \sum\limits_{t=1}^{T} \wght_{it}(x,h) \{ m_i(X_{it}) - m_i(x) \}. \label{def-Q-m}
\end{align}
In particular, it holds that 
\begin{align*}
Q_i(x,h) & = S_{i,2}(x,h) S_{i,0}(x,h) - S_{i,1}^2(x,h) \\
Q_i^{\varepsilon}(x,h) & = S_{i,2}(x,h) S_{i,0}^{\varepsilon}(x,h) - S_{i,1}(x,h) S_{i,1}^{\varepsilon}(x,h) \\
Q_i^m(x,h) & = S_{i,2}(x,h) S_{i,0}^m(x,h) - S_{i,1}(x,h) S_{i,1}^m(x,h).
\end{align*}
The uniform convergence rates of $Q_i(x,h)$, $Q_i^\varepsilon(x,h)$ and $Q_i^m(x,h)$ can be easily derived with the help of Lemma \ref{lemmaS1} and some additional straightforward arguments. Defining 
\begin{align*}
Q_i^*(x,h) & = \ex[S_{i,2}(x,h)] \ex[S_{i,0}(x,h)] - \ex[S_{i,1}(x,h)]^2 \\
Q_i^{m,*}(x,h) & = \ex[S_{i,2}(x,h)] \ex[S_{i,0}^m(x,h)] - \ex[S_{i,1}(x,h)] \ex[S_{i,1}^m(x,h)],
\end{align*}
we in particular obtain the following result. 
\begin{lemmaS}\label{lemmaS2}
Under \ref{C-mixing}, \ref{C-dens} and \ref{C-moments}--\ref{C-ker}, it holds that 
\begin{align}
\max_{1 \le i \le n} \max_{(x,h) \in \mathcal{G}_T} \sqrt{Th} \, \big| Q_i(x,h) - Q_i^*(x,h) \big| & = O_p \big( \sqrt{\gamma_{n,T}} \big) \label{lemmaS2-statement1} \\  
\max_{1 \le i \le n} \max_{(x,h) \in \mathcal{G}_T} \sqrt{Th} \, \big| Q_i^{\varepsilon}(x,h) \big| & = O_p \big( \sqrt{\gamma_{n,T}} \big) \label{lemmaS2-statement2} \\
\max_{1 \le i \le n} \max_{(x,h) \in \mathcal{G}_T} \sqrt{Th} \, \big| Q_i^m(x,h) - Q_i^{m,*}(x,h) \big| & = O_p \big( \sqrt{\gamma_{n,T}} \big) \label{lemmaS2-statement3} 
\end{align}
with $\gamma_{n,T} = \log n + \log T$. 
\end{lemmaS} 
\noindent In addition to $Q_i(x,h)$, $Q_i^\varepsilon(x,h)$ and $Q_i^m(x,h)$, we consider the kernel average 
\[ Q_i^{\text{fe}}(x,h) = \frac{1}{T} \sum\limits_{t=1}^{T} \wght_{it}(x,h) \big\{ \overline{\varepsilon}_t^{(i)} + \overline{m}_t^{(i)} \big\}, \]
whose uniform convergence rate is specified by the following lemma. 
\begin{lemmaS}\label{lemmaS3}
Under \ref{C-mixing}, \ref{C-dens} and \ref{C-moments}--\ref{C-ker}, it holds that 
\[ \max_{1 \le i \le n} \max_{(x,h) \in \mathcal{G}_T} \sqrt{Th} \, \big| Q_i^{\textnormal{fe}}(x,h) \big|  = O_p \big( \sqrt{\log n + \log T} \big). \]  
\end{lemmaS} 
\noindent \textbf{Proof of Lemma \ref{lemmaS3}.} Defining 
\[ S_{i,\ell}^{\text{fe}}(x,h) = \frac{1}{T} \sum\limits_{t=1}^{T}  \kernel_h(X_{it} - x) \Big( \frac{X_{it}-x}{h}\Big)^{\ell} Z_{it,T} \] 
with $Z_{it,T} = \overline{\varepsilon}_t^{(i)} + \overline{m}_t^{(i)}$, we can write $Q_i^{\text{fe}}(x,h) = S_{i,2}(x,h) S_{i,0}^{\text{fe}}(x,h) - S_{i,1}(x,h) S_{i,1}^{\text{fe}}(x,h)$. From \ref{C-mixing} and Theorem 5.1(a) in \cite{Bradley2005}, it follows that the collection of random variables $\mathcal{A}_{i,T} = \{ (X_{it}, Z_{it,T}): 1 \le t \le T \}$ is strongly mixing for any $i$ and $T$. In particular, the mixing coefficients $\alpha_{i,T}(\ell)$ of $\mathcal{A}_{i,T}$ are such that $\alpha_{i,T}(\ell) \le n \, \alpha(\ell)$, where the coefficients $\alpha(\ell)$ are defined in \ref{C-mixing} and decay exponentially fast to zero. According to this, the variables $(Z_{it,T},X_{it})$ satisfy condition \ref{P1}. Since the collection of random variables $\{ Z_{it,T}: 1 \le t \le T \}$ is independent from $\{ X_{it}: 1 \le t \le T \}$ for any $i$ under \ref{C-mixing}, it is straightforward to verify that the variables $(Z_{it,T},X_{it})$ fulfill condition \ref{P2} as well. Hence, we can apply Proposition \ref{propS1} to get that 
\[ \max_{1 \le i \le n} \max_{(x,h) \in \mathcal{G}_T} \sqrt{Th} \, \big| S_{i,\ell}^{\textnormal{fe}}(x,h) \big|  = O_p \big( \sqrt{\log n + \log T} \big). \]
With this and Lemma \ref{lemmaS1}, it is straightforward to complete the proof. \qed 
\vspace{10pt}

With the help of the kernel averages defined and analyzed above, the local linear kernel smoothers $\widehat{m}_{i,h}$ can be expressed as 
\[ \widehat{m}_{i,h}(x) - m_i(x) = \frac{Q_i^{\varepsilon}(x,h) + Q_i^m(x,h) - Q_i^{\text{fe}}(x,h)}{Q_i(x,h)} - \big\{ \overline{m}_i + \overline{\varepsilon}_i \big \} + \big\{ \overline{\overline{m}}^{(i)} + \overline{\overline{\varepsilon}}^{(i)} \big\}. \]
We now use this formulation to derive two different uniform expansions of the term $\sqrt{Th} \{ \widehat{m}_{i,h}(x) - m_i(x) \}$, which are required to prove different parts of Theorem \ref{theo-dist}.
\begin{propS}\label{propS2}
Let the conditions of Theorem \ref{theo-dist} be satisfied. Then it holds that 
\[ \sqrt{Th} \big\{ \widehat{m}_{i,h}(x) - m_i(x) \big\} = \sqrt{Th} \frac{Q_i^{m,*}(x,h)}{Q_i^*(x,h)} + R_i^{(a)}(x,h), \]
where the remainder $R_i^{(a)}(x,h)$ has the property that
\[ \max_{1 \le i \le n} \max_{(x,h) \in \mathcal{G}_T} \big| R_i^{(a)}(x,h) \big| = O_p \big( \sqrt{\log n + \log T} \big). \]
\end{propS}
\begin{propS}\label{propS3}
Under the conditions of Theorem \ref{theo-dist}, it holds that 
\[ \sqrt{Th} \big\{ \widehat{m}_{i,h}(x) - m_i(x) \big\} = \sqrt{Th^5} \frac{\kappa(x,h) m_i^{\prime\prime}(x)}{2} + R_i^{(b)}(x,h), \]
where we use the shorthand $\kappa(x,h) = \{ \kappa_2(x,h)^2 - \kappa_1(x,h) \kappa_3(x,h) \} / \{ \kappa_2(x,h) \kappa_0(x,h) - \kappa_1(x,h)^2 \}$ with $\kappa_\ell(x,h) = \int_{-x/h}^{(1-x)/h} u^\ell K(u) du$ and the remainder $R_i^{(b)}(x,h)$ is such that
\[ \max_{1 \le i \le n} \max_{(x,h) \in \mathcal{G}_T} \big| R_i^{(b)}(x,h) \big| = O_p \big( \sqrt{\log n + \log T} + \sqrt{T h_{\max}^7} \big). \]
\end{propS}

\noindent \textbf{Proof of Proposition \ref{propS2}.} Simple algebra yields that 
\[ \sqrt{Th} \big\{ \widehat{m}_{i,h}(x) - m_i(x) \big\} = \sqrt{Th} \frac{Q_i^{m,*}(x,h)}{Q_i^*(x,h)} + R_i^{(a)}(x,h), \]
where $R_i^{(a)}(x,h) = R_{i,1}^{(a)}(x,h) + \ldots + R_{i,6}^{(a)}(x,h)$ with 
\begin{align*}
R_{i,1}^{(a)}(x,h) & = \sqrt{Th} Q_i^{m,*}(x,h) \Big\{ \frac{1}{Q_i(x,h)} - \frac{1}{Q_i^*(x,h)} \Big \} \\
R_{i,2}^{(a)}(x,h) & = \sqrt{Th} \frac{Q_i^m(x,h) - Q_i^{m,*}(x,h)}{Q_i(x,h)} \\
R_{i,3}^{(a)}(x,h) & = \sqrt{Th} \frac{Q_i^{\varepsilon}(x,h)}{Q_i(x,h)} \\ 
R_{i,4}^{(a)}(x,h) & = -\sqrt{Th} \frac{Q_i^{\text{fe}}(x,h)}{Q_i(x,h)} 
\end{align*}
as well as $R_{i,5}^{(a)}(x,h) = -\sqrt{Th} \{ \overline{m}_i + \overline{\varepsilon}_i \}$ and $R_{i,6}^{(a)}(x,h) = \sqrt{Th} \{ \overline{\overline{m}}^{(i)} + \overline{\overline{\varepsilon}}^{(i)} \}$. To complete the proof, we show that  
\begin{equation}\label{propS2-remainders}
\max_{1 \le i \le n} \max_{(x,h) \in \mathcal{G}_T} \big| R_{i,\ell}^{(a)}(x,h) \big| = O_p \big( \sqrt{\log n + \log T} \big) 
\end{equation}
for $1 \le \ell \le 6$: By standard bias calculations, we obtain that 
\begin{align}
\max_{1 \le i \le n} \max_{(x,h) \in \mathcal{G}_T} \big|Q_i^{m,*}(x,h)\big| & = O(h_{\max}) \label{bias-Qim} \\
\max_{1 \le i \le n} \max_{(x,h) \in \mathcal{G}_T} \big| Q_i^*(x,h) - \big\{ \kappa_2(x,h) \kappa_0(x,h) - \kappa_1(x,h)^2 \big\} f_i^2(x) \big| & = O(h_{\max}), \label{bias-Qi}
\end{align}
where under our assumptions, the term $Q_i^{**}(x,h) = \{ \kappa_2(x,h) \kappa_0(x,h) - \kappa_1(x,h)^2 \} f_i^2(x)$ is bounded away from zero and infinity uniformly over $i$ and $(x,h)$, that is, $0 < c \le Q_i^{**}(x,h) \le C < \infty$ with some constants $c$ and $C$ that are independent of $i$ and $(x,h)$. With the help of these observations and Lemmas \ref{lemmaS2} and \ref{lemmaS3}, it is straightforward to derive \eqref{propS2-remainders} for $1 \le \ell \le 4$. Next, note that $\max_{1 \le i \le n} |\overline{\overline{m}}^{(i)}| \le \max_{1 \le i \le n} |\overline{m}_i|$ and $\max_{1 \le i \le n} |\overline{\overline{\varepsilon}}^{(i)}| \le \max_{1 \le i \le n} |\overline{\varepsilon}_i|$. Arguments similar to but simpler than those for Proposition \ref{propS1} yield that $\max_{1 \le i \le n} |\overline{m}_i| = O_p(\sqrt{\{\log n + \log T\}/T})$ and $\max_{1 \le i \le n} |\overline{\varepsilon}_i| = O_p(\sqrt{\{\log n + \log T\}/T})$. From this, \eqref{propS2-remainders} immediately follows for $\ell = 5$ and $\ell = 6$. \qed
\vspace{10pt}

\noindent \textbf{Proof of Proposition \ref{propS3}.} Straightforward calculations yield that 
\[ \sqrt{Th} \big\{ \widehat{m}_{i,h}(x) - m_i(x) \big\} = \sqrt{Th^5} \frac{\kappa(x,h) m_i^{\prime\prime}(x)}{2} + R_i^{(b)}(x,h), \] 
where $R_i^{(b)}(x,h) = R_{i,1}^{(b)}(x,h) + \ldots + R_{i,5}^{(b)}(x,h)$ with
\[ R_{i,1}^{(b)}(x,h) = \sqrt{Th} \Big\{ \frac{Q_i^m(x,h)}{Q_i(x,h)} - h^2 \frac{\kappa(x,h) m_i^{\prime\prime}(x)}{2} \Big\} \]
and $R_{i,\ell}^{(b)}(x,h) = R_{i,\ell+1}^{(a)}(x,h)$ for $2 \le \ell \le 5$. In order to prove Proposition \ref{propS3}, it suffices to show that 
\begin{align} 
\max_{1 \le i \le n} \max_{(x,h) \in \mathcal{G}_T} \big| R_{i,1}^{(b)}(x,h) \big| & = O_p \big( \sqrt{T h_{\max}^7} \big) + o_p \big( \sqrt{\log n + \log T} \big) \label{propS3-claim1} \\
\max_{1 \le i \le n} \max_{(x,h) \in \mathcal{G}_T} \big| R_{i,\ell}^{(b)}(x,h) \big| & = O_p \big( \sqrt{\log n + \log T} \big) \label{propS3-claim2} 
\end{align}
for $2 \le \ell \le 5$. \eqref{propS3-claim2} has already been verified in the proof of Proposition \ref{propS2}. To prove \eqref{propS3-claim1}, we make use of the following two facts:

\begin{enumerate}[label=(\alph*),leftmargin=0.75cm]

\item From Lemma \ref{lemmaS2} and \eqref{bias-Qi}, it follows that
\begin{equation} \label{propS3-claim1-fact1}
\max_{1 \le i \le n} \max_{(x,h) \in \mathcal{G}_T} \sqrt{Th} \big| Q_i(x,h) - Q_i^{**}(x,h) \big| = O_p \big( \sqrt{\log n + \log T} + \sqrt{T h_{\max}^3} \big)
\end{equation}
with $Q_i^{**}(x,h) = \{ \kappa_2(x,h) \kappa_0(x,h) - \kappa_1(x,h)^2 \} f_i^2(x)$. As already noted in the proof of Proposition \ref{propS2}, the term $Q_i^{**}(x,h)$ is bounded away from zero and infinity uniformly over $i$ and $(x,h)$.

\item A second-order Taylor expansion of $m_i$ yields that 
\begin{equation}\label{propS3-claim1-fact2a} 
\sqrt{Th} Q_i^m(x,h) = \sqrt{Th} Q_i^{m,**}(x,h) + R_i^m(x,h), 
\end{equation}
where 
\[ Q_i^{m,**}(x,h) = h^2 \frac{m_i^{\prime\prime}(x) f_i^2(x)}{2} \big[ \kappa_2(x,h)^2 - \kappa_1(x,h) \kappa_3(x,h) \big]. \] 
The remainder term $R_i^m(x,h)$ has the form $R_i^m(x,h) = R_{i,1}^m(x,h) + R_{i,2}^m(x,h)$, where
\begin{align*}
R_{i,1}^m(x,h) & = \sqrt{Th^5} \frac{m_i^{\prime\prime}(x)}{2} \Big\{ \big[ S_{i,2}(x,h)^2 - S_{i,1}(x,h) S_{i,3}(x,h) \big] \\ & \qquad - \big[ \kappa_2(x,h)^2 - \kappa_1(x,h) \kappa_3(x,h) \big] f_i^2(x) \Big \} \\
R_{i,2}^m(x,h) & = \frac{\sqrt{Th^5}}{2 T} \sum\limits_{t=1}^{T} \kernel_h(X_{it}-x) \Big[ S_{i,2}(x,h) - \Big(\frac{X_{it}-x}{h}\Big) S_{i,1}(x,h)\Big] \\ & \qquad \times \big\{ m_i^{\prime\prime}(\xi_{it}) - m_i^{\prime\prime}(x) \big\} \Big(\frac{X_{it}-x}{h}\Big)^2
\end{align*}
with $\xi_{it}$ denoting an intermediate point between $X_{it}$ and $x$. By Lemma \ref{lemmaS1} and standard bias calculations, we obtain that
\begin{equation}\label{propS3-claim1-fact2b}
\max_{1 \le i \le n} \max_{(x,h) \in \mathcal{G}_T} \big| R_{i,1}^m(x,h) \big| = O_p\big( h_{\max}^2 \sqrt{\log n + \log T} + \sqrt{Th_{\max}^7} \big). 
\end{equation}
As $m^{\prime\prime}_i$ is Lipschitz continuous by \ref{C-m}, we further get that $| R_{i,2}^m(x,h) | \le C \sqrt{Th^7} \linebreak \{ S_{i,2}(x,h)^2 + S_{i,1}^+(x,h) S_{i,3}^+(x,h) \}$. Applying Lemma \ref{lemmaS1} together with standard bias calculations to this upper bound, we can infer that 
\begin{equation}\label{propS3-claim1-fact2c}
\max_{1 \le i \le n} \max_{(x,h) \in \mathcal{G}_T} \big| R_{i,2}^m(x,h) \big| = O_p\big( h_{\max}^3 \sqrt{\log n + \log T} + \sqrt{Th_{\max}^7} \big).  
\end{equation}
Finally, by combining \eqref{propS3-claim1-fact2b} and \eqref{propS3-claim1-fact2c}, the remainder term $R_i^m(x,h)$ is seen to have the property that 
\begin{equation}\label{propS3-claim1-fact2d}
\max_{1 \le i \le n} \max_{(x,h) \in \mathcal{G}_T} \big| R_i^m(x,h) \big| = O_p\big( h_{\max}^2 \sqrt{\log n + \log T} + \sqrt{Th_{\max}^7} \big). 
\end{equation}

\end{enumerate}

\noindent We now proceed as follows: Simple algebra yields that 
\begin{align*}
 & \sqrt{Th} \Big( \frac{Q_i^m(x,h)}{Q_i(x,h)} - \frac{Q_i^{m,**}(x,h)}{Q_i^{**}(x,h)} \Big) \\*
 & \qquad = \frac{R_i^m(x,h)}{Q_i(x,h)} + \sqrt{Th} Q_i^{m,**}(x,h) \Big\{ \frac{1}{Q_i(x,h)} - \frac{1}{Q_i^{**}(x,h)} \Big\}. 
\end{align*} 
Since $Q_i^{m,**}(x,h)/Q_i^{**}(x,h) = h^2 \kappa(x,h) m_i^{\prime\prime}(x) / 2$, this implies that
\[ R_{i,1}^{(b)}(x,h) = \frac{R_i^m(x,h)}{Q_i(x,h)} + \sqrt{Th} Q_i^{m,**}(x,h) \Big\{ \frac{1}{Q_i(x,h)} - \frac{1}{Q_i^{**}(x,h)} \Big\}. \]
Using this representation of $R_{i,1}^{(b)}(x,h)$ together with \eqref{propS3-claim1-fact1}, \eqref{propS3-claim1-fact2d} and the fact that $Q_i^{**}(x,h)$ is bounded away from zero and infinity uniformly over $i$ and $(x,h)$, it is straightforward to verify \eqref{propS3-claim1}. \qed
\vspace{10pt}

The final result of this section is concerned with the normalization term 
\begin{equation}\label{def-nu}
\widehat{\nu}_{ij}(x,h) = \left\{ \frac{\widehat{\sigma}_{i,h}^2}{\widehat{f}_{i,h}(x)} + \frac{\widehat{\sigma}_{j,h}^2}{\widehat{f}_{j,h}(x)} \right\} s(x,h), 
\end{equation}
where $s(x,h) = \{ \int_{-x/h}^{(1-x)/h} \kernel^2(u) [ \kappa_2(x,h) - \kappa_1(x,h) u ]^2 du \} / \{ \kappa_0(x,h) \kappa_2(x,h) - \kappa_1(x,h)^2 \}^2$ with $\kappa_\ell(x,h) = \int_{-x/h}^{(1-x)/h} u^\ell \kernel(u) du$ for $0 \le \ell \le 2$, $\widehat{f}_{i,h}(x) = \{\kappa_0(x,h) T\}^{-1} \sum\nolimits_{t=1}^{T} \kernel_h(X_{it} - x)$ and $\widehat{\sigma}_{i,h}^2 = T^{-1} \sum\nolimits_{t=1}^{T} \{ \widehat{Y}_{it}^* - \widehat{m}_{i,h}(X_{it}) \}^2$.

\begin{propS}\label{propS4}
Let the conditions of Theorem \ref{theo-dist} be satisfied. Then there exist absolute constants $0 < c_\nu \le C_\nu < \infty$ such that 
\begin{align*}
\min_{1 \le i \le j \le n} \min_{(x,h) \in \mathcal{G}_T} \sqrt{\widehat{\nu}_{ij}(x,h)} & \ge c_\nu + o_p(1) \\
\max_{1 \le i \le j \le n} \max_{(x,h) \in \mathcal{G}_T} \sqrt{\widehat{\nu}_{ij}(x,h)} & \le C_\nu + o_p(1).
\end{align*}
\end{propS}

\noindent \textbf{Proof of Proposition \ref{propS4}.} The proposition is a straightforward consequence of the following three observations: 
\begin{enumerate}[label=(\alph*),leftmargin=0.75cm]

\item \label{fact1-propS4} Under our conditions, the term $s(x,h)$ is bounded away from zero and infinity uniformly over $(x,h)$, that is, $0 < c_s \le s(x,h) \le C_s < \infty$ for some absolute constants $c_s$ and $C_s$. 

\item \label{fact2-propS4} It holds that 
\[ \max_{1 \le i \le n} \max_{(x,h) \in \mathcal{G}_T} \big| \widehat{f}_{i,h}(x) - f_i(x) \big| = O_p \Big( \sqrt{\frac{\log n + \log T}{T h_{\min}}} + h_{\max} \Big), \]
where the densities $f_i$ are uniformly bounded away from zero and infinity by \ref{C-dens}.

\item \label{fact3-propS4} It holds that 
\begin{equation*}
\widehat{\sigma}_{i,h}^2 = \sigma_i^2 + b_i^\sigma + R_{i,h}^\sigma \quad \text{with} \quad \max_{1 \le i \le n} \max_{\{ h: (x,h) \in \mathcal{G}_T \}} |R_{i,h}^\sigma| = o_p(1), 
\end{equation*}
where $b_i^{\sigma} = \ex[ (\overline{m}_t^{(i)} + \overline{\varepsilon}_t^{(i)})^2 ]$ and the error variances $\sigma_i^2$ are uniformly bounded away from zero and infinity according to \ref{C-sigma}. Note that $0 \le b_i^\sigma \le C_b < \infty$ for some sufficiently large absolute constant $C_b$ and that $\max_{1 \le i \le n} b_i^{\sigma} = o(1)$ in the case that $n$ tends to infinity as $T \rightarrow \infty$.   

\end{enumerate}
\noindent Observation \ref{fact1-propS4} can be seen by straightforward arguments and \ref{fact2-propS4} follows from Lemma \ref{lemmaS1} together with standard bias calculations. In order to prove \ref{fact3-propS4}, we write $\widehat{\sigma}_{i,h}^2 = \sigma_i^2 + b_i^\sigma + R_{i,h}^\sigma$ with $R_{i,h}^\sigma = R_{i,h,1}^\sigma + \ldots + R_{i,h,5}^\sigma$, where
\begin{align*}
R_{i,h,1}^\sigma & = \frac{1}{T} \sum\limits_{t=1}^{T} \big\{ \varepsilon_{it}^2 -  \ex \big[ \varepsilon_{it}^2 \big] \big \} \\
R_{i,h,2}^\sigma & = \frac{1}{T} \sum\limits_{t=1}^T \big\{ (\overline{m}_t^{(i)} + \overline{\varepsilon}_t^{(i)})^2  - \ex \big[ (\overline{m}_t^{(i)} + \overline{\varepsilon}_t^{(i)})^2 \big] \big\} \\ 
R_{i,h,3}^\sigma & = \frac{1}{T} \sum\limits_{t=1}^T \big\{ \widehat{\Delta}_{i,h}(X_{it}) - (\overline{m}_i + \overline{\varepsilon}_i) + (\overline{\overline{m}}^{(i)} + \overline{\overline{\varepsilon}}^{(i)}) \big\}^2 \\
R_{i,h,4}^\sigma & = -\frac{2}{T} \sum\limits_{t=1}^T \big\{ \overline{m}_t^{(i)} + \overline{\varepsilon}_t^{(i)} \big\} \big\{ \widehat{\Delta}_{i,h}(X_{it}) - (\overline{m}_i + \overline{\varepsilon}_i) + (\overline{\overline{m}}^{(i)} + \overline{\overline{\varepsilon}}^{(i)}) \big\} \\
R_{i,h,5}^\sigma & = \frac{2}{T} \sum\limits_{t=1}^T \varepsilon_{it} \big\{ \widehat{\Delta}_{i,h}(X_{it}) - (\overline{m}_i + \overline{\varepsilon}_i) - (\overline{m}_t^{(i)} + \overline{\varepsilon}_t^{(i)}) + (\overline{\overline{m}}^{(i)} + \overline{\overline{\varepsilon}}^{(i)}) \big\}
\end{align*}
with the shorthand $\widehat{\Delta}_{i,h}(X_{it}) = m_i(X_{it}) - \widehat{m}_{i,h}(X_{it})$. A simplified version of Proposition \ref{propS1} yields that 
\begin{align}
\max_{1 \le i \le n} \Big| \frac{1}{T} \sum\limits_{t=1}^{T} \big\{ \varepsilon_{it}^2 - \ex \big[ \varepsilon_{it}^2 \big] \big\} \Big| & = O_p \Big( \sqrt{\frac{\log n + \log T}{T}} \Big). \label{propS4-fact3-aux2}
\end{align}
By \ref{C-mixing} and Theorem 5.1(a) in \cite{Bradley2005}, the collection of random variables $\mathcal{A}_{i,T} = \{ (\varepsilon_{it}, \overline{\varepsilon}_t^{(i)}, \overline{m}_t^{(i)}): 1 \le t \le T \}$ is strongly mixing for any $i$ and $T$, where the mixing coefficients $\alpha_{i,T}(\ell)$ of $\mathcal{A}_{i,T}$ are such that $\alpha_{i,T}(\ell) \le n \, \alpha(\ell)$ with $\alpha(\ell)$ decaying to zero exponentially fast. For this reason, we can once again apply a simplified version of  Proposition \ref{propS1} to obtain that
\begin{align} 
\max_{1 \le i \le n} \Big| \frac{1}{T} \sum\limits_{t=1}^T \varepsilon_{it} \big(\overline{m}_t^{(i)} + \overline{\varepsilon}_t^{(i)}\big) \Big| & = O_p \Big( \sqrt{\frac{\log n + \log T}{T}} \Big) \label{propS4-fact3-aux3a} \\
\max_{1 \le i \le n} \Big| \frac{1}{T} \sum\limits_{t=1}^T \Big\{ (\overline{m}_t^{(i)} + \overline{\varepsilon}_t^{(i)})^2  - \ex \big[ (\overline{m}_t^{(i)} + \overline{\varepsilon}_t^{(i)})^2 \big] \Big\} \Big| & = O_p \Big( \sqrt{\frac{\log n + \log T}{T}} \Big). \label{propS4-fact3-aux3}
\end{align}
Moreover, slightly modifying the proof of Proposition \ref{propS3}, we can infer that
\begin{equation}\label{propS4-fact3-aux4} 
\max_{1 \le i \le n} \max_{(x,h) \in \mathcal{G}_T} \big| \widehat{\Delta}_{i,h}(x) \big| = O_p \Big( \sqrt{\frac{\log n + \log T}{T h_{\min}}} + h_{\max}^2 \Big). 
\end{equation}
Finally, as already seen in the proof of Proposition \ref{propS2}, 
\begin{align}
\max_{1 \le i \le n} \big|\overline{m}_i + \overline{\varepsilon}_i\big| & = O_p \Big( \sqrt{\frac{\log n + \log T}{T}} \Big) \label{propS4-fact3-aux5} \\
\max_{1 \le i \le n} \big|\overline{\overline{m}}^{(i)} + \overline{\overline{\varepsilon}}^{(i)}\big| & = O_p \Big( \sqrt{\frac{\log n + \log T}{T}} \Big). \label{propS4-fact3-aux6} 
\end{align} 
With the help of \eqref{propS4-fact3-aux2}--\eqref{propS4-fact3-aux6}, it is not difficult to infer that 
\begin{equation}
\max_{1 \le i \le n} \max_{\{ h: (x,h) \in \mathcal{G}_T \}} |R_{i,h,\ell}^\sigma| = o_p(1) 
\end{equation} 
for $1 \le \ell \le 5$, which implies \ref{fact3-propS4}. \qed

\section{Proof of Theorem \ref{theo-dist}}\label{secS-theo-dist}

\textbf{Proof of (\ref{theo-dist-stat1}).} From Proposition \ref{propS2}, it follows that
\begin{align*}
\sqrt{Th} & \big\{ \widehat{m}_{i,h}(x) - \widehat{m}_{j,h}(x) \big\} \\
 & = \sqrt{Th} \big\{ m_i(x) - m_j(x) \big\} \\
 & \qquad + \sqrt{Th} \Big \{ \frac{Q_i^{m,*}(x,h)}{Q_i^*(x,h)} - \frac{Q_j^{m,*}(x,h)}{Q_j^*(x,h)} \Big\} + R_{ij}(x,h), 
\end{align*}
where $\max_{1 \le i \le j \le n} \max_{(x,h) \in \mathcal{G}_T} |R_{ij}(x,h)| = O_p(\sqrt{\log n + \log T})$. Since $Q_i^{m,*}(x,h) = Q_j^{m,*}(x,h)$ and $Q_i^*(x,h) = Q_j^*(x,h)$ for any two time series $i$ and $j$ in the same group $G_k$ under our conditions, this implies that 
\begin{equation}
\max_{1 \le k \le K_0} \max_{i,j \in G_k} \max_{(x,h) \in \mathcal{G}_T} \sqrt{Th} \big| \widehat{m}_{i,h}(x) - \widehat{m}_{j,h}(x) \big| = O_p \big( \sqrt{\log n + \log T} \big). 
\end{equation}
Moreover, by Proposition \ref{propS4}, 
\[ \min_{1 \le i \le j \le n} \min_{(x,h) \in \mathcal{G}_T} \sqrt{\widehat{\nu}_{ij}(x,h)} \ge c_{\nu} + o_p(1), \] 
where $c_{\nu} > 0$ is a sufficiently small absolute constant. As a result, we arrive at
\begin{align*} 
\max_{1 \le k \le K_0} \max_{i,j \in G_k} \widehat{d}_{ij}
 & \le \max_{1 \le k \le K_0} \max_{i,j \in G_k} \Big\{ \max_{(x,h) \in \mathcal{G}_T} | \widehat{\psi}_{ij}(x,h) | \Big\} \\
 & \le \frac{\displaystyle{ \max_{1 \le k \le K_0} \max_{i,j \in G_k} \max_{(x,h) \in \mathcal{G}_T} \sqrt{Th} \, \big| \widehat{m}_{i,h}(x) - \widehat{m}_{j,h}(x) \big| }}{\displaystyle{ \min_{1 \le i \le j \le n} \min_{(x,h) \in \mathcal{G}_T} \sqrt{\widehat{\nu}_{ij}(x,h)}}} \\
 & = O_p \big( \sqrt{\log n + \log T} \big),
\end{align*}
which completes the proof. \qed
\vspace{10pt}

\noindent \textbf{Proof of (\ref{theo-dist-stat2}).} By Proposition \ref{propS3}, it holds that
\begin{align*}
\sqrt{Th} & \big\{ \widehat{m}_{i,h}(x) - \widehat{m}_{j,h}(x) \big\} \\
 & = \sqrt{Th} \big\{ m_i(x) - m_j(x) \big\} \\
 & \qquad + \sqrt{Th^5} \frac{\kappa(x,h)}{2} \big\{ m_i^{\prime\prime}(x) - m_j^{\prime\prime}(x) \big\} + R_{ij}(x,h),
\end{align*}
where $\max_{1 \le i \le j \le n} \max_{(x,h) \in \mathcal{G}_T} |R_{ij}(x,h)| = O_p(\sqrt{\log n + \log T} + \sqrt{Th_{\max}^7})$. With the help of this expansion, we can infer that    
\begin{align*} 
 & \min_{1 \le k < k^\prime \le K_0} \min_{\substack{i \in G_k, \\ j \in G_{k^\prime}}} \max_{(x,h) \in \mathcal{G}_T} \sqrt{Th} \, \big| \widehat{m}_{i,h}(x) - \widehat{m}_{j,h}(x) \big| \\*
 & \qquad \ge \min_{1 \le k < k^\prime \le K_0} \min_{\substack{i \in G_k, \\ j \in G_{k^\prime}}} \max_{(x,h) \in \mathcal{G}_T} \sqrt{Th} \, \big| m_i(x) - m_j(x) \big| \\
 & \qquad \qquad - \max_{1 \le i \le j \le n} \max_{(x,h) \in \mathcal{G}_T} \sqrt{Th^5} \, \frac{|\kappa(x,h)|}{2} \, \big| m_i^{\prime\prime}(x) - m_j^{\prime\prime}(x) \big| \\
 & \qquad \qquad - \max_{1 \le i \le j \le n} \max_{(x,h) \in \mathcal{G}_T} \big| R_{ij}(x,h) \big| \\
 & \qquad = \min_{1 \le k < k^\prime \le K_0} \min_{\substack{i \in G_k, \\ j \in G_{k^\prime}}} \max_{(x,h) \in \mathcal{G}_T} \sqrt{Th} \, \big| m_i(x) - m_j(x) \big| \\
 & \qquad \qquad + O_p \big( \sqrt{T h_{\max}^5} + \sqrt{\log n + \log T} \big) \\[0.2cm]
 & \qquad \ge c \sqrt{T h_{\max}} + o_p \big( \sqrt{T h_{\max}} \big),
\end{align*}  
where $c > 0$ is a sufficiently small absolute constant. Moreover, by Proposition \ref{propS4}, 
\[ \max_{1 \le i \le j \le n} \max_{(x,h) \in \mathcal{G}_T} \sqrt{\widehat{\nu}_{ij}(x,h)} \le C_\nu + o_p(1) \] 
with $C_\nu > 0$ being an absolute constant that is chosen sufficiently large. As a consequence, we get that 
\begin{align} 
 & \min_{1 \le k < k^\prime \le K_0} \min_{\substack{i \in G_k, \\ j \in G_{k^\prime}}} \Big\{ \max_{(x,h) \in \mathcal{G}_T} | \widehat{\psi}_{ij}(x,h) | \Big\} \nonumber \\*[0.1cm]
 & \quad \ge \frac{\displaystyle{ \min_{1 \le k < k^\prime \le K_0} \, \min_{i \in G_k, \, j \in G_{k^\prime}} \max_{(x,h) \in \mathcal{G}_T} \sqrt{T h} \, \big| \widehat{m}_{i,h}(x) - \widehat{m}_{j,h}(x) \big|}}{\displaystyle{ \max_{1 \le i \le j \le n} \max_{(x,h) \in \mathcal{G}_T} \sqrt{\widehat{\nu}_{ij}(x,h)}} } \nonumber \\[0.1cm]
 & \quad \ge c_0 \sqrt{T h_{\max}} + o_p \big( \sqrt{T h_{\max}} \big) \label{theo-dist-eq2-step}
\end{align}
with some sufficiently small absolute constant $c_0$. Since $\lambda(2 h_{\min}) = O(\sqrt{\log T})$ by the conditions on the bandwidth $h_{\min}$ in \ref{C-h}, we finally obtain that 
\begin{align*}
\min_{1 \le k < k^\prime \le K_0} \min_{\substack{i \in G_k, \\ j \in G_{k^\prime}}} \widehat{d}_{ij} 
 & \ge \min_{1 \le k < k^\prime \le K_0} \min_{\substack{i \in G_k, \\ j \in G_{k^\prime}}} \Big\{ \max_{(x,h) \in \mathcal{G}_T} | \widehat{\psi}_{ij}(x,h) | \Big\} - \lambda(2 h_{\min}) \\
 & = \min_{1 \le k < k^\prime \le K_0} \min_{\substack{i \in G_k, \\ j \in G_{k^\prime}}} \Big\{ \max_{(x,h) \in \mathcal{G}_T} | \widehat{\psi}_{ij}(x,h) | \Big\} + O(\sqrt{\log T}) \\
 & \ge c_0 \sqrt{T h_{\max}} + o_p \big( \sqrt{T h_{\max}} \big), 
\end{align*}
the last line following from \eqref{theo-dist-eq2-step}. \qed